\documentclass[12pt]{article}
\usepackage{amssymb,epsfig}

\usepackage[french]{babel}
\usepackage[T1]{fontenc}
\usepackage{color}

\newtheorem{defi}{D\'efinition}[section]
\newtheorem{prop}[defi]{Proposition}
\newtheorem{theo}[defi]{Th\'eor\`eme}
\newtheorem{conj}[defi]{Conjecture}
\newtheorem{lemm}[defi]{Lemme}
\newtheorem{coro}[defi]{Corollaire}
\newtheorem{rema}[defi]{Remarque}
\newtheorem{exem}[defi]{Exemple}
\newtheorem{exer}{Exercice}

\newcommand{\bdefi}{\begin{defi}}
\newcommand{\edefi}{\end{defi}}
\newcommand{\bprop}{\begin{prop}}
\newcommand{\eprop}{\end{prop}}
\newcommand{\btheo}{\begin{theo}}
\newcommand{\etheo}{\end{theo}}
\newcommand{\blemm}{\begin{lemm}}
\newcommand{\brema}{\begin{rema}}
\newcommand{\erema}{\end{rema}}
\newcommand{\bexer}{\begin{exer}}
\newcommand{\eexer}{\end{exer}}
\newcommand{\bconj}{\begin{conj}}
\newcommand{\econj}{\end{conj}}
\newcommand{\elemm}{\end{lemm}}
\newcommand{\bcoro}{\begin{coro}}
\newcommand{\ecoro}{\end{coro}}
\newcommand{\bexem}{\begin{exem}}
\newcommand{\eexem}{\end{exem}}
\newcommand{\dem}{\noindent{\bf Preuve. }}
\newcommand{\rem}{\noindent{\bf Remarque. }}

\newcommand{\A}{{\cal A}}

\renewcommand{\P}{{\cal P}}
\newcommand{\G}{{\cal G}}

\renewcommand{\S}{{\cal S}}

\usepackage{amssymb}

\newcommand{\maths}[1]{{\mathbb #1}}  

\newcommand{\RR}{\maths{R}}
\newcommand{\NN}{\maths{N}}

\newcommand{\QQ}{\maths{Q}}

\newcommand{\FF}{\maths{F}}

\newcommand{\ZZ}{\maths{Z}}

\newcommand{\TT}{\maths{T}}

\newcommand{\ra}{\rightarrow}

\newcommand{\ov}[1]{{\overline{#1}}} 
\newcommand{\wt}[1]{{\widetilde{#1}}}
\newcommand{\wh}[1]{{\widehat{#1}}}

\newcommand{\ga}{\gamma}
\newcommand{\Ga}{\Gamma}

\newcommand{\cqfd}{\hfill$\Box$}
\newcommand{\bac}{{\backslash\!\backslash}}

\newcounter{fig}

\usepackage{epsfig}

\usepackage{color}

\title{Sur le codage du flot g\'eod\'esique dans un arbre}

\author{Anne Broise-Alamichel \and Fr\'ed\'eric Paulin}
\date{}
\begin{document}
\maketitle

\noindent {\small
\begin{tabular}{l}
Laboratoire de Math\'ematique UMR 8628 CNRS\\
Equipe de Topologie et Dynamique (B\^at. 425)\\
Universit\'e Paris-Sud \\
91405 ORSAY Cedex, FRANCE.\\
{\it e-mail: Anne.Broise@math.u-psud.fr}
\end{tabular}
\\
   \mbox{}
\\
   \mbox{}
\\
\begin{tabular}{l}
D\'epartement de Math\'ematique et Applications, UMR 8553 CNRS\\
Ecole Normale Sup\'erieure\\
45 rue d'Ulm\\
75230 PARIS Cedex 05, FRANCE\\
{\it e-mail: Frederic.Paulin@ens.fr}
\end{tabular}
}

\begin{center}
{\bf R\'esum\'e}
\end{center}
\'Etant donn\'e un arbre $T$ et un groupe $\Ga$ d'automorphismes de
$T$, nous \'etudions les propri\'et\'es markoviennes du flot
g\'eod\'esique sur le quotient de l'espace des g\'eod\'esiques de $T$
par $\Ga$. Par exemple, quand $T$ est l'arbre de Bruhat-Tits d'un
groupe alg\'ebrique lin\'eaire connexe semi-simple $\underline{G}$ de
rang 1 au dessus d'un corps local non archim\'edien $\wh K$ et si
$\Ga$ est r\'eseau (\'eventuellement non uniforme) dans
$\underline{G}(\wh K)$, nous montrons que l'action des puissances
paires de la transformation g\'eod\'esique est Bernoulli d'entropie
finie. Sous des hypoth\`eses g\'en\'erales b\'enignes, nous montrons
que si le flot g\'eod\'esique est m\'elangeant pour une mesure de
probabilit\'e de Patterson-Sullivan-Bowen-Margulis, alors il est
l\^achement Bernoulli.

\begin{center}
{\bf Abstract}
\end{center}
  Given a tree $T$ and a group $\Ga$ of automorphisms of $T$, we study
  the markovian properties of the geodesic flow on the quotient by
  $\Ga$ of the space of geodesics of $T$. For instance, when $T$ is
  the Bruhat-Tits tree of a semi-simple connected algebraic group
  $\underline{G}$ of rank one over a non archimedian local field $\wh
  K$, and $\Ga$ is a (possibly non uniform) lattice in
  $\underline{G}(\wh K)$, we prove that the type preserving geodesic
  flow is Bernoulli with finite entropy. Under some mild assumptions, we prove that if the
  quotient geodesic flow is mixing for a probability
  Patterson-Sullivan-Bowen-Margulis measure, then it is loosely
  Bernoulli.  \footnote{ {\bf AMS codes:} 20 G 25, 20 E 08, 37 A 25.
    {\bf Keywords:} group actions on trees, Bruhat-Tits trees,
    geodesic flow, coding, Bernoulli shifts.  }

\newpage

\section{Introduction}
\label{sec:Intro}

Soit $T$ un arbre localement fini, $\Ga$ un sous-groupe discret
d'automorphismes de $T$, $\G T$ l'espace des g\'eod\'esiques de $T$
(i.e.~des isom\'etries $\ell:\RR\ra T$ d'origine $\ell(0)$ un sommet
de $T$), et $\wt\varphi:\G T\ra \G T$ la transformation
g\'eod\'esique sur $\G T$, d\'efinie par $\ell\mapsto \{t\mapsto
\ell(t+1)\}$. 

Le but de cet article est d'\'etudier la dynamique symbolique de la
transformation g\'eod\'esique quotient $\varphi:\Ga\backslash\G T\ra
\Ga\backslash \G T$ de $\wt\varphi$, pour obtenir des propri\'et\'es
ergodiques plus fines que celles obtenues dans \cite{BM,Rob}.  Il ne
s'agit pas de se restreindre au cas o\`u $\Ga$ est un r\'eseau
uniforme, qui est bien connu et bien plus \'el\'ementaire (voir par
exemple les r\'ef\'erences \cite{Coo,CP}, qui s'int\'eressent au cas
plus g\'en\'eral des groupes hyperboliques). Nous nous int\'eressons
au contraire au cas des r\'eseaux non uniformes (voir le livre
\cite{BL} pour avoir une id\'ee de la richesse des exemples)~; en
g\'en\'eral on ne peut pas se d\'ebarasser de la torsion par passage
\`a un sous-groupe d'indice fini, ceci est un probl\`eme crucial en ce
qui concerne le codage.  En supposant que la mesure de
Patterson-Sullivan ne charge pas les ensembles de points fixes
d'\'el\'ements elliptiques non triviaux, un premier r\'esultat de
codage (voir paragraphe \ref{sec:codageneuf}) est le suivant (voir
paragraphe \ref{sec:rappels} pour des rappels de d\'efinitions),

\btheo \label{theo:introneuf} %
Soit ${\wt \mu}_{\mbox{\tiny BM}}$ une mesure de
(Patterson-Sullivan)-Bowen-Margulis pour $\Ga$ sur $\G T$. Supposons
que le syst\`eme dynamique mesur\'e quotient de $(\G T,\wt\varphi,{\wt
  \mu}_{\mbox{\tiny BM}})$ par $\Ga$ soit de probabilit\'e et
m\'elangeant. Alors il est l\^achement Bernoulli.  \etheo

Voir \cite{BM,Rob} (ou le paragraphe \ref{sec:finibowmarg}) pour de
grandes classes d'exemples o\`u les conditions de finitude de la
mesure et de m\'elange sont v\'erifi\'ees.  Dans le cadre
alg\'ebrique, nous am\'eliorons encore ce r\'esultat, de la mani\`ere suivante.

\medskip Soit $\wh K$ un corps local, $\underline{G}$ un groupe
alg\'ebrique lin\'eaire connexe semi-simple, d\'efini sur ce corps,
$\underline{S}$ un tore $\wh K$-d\'eploy\'e ma\-xi\-mal, et $\Ga$ un
r\'eseau de $G=\underline{G}(\wh K)$. Les propri\'et\'es dynamiques et
ergodiques de l'action de $S=\underline{S}(\wh K)$ par translations
\`a droite sur l'espace quotient $\Ga\backslash G$ font actuellement
l'objet de nombreuses \'etudes (voir par exemple
\cite{Mar1,Zim,Mar2,Tom,LW}).  Nous nous int\'eresserons dans cet
article au cas o\`u $\underline{S}$ est de $\wh K$-rang $1$ et $\wh K$
est non archim\'edien, surtout dans la situation peu \'etudi\'ee o\`u
$\Ga$ est non uniforme (l'existence d'un tel $\Ga$ implique que $\wh
K$ est isomorphe \`a un corps de s\'eries formelles de Laurent sur un
corps fini).  Pour $\wh K=\FF_q((X^{-1}))$, $\underline{G}={\rm
  PGL}_2$, $\underline{S}$ le sous-groupe diagonal et $\Ga={\rm
  PGL}_2(\FF_q[X])$, la situation a \'et\'e compl\`etement d\'ecrite
dans \cite{BP}, en termes arithm\'etiques.

Si $M$ est le sous-groupe compact maximal de $S$, il revient presqu'au
m\^eme (voir par exemple \cite{Moz,LP}) d'\'etudier l'action par
translations \`a droite du groupe $S/M$ sur l'espace $\Ga\backslash
G/M$.  Celle-ci s'interpr\`ete en termes d'actions de groupes sur des
arbres, de la mani\`ere suivante. Soit $\TT$ l'arbre de Bruhat-Tits
\cite{BT} de $(\underline{G},\wh K)$ (biparti, de sommets bleus ou
verts). Alors $G$ agit transitivement (par translation au but) sur le
sous-espace $\G_0 \TT$ de $\G \TT$ form\'e des g\'eod\'esiques
d'origine un sommet vert, et $M$ est le stabilisateur d'un point de
$\G_0 \TT$.  L'action \`a droite de $S/M$, qui est isomorphe \`a $\ZZ$,
sur $G/M$, qui s'identifie \`a $\G_0 \TT$, correspond \`a l'action des
puissances paires de la transformation g\'eod\'esique.

Lorsque $\Ga$ est uniforme, il est connu (voir par exemple \cite{CP})
que l'action de $S/M$ sur $\Ga\backslash G/M$ est Bernoulli pour la
mesure naturelle sur $\Ga\backslash G/M$ venant de la mesure de Haar
sur $G$ (voir par exemple \cite{HK} pour les d\'efinitions et rappels
de th\'eorie ergodique). Nous g\'en\'eralisons ce r\'esultat au cas non
uniforme.

\btheo\label{theo:intro} %
Pour tout r\'eseau $\Ga$ de $G=\underline{G}(\wh K)$, l'action par
translations \`a droite de $S/M$ sur $\Ga\backslash G/M$ est Bernoulli
d'entropie finie.
\etheo

Nous montrons en fait un r\'esultat (voir le th\'eor\`eme
\ref{theo:geodBernou}) valable pour de nombreux sous-groupes
g\'eom\'etriquement finis d'automorphismes d'arbres localement finis
au sens de \cite{Rob,Pau}.

Plus g\'en\'eralement, \'etant donn\'e un arbre $T$ et un sous-groupe
d'automorphismes $\Ga$ de $T$, nous nous int\'eresserons au codage du
flot g\'eod\'esique sur $\Ga\backslash \G T$. Nous donnons dans la
partie \ref{sec:exappendice} des codages intrins\`eques, au sens o\`u
ils n'utilisent que la structure de graphe de groupes quotient (au
sens de \cite{Ser}) de $T$ par $\Ga$. Les propri\'et\'es canoniques de
cette construction devraient \^etre utiles (voir par exemple
\cite{LP}).  Ces codages markoviens (sur des alphabets
\'eventuellement infinis) sont obtenus pour le cas d'actions
$k$-acylindriques au sens de Sela \cite{Sel} de n'importe quel groupe
$\Ga$ sur n'importe quel arbre simplicial $T$ (en particulier sans
supposer $T$ localement fini, et sans supposer finis les
stabilisateurs de sommets dans $\Ga$), voir le th\'eor\`eme
\ref{prop:flotgeod_ordqqc}.

Les r\'eseaux non uniformes du th\'eor\`eme \ref{theo:intro}
n'agissent pas de mani\`ere acylindrique sur leur arbre de
Bruhat-Tits, mais nous montrons dans la partie \ref{subseq:actacy}
comment modifier ces actions pour les rendre acylindriques.  En
particulier, l'action modifi\'ee de l'action de ${\rm
  PGL}_2(\FF_q[X])$ sur l'arbre de Bruhat-Tits de $({\rm
  PGL}_2,\FF_q((X^{-1})))$ est $5$-acylindrique, et donne lieu \`a un
codage par la m\'ethode g\'en\'erale, qui est tr\`es proche du codage
particulier obtenu dans \cite{BP}.

La partie \ref{sec:exappendice} de cet article a \'et\'e \'ecrite
avant la partie 4.1 de \cite{Pau}, o\`u le second auteur \'etudie
d'autres propri\'et\'es dynamiques du flot g\'eod\'esique sur un
arbre, et en particulier certains arguments de la partie 4.1 de
\cite{Pau} ont \'et\'e inspir\'es de ceux de la partie
\ref{sec:exappendice}, et pas inversement. Il faut remarquer que
lorsque l'on autorise de la torsion dans les r\'eseaux, le flot
g\'eod\'esique n'est pas a priori markovien. C'est pr\'ecis\'ement
pour obtenir un caract\`ere markovien (et donc un codage par une
dynamique symbolique) que nous avons introduit un ``flot
g\'eod\'esique d'ordre $k$'' sur un arbre dans la partie
\ref{sec:exappendice}.

\medskip
\noindent {\small {\it Remerciements : } Nous remercions
  J.-P.~Thouvenot pour son aide pr\'ecieuse, en par\-ti\-cu\-lier
  concernant les r\'ef\'erences, ainsi que S.~Mozes et F.~Ledrappier.
  Nous remercions le rapporteur anonyme d'une version pr\'ec\'edente
  de cet article, certains de ses commentaires nous ont permis de
  d\'emontrer le th\'eor\`eme \ref{theo:introneuf}.}

\section{Notations et rappels}
\label{sec:rappels}

Nous renvoyons \`a \cite{Ser,Coo,Pau,Rob} pour des preuves et
compl\'ements concernant cette partie.  Pour toute action d'un groupe
$\Ga$ sur un ensemble, nous notons $\Ga_x$ le stabilisateur d'un point
$x$.  Par boule d'un espace m\'etrique, nous entendons boule ferm\'ee.

\subsection{Graphes de groupes et flot g\'eod\'esique sur un arbre}
\label{subsec:graphgroupflogeod}

Si $X$ est un graphe, on note $VX$ l'ensemble de ses sommets et $EX$
l'ensemble de ses ar\^etes. Pour toute ar\^ete $e$, on d\'esigne par
$o(e)$ son sommet origine, $t(e)$ son sommet terminal et
$\overline{e}$ son ar\^ete oppos\'ee. Les longueurs d'ar\^etes (des
r\'ealisations g\'eom\'etriques) sont suppos\'ees \'egales \`a $1$.

On appelle {\it graphe de groupes}, et on note $(X,G_\ast)$, la
donn\'ee des objets suivants~:
\begin{itemize}
\item un graphe $X$ (suppos\'e connexe dans la suite);
\item pour tout sommet $v$ de $X$, un groupe $G_v$;
\item pour toute ar\^ete $e$ de $X$, un groupe $G_e$, tel que
    $G_e=G_{\overline{e}}$;
\item pour toute ar\^ete $e$ de $X$, un morphisme injectif
    $\rho_e:G_e\ra G_{t(e)}$.
\end{itemize}

\medskip %
Par exemple, si $\Ga$ est un groupe d'automorphismes (sans inversion)
d'un arbre $T$, alors le graphe $X=\Gamma\backslash T$ est muni d'une
structure de graphe de groupes, appel\'ee {\it graphe de groupes
  quotient} et not\'ee $\Ga\bac T$. On proc\`ede ainsi pour la
construire. On fixe un relev\'e $\wt v$ dans $T$ de chaque sommet $v$
de $X$, un relev\'e $\wt e$ dans $T$ de chaque ar\^ete $e$ de $X$, on
impose que $\overline{\wt e}=\wt{\overline{e}}$ et on fixe un
\'el\'ement $g_e$ de $\Gamma$ tel que $g_e\wt {t(e)}=t(\wt e)$. On
d\'efinit alors $G_e$ et $G_v$ comme les fixateurs dans $\Gamma$ de
$\wt e$ et $\wt v$.  Alors $\rho_e:G_e\ra G_{t(e)}$ est d\'efinie
comme la restriction \`a $G_e$ de la conjugaison par $g_e^{-1}$.  Le
graphe de groupes quotient $\Ga\bac T$ ne d\'epend pas (\`a
isomorphisme de graphes de groupes pr\`es), du choix des $\wt e$, $\wt
v$ et $g_e$ (voir \cite{Ser} pour tout compl\'ement). Si $T$ est
localement fini et $\Ga$ est discret, alors les groupes $G_e$ et $G_v$
sont finis. Si de plus $\Ga\backslash T$ est fini, alors $\Ga\bac T$
est un graphe (connexe) fini de groupes finis.

\medskip %
Soit $T$ un arbre simplicial, muni de sa topologie faible.  L'{\it
  espace des g\'eod\'esiques} de $T$ est l'espace $\G T$ des
applications simpliciales injectives de $\RR$ dans $T$ (avec $\RR$
muni de sa structure simpliciale usuelle d'ensemble de sommets $\ZZ$),
muni de la topologie compacte-ouverte.  Comme $T$ est un arbre, la
condition d'injectivit\'e est \'equivalente \`a la condition
d'injectivit\'e locale. Notons ${\rm Aut}(T)$ son groupe
d'automorphismes sans inversion. Il est localement compact pour la
topologie compacte-ouverte si $T$ est localement fini.

Le groupe $\ZZ$ agit sur $\G T$ par translations \`a la source
$(n,f)\mapsto\{x\mapsto f(x+n)\}$. Le groupe ${\rm Aut}(T)$ agit par
hom\'eomorphismes sur $\G T$ par composition au but $(\ga,f) \mapsto
\{x\mapsto \ga f(x)\}$. Ces deux actions commutent.  Appelons {\it
  transformation g\'eod\'esique sur $\G T$} l'application $\wt\varphi:
\G T\ra \G T$ d\'efinie par $\ell\mapsto \{t\mapsto\ell(t+1)\}$.
Appelons {\it renversement du temps sur $\G T$} l'application
$\wt\tau:\G T\ra \G T$ d\'efinie par $\ell\mapsto \{t \mapsto
\ell(-t)\}$.

Soit $\Ga$ un sous-groupe de ${\rm Aut}(T)$. On munit les quotients
$\Ga \backslash T$ et $\Ga\backslash\G T$ de la topologie quotient. On
note $\pi:T\ra \Ga\backslash T$ et $\pi':\G T\ra \Ga\backslash\G T$
les projections canoniques.  L'application $\wt\varphi$ induit une
application continue $\varphi:\Ga\backslash\G T\ra \Ga\backslash\G T$,
aussi appel\'ee la {\it transformation g\'eod\'esique} sur
$\Ga\backslash\G T$.  L'application $\wt\tau$ induit une application
$\tau:\Ga \backslash \G T\ra \Ga\backslash\G T$, que nous appelerons
{\it renversement du temps sur $\Ga\backslash\G T$}.

\subsection{$\!$Groupes g\'eom\'etriquement finis et mesure de
  Bowen-Margulis}
\label{subsec:groupgeofinmesbowmar}

Soit $T$ un arbre localement fini. Notons $T\cup\partial T$ la
compactification par l'espace des bouts de $T$, et $\partial_2 T$ le
produit $\partial T\times\partial T$ priv\'e de sa diagonale.
Rappelons que toute ar\^ete $e$ de $T$ d\'efinit une partition en deux
parties $\partial_e T$ et $^c\partial_e T$ de $\partial T$, de sorte
que toute droite g\'eod\'esique d'origine dans $\partial_e T$ et
d'extr\'emit\'e dans $^c\partial_e T$ parcourt $e$ suivant
l'orientation de $e$. Nous dirons que $T$ est {\it uniforme} s'il
existe un sous-groupe discret dans ${\rm Aut}(T)$ tel que le graphe
$\Ga\backslash T$ soit fini. Par exemple, un arbre r\'egulier ou
bi-r\'egulier est uniforme.

Notons $x_0$ un point base de $T$. L'{\it entropie volumique} de $T$,
qui ne d\'epend pas de $x_0$, est
$$\delta_T=\limsup_{n\ra\infty} 
\frac{1}{n}\log {\rm Card}(B(x_0,n)\cap VT)\;.
$$

\medskip
Soit $\Ga$ un sous-groupe discret de ${\rm Aut}(T)$. En particulier,
l'action de $\Ga$ sur $\G T$ est proprement discontinue (mais pas
forc\'ement libre en g\'en\'eral), et donc $\Ga\backslash\G T$ est
localement compact.

Le groupe $\Ga$ est dit {\it non \'el\'ementaire} s'il ne pr\'eserve ni
point ni paire de points de $T\cup \partial T$. Il existe alors un
unique plus petit sous-arbre $\Ga$-invariant non vide, not\'e
$T_{\Ga,{\rm min}}$.

On appelle {\it rayon cuspidal de groupes} un graphe de groupes finis
$(R,G_*)$ avec $R$ un rayon, de suite des ar\^etes cons\'ecutives
$(e_n)_{n\in\NN}$ orient\'ees vers le bout de $R$, tel que pour tout $n$
dans $\NN-\{0\}$, le morphisme $G_{e_n}\ra G_{o(e_n)}$ soit surjectif.
Le groupe $\Gamma$ est dit {\it g\'eom\'etriquement fini} s'il est non
\'el\'ementaire et si le graphe de groupes quotient $\Ga\bac T_{\Ga,{\rm
    min}}$ est r\'eunion d'un graphe fini de groupes finis et, recoll\'es
en leurs extr\'emit\'es, d'un nombre fini de rayons cuspidaux de groupes.
Voir \cite{Pau} pour l'\'equivalence avec la d\'efinition dynamique
usuelle (comme dans \cite{Rob}), et des d\'eveloppements.

Nous renvoyons par exemple \`a \cite{BH} pour la d\'efinition des
horoboules (ferm\'ees par d\'efaut) dans un espace m\'etrique
g\'eod\'esique CAT$(0)$. Comme montr\'e dans \cite{Pau}, la pr\'eimage
des rayons cuspidaux maximaux dans $T_{\Ga,{\rm min}}$ forme alors une
famille disjointe $\Ga$-invariante maximale d'horoboules ouvertes. Les
points \`a l'infinis de ces horoboules, qui sont donc les
extr\'emit\'es des rayons g\'eod\'esiques relevant les rayons
cuspidaux, seront appel\'es les {\it points paraboliques born\'es} de
$\Ga$ (voir \cite{Rob,Pau} pour l'explication dynamique).

Par exemple, soit $\wh K$ un corps local non archim\'edien,
$\underline{G}$ un groupe alg\'ebrique lin\'eaire connexe semi-simple
sur $\wh K$, de $\wh K$-rang $1$. Soit $\Ga$ un r\'eseau de
$G=\underline{G}(\wh K)$. Soit $T$ l'arbre de Bruhat-Tits de
$(\underline{G},\wh K)$. Alors, par un th\'eor\`eme de A.~Lubotzky
\cite{Lub}, l'action de $\Ga$ sur $T$ est g\'eom\'etriquement finie et
$T=T_{\Ga,{\rm min}}$.

\medskip Soit $\Ga$ un sous-groupe g\'eom\'etriquement fini de ${\rm
  Aut}(T)$. Nous dirons que $\Ga$ {\it poss\`ede la propri\'et\'e de
  Selberg} s'il admet un sous-groupe d'indice fini, dont tout
\'el\'ement de torsion est conjugu\'e \`a un \'el\'ement du groupe
d'un sommet int\'erieur d'un sous-rayon cuspidal de $\Ga\bac T$. Par
exemple, c'est vrai si $\Ga\bac T$ n'a pas de sous-rayon cuspidal
(voir \cite{Ser}).  Par le lemme de Selberg \cite{Alp}, c'est aussi
vrai pour $\Ga$ un r\'eseau de $G=\underline{G}(\wh K)$ agissant sur
l'arbre de Bruhat-Tits de $(\underline{G},\wh K)$, avec $\wh K$ un
corps local non archim\'edien et $\underline{G}$ un groupe
alg\'ebrique lin\'eaire connexe semi-simple sur $\wh K$, de $\wh
K$-rang $1$. Cette propri\'et\'e est aussi v\'erifi\'ee si tout
stabilisateur de point parabolique born\'e dans $\Ga$ est
r\'esiduellement fini, car on peut alors utiliser le r\'esultat
susnomm\'e de \cite{Ser} pour enlever la torsion sur le graphe priv\'e
de ses rayons cuspidaux, et recoller des rayons cuspidaux
correspondant \`a des sous-groupes d'indice fini des stabilisateurs de
points paraboliques born\'es. Rappelons qu'un groupe est {\it
  r\'esiduellement fini} si l'intersection de ses sous-groupes
d'indice fini est r\'eduite \`a l'\'el\'ement neutre.

\bigskip 
L'{\it exposant critique} $\delta=\delta_\Ga$ de $\Ga$, qui ne d\'epend
pas de $x_0$, est l'\'el\'ement de $[0,+\infty]$ tel que la {\it s\'erie de
  Poincar\'e} de $\Ga$
$$P(s)=P_{\Ga,x_0}(s)=\sum_{\ga\in\Ga} e^{-s d(x_0,\ga x_0)}$$
converge pour $s>\delta$ et diverge pour $s<\delta$. Le groupe $\Ga$
est {\it de type divergent} si sa s\'erie de Poincar\'e $P(s)$ diverge
pour $s=\delta$. 

Si $\Ga$ est de type divergent d'exposant critique fini non nul, alors
il existe (voir \cite{Coo}) une famille $(\mu_x)_{x\in VT}$ de mesures
finies sur $\partial T$, appel\'ee {\it mesure de Patterson-Sullivan},
unique \`a scalaire multiplicatif pr\`es, de supports $\partial
T_{\Ga,{\rm min}}$, telle que
\begin{itemize}
\item $\forall\;\ga\in\Ga\;,\;\;\ga_*\mu_x=\mu_{\ga x}$,
\item $\forall\;x,y\in VT,\;\forall\;\xi\in\partial
  T\;,\;\;\frac{d\mu_x}{d\mu_y}(\xi)=e^{-\delta \beta_\xi(x,y)}$, 
\end{itemize}
o\`u $\beta_\xi(x,y)= d(x,z)-d(z,y)$ pour tout sommet $z$
suffisamment proche de $\xi$.
  
Pour toute g\'eod\'esique $\ell$, notons $\ell_-,\ell_+$ les points de
$\partial T$ origine et extr\'emit\'e de $\ell$. L'application $\G
T\ra\partial_2T\times\ZZ$, qui \`a $\ell$ associe $(\ell_-,\ell_+,t)$
avec $t$ la distance alg\'ebrique sur $\ell$ entre $\ell(0)$ et le
point de $\ell$ le plus proche de $x_0$, est un hom\'eomorphisme. Ce
param\'etrage d\'epend du point base $x_0$. Si $\Ga$ est de type
divergent d'exposant critique fini non nul, on d\'efinit une mesure
$\wt{m}_{\mbox{\tiny BM}}$ sur $\G T$, appel\'ee {\it mesure de
  Bowen-Margulis}, par
$$d\wt m_{\mbox{\tiny BM}}(\ell_-,\ell_+,t)=
\frac{d\mu_{x_0}(\ell_-)d\mu_{x_0}(\ell_+)dt}
{d_{x_0}(\ell_-,\ell_+)^{2\delta}}\;,$$ o\`u $d_{x_0}$ est la distance
sur $\partial T$ d\'efinie par $d_{x_0}(\ell_-,\ell_+)=e^{-u}$ avec
$u$ la longueur de l'intersection des rayons g\'eod\'esiques issus de
$x_0$ convergeant vers $\ell_-$ et $\ell_+$. (L'origine de cette mesure
remonte aussi aux travaux de Patterson-Sullivan, mais nous
pr\'ef\'erons donner des noms diff\'erents \`a deux mesures
diff\'erentes, la mesure de Patterson-Sullivan qui vit sur $\partial
T$ et celle de Bowen-Margulis qui vit sur $\G T$.) La mesure $\wt
m_{\mbox{\tiny BM}}$ ne d\'epend pas de $x_0$. Elle est invariante par
la transformation g\'eod\'esique sur $\G T$ et par $\Ga$. Elle induit
donc une mesure $m_{\mbox{\tiny BM}}$ sur $\Ga\backslash\G T$,
appel\'ee {\it mesure de Bowen-Margulis sur $\Ga\backslash\G T$}. Le
support de $m_{\mbox{\tiny BM}}$ est $\Ga\backslash\G T_{\Ga,{\rm
    min}}$.  Lorsque $\Ga$ est cocompact, cette mesure est la mesure
d'entropie maximale pour $\varphi$ (voir \cite{CP,Kai,Bou}).  Lorsque
$T$ est un arbre de Bruhat-Tits comme dans l'introduction, alors la
restriction de $\wt m_{\mbox{\tiny BM}}$ \`a $\G_0 T$ (d\'efini dans
l'introduction) s'identifie (\`a un scalaire multiplicatif pr\`es)
avec l'image dans $G/M$ de la mesure de Haar de $G$.

\section{Finitude de la mesure et m\'elange}
\label{sec:finibowmarg}

Les r\'esultats de ce paragraphe d\'ecoulent essentiellement de
r\'esultats connus (voir \cite{BM,Rob}). En particulier, les trois
premi\`eres assertions de la proposition suivante d\'ecoulent du
Corollary 6.5 de \cite{BM}, l'avant-derni\`ere de la Proposition 7.3
de \cite{BM} et la derni\`ere assertion, en utilisant l'avant
derni\`ere et la troisi\`eme, d\'ecoule de \cite{Rob} (adaptant des
id\'ees de \cite{DOP}).  Nous ne donnons la preuve regroup\'ee que par
souci de compl\'etude.

\bprop\label{prop:finitudebowmarg} %
Soit $T$ un arbre localement fini et $\Ga$ un sous-groupe
g\'e\-o\-m\'e\-tri\-que\-ment fini de ${\rm Aut}(T)$. Si $T_{\Ga, \rm
  min}$ est uniforme, alors
\begin{itemize}
\item $\delta_\Ga$ est fini, non nul;
\item $\delta_\Ga=\delta_{T_{\Gamma,{\rm min}}}$;
\item $\Ga$ est de type divergent;
\item pour tout point parabolique born\'e $\xi$ de $\partial T$ et tout
  point $y_0$ dans $VT$, il existe une constante $c\geq 1$ telle que
  $\frac{1}{c}\; e^{n\delta_{\Gamma}}\leq {\rm Card~} \big(B(y_0,2n)\cap
  \Ga_\xi \,y_0\big)\leq c \;e^{n\delta_{\Gamma}}$ pour tout $n$ dans $\NN$;
\item pour tout point parabolique born\'e $\xi$ de $\partial T$, on a
  $\delta_{\Ga_\xi}=\delta_\Ga/2$;
\item la mesure de Bowen-Margulis sur $\Ga\backslash\G T$ est finie.
\end{itemize}
\eprop

\noindent{\bf Remarques.} (1)
L'hypoth\`ese que $T_{\Ga, \rm min}$ est uniforme ne peut \^etre omise.
En effet, pour tout $n$ dans $\NN$, posons $q_n=2^{2^n}$ et
$\Ga_n=(\ZZ/2\ZZ)^{q_n}$. Notons que $\Ga_n$ s'injecte naturellement
dans $\Ga_{n+1}$. Consid\'erons le graphe de groupes suivant (de type
Nagao au sens de \cite{BL}), qui est g\'eom\'etriquement fini (et de
volume fini)~:
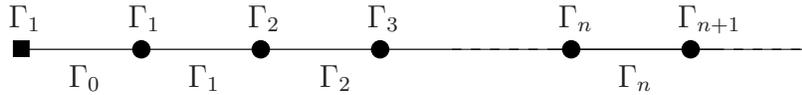
\begin{center}
\input{fig_nagaocroa.pstex_t}
\end{center}
Soit $T$ le rev\^etement universel de ce graphe de groupes (qui est un
arbre non uniforme, car de valences non uniform\'ement born\'ees) et
$\Ga$ son groupe fondamental, pour des choix indiff\'erents de points
bases (voir \cite{Ser}).  Si $x_0$ est un sommet de $T$, pr\'eimage de
l'origine de ce rayon de groupes, alors la boule de rayon $2n$ et de
centre $x_0$ contient au moins $q_n/q_{n-1}-1=q_{n-1}-1$ points, donc
$$\delta_\Ga\geq \limsup_{n\ra\infty} \frac{1}{2n}\log
(q_{n-1}-1)=+\infty\;,$$
et l'exposant critique de $\Ga$ est infini.

\medskip 
(2) Soit $T$ un arbre localement fini et $\Ga'$ un
sous-groupe discret de ${\rm Aut~}T$. Alors $\delta_{\Ga'}\leq
\delta_{T}$, car
$$P_{\Ga',x_0}(s)=\big({\rm Card~}\Ga'_{x_0}\big)\;\sum_{y\in\Ga x_0}
e^{-s d(x_0,y)}\leq\big({\rm Card~}\Ga'_{x_0}\big)\;\sum_{y\in VT}
e^{-s d(x_0,y)}\;,$$
qui converge si $s>\delta_{T}$. Si $\Ga'$ est non
\'el\'ementaire, alors $\delta_{\Ga'}$ est non nul, car alors $\Ga'$
contient au moins un groupe libre de rang $2$ de Schottky (voir par
exemple \cite{Lub}).  Si les valences de $T$ sont uniform\'ement
born\'ees, disons par $q+1$, alors $\delta_{\Ga'}$ est fini, car
$\delta_{T}\leq\log q$.

\medskip
\noindent{\it D\'emonstration de la proposition \ref{prop:finitudebowmarg}.}
Quitte \`a remplacer $T$ par $T_{\Ga,{\rm min}}$, nous pouvons supposer
que $T=T_{\Ga,{\rm min}}$.

Si $T$ est un arbre uniforme ayant au moins trois bouts, il est bien
connu (voir par exemple \cite{Coo,Bou,Rob}) que $\delta_{T}$ est fini
et non nul, qu'il existe une constante $c_1\geq 1$ telle que pour tout
$x$ dans $VT$ et $n$ dans $\NN-\{0\}$, pour toute composante connexe
$C$ de $T-\{x\}$,
$$\frac{1}{c_1}\;e^{n\delta_{T}}\leq {\rm Card~} \left(B(x,n)\cap
  VT\cap C\right)\leq c_1 \;e^{n\delta_{T}}\;,\;\;\;\;(*)$$
et donc que la s\'erie $\sum_{y\in VT} e^{-s d(x_0,y)}$ diverge en
$s=\delta_{T}$. En particulier, la premi\`ere assertion de la
proposition \ref{prop:finitudebowmarg} d\'ecoule de la seconde remarque
ci-dessus.

\medskip %
Soit $\xi$ un point parabolique born\'e de $\partial T$. Soit
$(x_i)_{i\in\NN}$ la suite des sommets cons\'ecutifs d'un rayon
g\'eod\'esique d'extr\'emit\'e $\xi$, se projetant sur un rayon
cuspidal dans $\Ga\backslash T$. Alors l'intersection avec $\Ga_\xi
\,x_0$ de la sph\`ere $S(x_0,2n)$ est la r\'eunion des $C\cap
S(x_n,n)$, o\`u $C$ est une composante connexe de $T-\{x_n\}$ qui ne
contient ni $x_{n-1}$ ni $x_{n+1}$. Comme $T$ est de valences
uniform\'ement born\'ees, la quatri\`eme assertion de la proposition
\ref{prop:finitudebowmarg} d\'ecoule de $(*)$. Cette quatri\`eme
assertion implique en particulier que $\delta_{\Ga_\xi}=\delta_T/2$.

\medskip %
Si $T_0$ est le graphe obtenu en enlevant \`a l'arbre $T$ la r\'eunion
de la famille disjointe $\Ga$-invariante maximale d'horoboules
ouvertes, alors $VT_0$ ne contient qu'un nombre fini d'orbites sous
$\Ga$. Par cons\'equent, la s\'erie de Poincar\'e de $\Ga$ diverge si
et seulement si la s\'erie $\sum_{y\in VT_0} e^{-s d(x_0,y)}$ diverge.

Montrons qu'il existe une constante $c_2\geq 1$ telle que, pour toute
horosph\`ere $H$ dans $T$, passant par un sommet de $T$, bord d'une
horoboule $H\!B$ ne contenant pas $x_0$, et pour tout $n$ dans $\NN$,
$${\rm Card~}(B(x_0,n)\cap H) \leq {\rm Card~}(B(x_0,n)\cap H\!B \cap
VT) \leq c_2\;{\rm Card~}(B(x_0,n)\cap H)\;.\;\;\;\;(**)$$ En effet,
soit $\xi_0$ le point \`a l'infini de $H\!B$. On peut supposer que
l'int\'erieur de $H\!B$ rencontre $B(x_0,n)$.  Notons $p$ la distance
de $x_0$ \`a $H\!B$, et $y$ le point du rayon g\'eod\'esique entre
$x_0$ et $\xi_0$, \`a distance $\frac{n+p}{2}$ de $x_0$. Alors
$B(x_0,n)\cap H\!B= B(y,\frac{n-p}{2})$, et si $B$ est la r\'eunion
des composantes connexes de $B(y,\frac{n-p}{2})-\{y\}$ rencontrant
$H$, alors $B(x_0,n)\cap H=B\cap S(y,\frac{n-p}{2})$. Or par $(*)$, il
existe une constante $c_3\geq 1$, ne d\'ependant pas de $y,n,p$ telle
que
$${\rm Card~}\left(B\cap S(y,\frac{n-p}{2})\right)\geq
\frac{1}{c_3}\;{\rm Card~}\left(B(y,\frac{n-p}{2})\cap VT\right)\;.$$
L'affirmation $(**)$ s'en d\'eduit.

Il d\'ecoule de $(**)$ que la s\'erie $\sum_{y\in VT_0} e^{-s
  d(x_0,y)}$ diverge si et seulement si la s\'erie $\sum_{y\in VT}
e^{-s d(x_0,y)}$ diverge. Ceci montre que $\delta_\Ga=\delta_{T}$ et
que $\Ga$ est de type divergent.

\medskip %
Pour montrer que la mesure de Bowen-Margulis est finie, d'apr\`es le
th\'eor\`eme B de \cite{DOP}, comme remarqu\'e dans
\cite[th\'eo.~1.11]{Rob}, il suffit de montrer que $\Ga$ est de type
divergent, ce que nous venons de faire, et que pour tout point
parabolique born\'e $\xi$ dans $\partial T$, l'exposant critique
$\delta_{\Ga_\xi}$ est strictement inf\'erieur \`a $\delta_{\Ga}$.  Or
$\delta_{\Ga_\xi}=\delta_{T}/2= \delta_{\Ga}/2<\delta_{\Ga}$, ce qui
montre le r\'esultat.  \cqfd

\medskip %
Dans la suite de cet article, si les hypoth\`eses de la proposition
\ref{prop:finitudebowmarg} sont v\'erifi\'ees, nous supposerons,
quitte \`a normaliser, que $m_{\mbox{\tiny BM}}$ est une mesure de
probabilit\'e.

\bigskip Apr\`es la propri\'et\'e de finitude de la mesure de
Bowen-Margulis, regardons celle de m\'elange.  Soit $T$ un
arbre localement fini et $\Ga$ un sous-groupe discret non
\'el\'ementaire de ${\rm Aut}(T)$.

M\^eme lorsque $\Ga$ est cocompact, la transformation g\'eod\'esique
$\varphi$ sur $\Ga\backslash \G T$ n'est pas forc\'ement
m\'elangeante.  Remarquons par exemple que si $T^{(1)}$ est la
premi\`ere subdivision barycentrique de $T$, alors $\Ga$ est encore un
sous-groupe de ${\rm Aut}(T^{(1)})$, mais la transformation
g\'eod\'esique sur $\Ga\backslash \G T^{(1)}$ n'est pas m\'elangeante.
De plus, la propri\'et\'e de m\'elange de $\varphi$ n'est pas
invariante par passage \`a un sous-groupe d'indice fini. Par exemple,
la transformation g\'eod\'esique pour le bouquet de deux cercles est
m\'elangeante, mais pas celle pour son rev\^etement connexe \`a deux
feuillets o\`u aucun des deux cercles n'est relevable.

Il existe un crit\`ere assez pratique pour v\'erifier que la
transformation g\'eod\'esique est m\'elangeante.

Notons $L_\Ga$ le sous-groupe de $\ZZ$ engendr\'e par les distances de
translation des \'el\'ements de $\Ga$ dans $T$. Par exemple, si $T$
est l'arbre de Bruhat-Tits d'un groupe alg\'ebrique lin\'eaire connexe
semi-simple $\underline{G}$ de rang $1$ sur un corps local non
ar\-chi\-m\'e\-dien $\wh K$, et si $\Ga$ est un r\'eseau de
$\underline{G}(\wh K)$, alors le groupe $L_\Ga$ vaut $2\ZZ$.

Si $T$ n'a pas de sous-arbre invariant non vide et propre, et n'a pas
de sommet de valence $2$, alors $L_\Ga$ vaut $\ZZ$ ou $2\ZZ$ (ces
hypoth\`eses sont b\'enignes, car on peut toujours passer au
sous-arbre invariant minimal $T_{\Ga,{\rm min}}$, et lui enlever les
sommets de valence $2$, et elles sont pr\'eserv\'ees par passage ˆ un
sous-groupe d'indice fini de $\Ga$). Une des mani\`eres de d\'emontrer
cette affirmation est d'introduire le sous-groupe $\Lambda_\Gamma$ de
$\ZZ$ engendr\'e par les distances entre sommets de valence au moins
$3$ de $T$, qui est donc \'egal \`a $\ZZ$ sous notre hypoth\`ese, et
de remarquer avec \cite[page 564]{GL} que
$$2\Lambda_\Gamma\subset L_\Ga\subset \Lambda_\Ga\;.$$
Remarquons quand m\^eme que la suppression des sommets de valence $2$
peut avoir un certain effet dans le cas alg\'ebrique. Soit $T$ l'arbre
de Bruhat-Tits de $(\underline{G},\wh K)$ comme dans l'introduction,
et soit $\Ga$ un r\'eseau de $\underline{G}(\wh K)$. Supposons qu'il
existe une (alors unique) $\underline{G}(\wh K)$-orbite de sommets de
valence $2$ (ce qui est le cas par exemple pour $\underline{G}={\rm
  PGL}_2$). Alors apr\`es suppression de ces sommets, le groupe
$L_\Ga$ devient $\ZZ$, et l'action de $S/M$ sur $\Ga\backslash G/M$
s'identifie maintenant exactement avec l'action de la transformation
g\'eod\'esique.

Si $L_\Ga=2\ZZ$, fixons $x_0$ un point base de $T$, et notons $\G_0 T$
le sous-espace de $\G T$ form\'e des g\'eod\'esiques dont l'origine
est \`a distance paire de $x_0$. Il est facile de voir que $\G_0 T$
est invariant par $\Ga$ et par $\varphi^2$. Pour l'invariance par
$\Ga$, on remarque que si $A_\ga$ est l'axe de translation ou
l'ensemble des points fixes d'un \'el\'ement $\ga$ dans $\Ga$ de
distance de translation $\lambda(\ga)$, alors $d(x,\ga
x)=2d(x,A_\ga)+\lambda(\ga)$ (voir par exemple \cite{Ser}), donc
$\lambda(\ga)$ est impair si $d(x,\ga x)$ l'est. Pour l'invariance par
$\varphi^2$, on remarque que pour $x,y,z$ trois points d'un arbre, si
$d(x,y)$ et $d(y,z)$ sont pairs, alors $d(x,z)$ l'est. Si
$L_\Ga=2\ZZ$, nous munirons $\Ga\backslash\G_0 T$ de la restriction de
la mesure de Bowen-Margulis, normalis\'ee pour \^etre de
probabilit\'e.

Le r\'esultat suivant d\'ecoule alors de \cite[Theo.~3.1]{Rob} (en
fait d'une version discr\`ete de ce th\'eor\`eme). Lorsque $\Ga$ est
un r\'eseau de $G=\underline{G}(\wh K)$ agissant sur l'arbre de
Bruhat-Tits de $(\underline{G},\wh K)$, avec $\wh K$ un corps local
non archim\'edien et $\underline{G}$ un groupe alg\'ebrique lin\'eaire
connexe semi-simple sur $\wh K$, de $\wh K$-rang $1$, le second
\'enonc\'e d\'ecoule aussi du th\'eor\`eme de Howe-Moore \cite{Zim},
par les rappels de l'introduction.

\bprop\label{prop:melang} 
(1) Si $\Ga$ est g\'eom\'etriquement fini, si $T_{\Ga,{\rm min}}$
est uniforme et si $L_\Ga=\ZZ$, alors la transformation g\'eod\'esique
$\varphi$ est m\'elangeante sur $\Ga\backslash \G T$.

(2) Si $\Ga$ est g\'eom\'etriquement fini, si $T_{\Ga,{\rm min}}$ est
uniforme, sans sommet de valence $2$, et si $L_\Ga=2\ZZ$, alors le
carr\'e $\varphi^2$ de la transformation g\'eod\'esique est
m\'elangeant sur $\Ga\backslash \G_0 T$.  \cqfd 
\eprop

\section{Un premier codage g\'en\'eral}
\label{sec:codageneuf}

Fixons-nous les donn\'ees suivantes:
\begin{itemize}
\item[$\bullet$] $T$ un arbre localement fini, 
\item[$\bullet$] $\Ga$ un sous-groupe discret non
\'el\'ementaire de ${\rm Aut}(T)$, 
\item[$\bullet$] $(\mu_x)_{x\in VT}$ une mesure de Patterson-Sullivan
  pour $\Ga$ de dimension $\delta$, o\`u $0<\delta<+\infty$, 
telle que pour tout \'el\'ement elliptique $\ga$ de $\Ga-\{e\}$,
l'ensemble des points fixes de $\ga$ dans $\partial T$ soit de
mesure nulle pour $\mu_x$, 
\item[$\bullet$] $\G T$ l'espace des g\'eod\'esiques de $T$ et
  $\wt\varphi:\G T\ra\G T$ la transformation g\'eod\'e\-si\-que,
\item[$\bullet$] $\wt\psi=\wt\varphi\,^N$ avec $N\geq 1$, et 
\item[$\bullet$] $V'T$ une partie $\Ga$-invariante de $VT$ telle que,
  pour tout $\ell$ dans $$\G'T=\{\ell\in\G T\;:\; \ell(0)\in V'T\}$$
  et pour tout $t$ dans $\mathbb{Z}$, le sommet $\ell(t)$ appartient
  \`a $V'T$ si et seulement s'il existe $n$ dans $\mathbb{Z}$ tel que
  $\wt\psi^n\ell(0)=\ell(t)$. En particulier $\G'T$ est une partie
  mesurable $\Ga$-invariante de $\G T$, invariante par $\wt\psi$.
\end{itemize}  

\medskip
Notons ${\wt \mu}_{\mbox{\tiny BM}}$
la restriction \`a $\G'T$ de la mesure de Bowen-Margulis de $\Ga$
associ\'ee \`a $(\mu_x)_{x\in VT}$, $\psi:\Ga\backslash\G' T\ra
\Ga\backslash\G' T$ l'application induite par $\wt\psi$, et
${\mu}_{\mbox{\tiny BM}}$ la mesure sur l'espace quotient
$\Ga\backslash\G' T$ induite par ${\wt \mu}_{\mbox{\tiny BM}}$.

\bigskip La condition de mesure nulle de l'ensemble des points fixes
des \'el\'ements elliptiques non triviaux est tr\`es souvent v\'erifi\'ee, et
nous pensons qu'elle l'est toujours si $\Ga$ est un r\'eseau (uniforme
ou non) de ${\rm Aut}(T)$, avec $T$ uniforme. Donnons ci-dessous
quelques arguments pour d'une part \'etayer cet espoir, et d'autre part
montrer que le probl\`eme n'est pas compl\`etement trivial. La premi\`ere
remarque donne une condition suffisante, la seconde montre que les
ensembles de points fixes ne peuvent pas \^etre outranci\`erement gros,
mais la derni\`ere montre qu'ils peuvent quand m\^eme \^etre assez gros. 

Nous noterons ${\rm Fix}_T(\ga)$ et ${\rm Fix}_{\partial T}(\ga)$ les
ensembles de points fixes dans $T$ et $\partial T$ respectivement
d'une isom\'etrie $\ga$ de $T$, et $\delta_\ga=\delta_{{\rm
    Fix}_T(\ga)}$ l'entropie volumique de ${\rm Fix}_T(\ga)$.
Rappelons que ${\rm Fix}_T(\ga)$ est un sous-arbre de $T$
(\'eventuellement vide).

\brema \label{rem:condnecmesfixellnul} %
Si $T$ est uniforme, si $\Ga$ est un r\'eseau (uniforme ou non)
de ${\rm Aut}(T)$ d'exposant critique \'egal ˆ $\delta$, si $\ga\in\Ga$
v\'erifie $\delta_\ga<\delta$, alors $\mu_x({\rm Fix}_{\partial
  T}(\ga))=0$.  \erema

\dem Par unicit\'e (voir par exemple \cite{BM}), la mesure $\mu_x$ est
la mesure de Hausdorff de la distance visuelle $d_x$ (voir par exemple
\cite{Coo,BM}) sur $\partial T$, et $\delta$ est la dimension de
Hausdorff de $d_x$.  Donc si la dimension de Hausdorff $\delta'_\ga$
de ${\rm Fix}_{\partial T}(\ga)$ pour $d_x$ est strictement inf\'erieure
ˆ $\delta$, alors $\mu_x({\rm Fix}_{\partial T}(\ga))=0$. Or
$\delta'_\ga\leq \delta_\ga$ (ce qui montre le r\'esultat), car il est
classique que pour tout arbre localement fini $T'$, la dimension de
Hausdorff $\delta''$ de son bord pour sa distance visuelle est
inf\'erieure ˆ son entropie volumique $\delta'$.  (En effet, soit $x'_0$
dans $T'$. Pour tout $s>\delta'$, soit $(n_k)_{k\in\NN}$ une suite
strictement croissante dans $\NN$ telle que ${\rm Card~}B(x'_0,n_k)\leq
e^{sn_k}$. Soit $(\xi_i)_{1\leq i\leq p_k}$ une partie finie minimale
de $\partial T$ telle que $\partial T\subset \bigcup_{i=1}^{p_k}
B(\xi_i,e^{-n_k})$. Soit $x_i$ le point ˆ distance $n_k$ de $x'_0$ sur
le rayon g\'eod\'esique de $x'_0$ ˆ $\xi_i$.  Alors par minimalit\'e les
$x_i$ sont deux ˆ deux disjoints, donc sont au nombre de $e^{sn_k}$ au
plus.  Par cons\'equent, si $\mu_{s,e^{-n_k}}(\partial T)$ est la borne
inf\'erieure, sur tous les recouvrements finis de $\partial T$ par des
boules $B_j$ de rayon $r_j$ au plus $e^{-n_k}$, des nombres $\sum_j
r_j^s$, alors
$$\mu_{s,e^{-n_k}}(\partial T)\leq \sum_{i=1}^{p_k} e^{-sn_k}\leq 1<+\infty\;.$$
par d\'efinition de la dimension de Hausdorff, on a donc $\delta''\leq
s$, d'o\`u $\delta''\leq \delta'$).  
\cqfd

\medskip Par exemple, si $\Ga$ est un r\'eseau g\'eom\'etriquement fini de
$T$, dont tout \'el\'ement de torsion est conjugu\'e ˆ un \'el\'ement d'un
groupe de sommet d'un rayon cuspidal de $\Ga\bac T$, alors pour tout
$\ga$ dans $\Ga-\{0\}$, l'ensemble ${\rm Fix}_T(\ga)$ est contenu dans
une horoboule de $T$, et donc son entropie volumique est au plus
$\delta/2$, par la proposition \ref{prop:finitudebowmarg}, donc
$\mu_x({\rm Fix}_{\partial T}(\ga))$ est bien nulle par la remarque
pr\'ec\'edente.

\brema 
Si $T$ est minimal, alors aucun \'el\'ement non trivial de $\Ga$ ne
fixe point par point un ouvert non vide de $\partial T$.
\erema

\dem %
Supposons sinon que $\alpha$ soit un tel \'el\'ement, et que $U$ soit
un tel ouvert. Pour tout $\epsilon>0$ et tout $\xi$ dans $\partial T$,
par densit\'e dans $\partial T_{\Ga,{\rm min}}\times\partial T_{\Ga,{\rm
    min}}$ des couples de points fixes d'\'el\'ements hyperboliques
(rappelons que $\Ga$ est un sous-groupe discret non \'el\'ementaire de
$T$), il existe un \'el\'ement hyperbolique $\ga$ dans $\Ga$ ayant son
point fixe r\'epulsif ˆ distance au plus $\epsilon$ de $\xi$, et l'autre
point fixe dans $U$.  Soit $x$ un point de l'axe de translation de
$\ga$. Pour tout entier $n$ assez grand, $\ga^nx$ est fix\'e par
$\alpha$ (car celui-ci fixe $U$ point par point). Donc pour tout $n$
assez grand, $\ga^{-n}\alpha\ga^nx=x$.  Comme $\Ga$ est discret, il
existe des entiers distincts $n$ et $m$ tels que $\ga^{-n}\alpha\ga^n=
\ga^{-m}\alpha\ga^m$. Donc, en posant $p=n-m\neq 0$, l'\'el\'ement
elliptique $\alpha$ commute avec l'\'el\'ement hyperbolique $\ga^p$. Par
cons\'equent $\alpha$ fixe l'axe de translation de $\ga^p$, qui est
l'axe de translation de $\ga$, donc $\alpha$ admet un point fixe ˆ
distance au plus $\epsilon$ de $\xi$.  Donc l'ensemble des points
fixes de $\alpha$, qui est ferm\'e et dense dans $\partial T$, est \'egal
ˆ $\partial T$.  Par cons\'equent $\alpha$ vaut l'identit\'e sur $\partial
T$, donc sur $T$, ce qui contredit le fait que $\alpha$ soit non
trivial.

\medskip %
Remarquons toutefois que contrairement au cas des r\'eseaux
sans torsion, dont tout \'el\'ement diff\'erent de l'identit\'e n'a que deux
points fixes dans l'espace des bouts de l'arbre, il existe des r\'eseaux
ayant des \'el\'ements (elliptiques) non triviaux $\ga$ dont l'ensemble
des points fixes est un espace de Cantor, avec $\delta_\ga$
arbitrairement prescrit.

\brema %
Il existe une partie $\S$ dense dans $]0,+\infty[$, telle que pour
tout $s$ dans $\S$, il existe un triplet $(T,\Ga,\ga)$, avec $T$ un
arbre localement fini, $\Ga$ un r\'eseau uniforme de $T$ et $\ga$ un
\'el\'ement non trivial de $\Ga,$ tel que $\delta_\ga=s$.  \erema

\dem Soit $\S$ l'ensemble des $s>0$ tels qu'il existe un r\'eseau
uniforme d'un arbre localement fini dont l'entropie volumique est $s$.
Il est bien connu que $\S$ est une partie (d\'enombrable) dense de
$]0,+\infty[$ (par exemple, l'entropie volumique de la $p$-\`eme
subdivision barycentrique de l'arbre r\'egulier de valence $q+1$ est
$\frac{\log q}{p}$).

Soit $s$ un \'el\'ement de $\S$, et $\Ga_0$ un r\'eseau uniforme d'un arbre
localement fini $T_0$, d'entropie volumique $s$.  Notons $A_1,A_2$
deux copies du groupe $\ZZ/2\ZZ$. Consid\'erons le graphe de groupes
$\G$, obtenu ˆ partir du graphe de groupes $\Ga_0\bac T_0$ en y
remplacant chaque groupe d'ar\^ete ou de sommet par son produit direct
avec $A_1$, avec monomorphismes \'evidents, en rajoutant au graphe
$\Ga_0\backslash T_0$ une ar\^ete partant d'un sommet quelconque et
d'extr\'emit\'e libre, le groupe de la nouvelle ar\^ete \'etant trivial, et
le groupe du nouveau sommet \'etant $A_2$. Notons $T$ l'arbre de
Bass-Serre de $\G$, $\Ga$ le groupe fondamental de $\G$ pour un choix
indiff\'erent de point base dans $\Ga_0\backslash T_0$ (voir par exemple
\cite{Ser}), et $\ga$ l'\'el\'ement de $\Ga$ correspondant ˆ l'\'el\'ement non
trivial de $A_1$.  Alors $\Ga$ (isomorphe au produit libre
$(\Ga_0\times A_1)*A_2$) est un r\'eseau uniforme de $T$, $\ga$ est un
\'el\'ement elliptique non trivial de $\Ga$, l'ensemble ${\rm
  Fix}_{\partial T}(\ga)$ est un espace de Cantor dans $\partial T$
(qui est l'espace des bouts du sous-arbre ${\rm Fix}_T(\ga)$ de $T$,
isomorphe ˆ $T_0$, invariant non vide minimal par le groupe
fondamental du sous-graphe de groupes $\Ga_0\bac T_0$), et l'entropie
volumique $\delta_\ga$ de ${\rm Fix}_T(\ga)$ est \'egale au nombre
prescrit $s$.  \cqfd

\medskip %
Bien s\^ur, dans cet exemple, $\delta_\Ga>\delta_\ga$. Nous ne savons
pas s'il est possible de trouver des triplets $(T,\Ga,\ga)$ (disons en
imposant une borne fix\'ee sur les valences des sommets de $T$) tels que
la diff\'erence $\delta_\Ga-\delta_\ga$ soit arbitrairement petite.

\bigskip %
Revenons au cadre initial de cette partie. Dans toute la suite de
cette partie, nous supposons que le syst\`eme dynamique mesur\'e
$(\Ga\backslash\G' T,\psi,{\mu}_{\mbox{\tiny BM}})$ est de
probabilit\'e et m\'elangeant.

Par exemple, dans le cas o\`u $\Gamma$ est g\'eom\'etriquement fini et
$T_{{\Gamma,{\rm min}}} $ est uniforme, la proposition
\ref{prop:finitudebowmarg} assure que $\mu_{\mbox{\tiny
    BM}}$ est une mesure de probabilit\'e. La proposition
\ref{prop:melang} donne des conditions sur $\Gamma$ assurant que
$(\Ga\backslash\G' T,\psi,{\mu}_{\mbox{\tiny BM}})$ est m\'elangeant,
d'une part pour $N=1$ et $V'T=V T$, et d'autre part pour $N=2$ et
$V'T$ le sous-ensemble des sommets ˆ distance paire d'un sommet donn\'e
de $T$ .

D'abord, nous introduisons les notations qui vont permettre
d'\'enoncer le r\'esultat principal de cette partie.

Pour $x$ dans $V'T$, notons $\G'T(x)$ le sous-espace mesurable
$\Ga$-invariant et $\wt\psi$-invariant des \'el\'ements $\ell$ de
$\G'T$ tels que, pour une infinit\'e de temps positifs et de temps
n\'egatifs $t$, le point $\ell(t)$ appartienne \`a l'orbite $\Ga x$
(ou de mani\`ere \'equivalente par l'hypoth\`ese sur $V'T$, l'origine
de la g\'eod\'esique $\wt\psi\,^t\ell$ appartienne \`a $\Ga x$).  Par
passage au quotient et par ergodicit\'e de $(\Ga\backslash\G'
T,\psi,{\mu}_{\mbox{\tiny BM}})$, il existe un point $x_0$ dans $V'T$
tel que $\wt X=\G'T(x_0)$ soit de mesure pleine dans $\G'T$. On
appelle $X$ l'espace topologique quotient $\Ga\backslash\wt X$, et
$\pi:\wt X\ra X$ la projection canonique. On note de la m\^eme
mani\`ere les restrictions de $\wt\psi$ et $\wt\mu_{\mbox{\tiny BM}}$
\`a $\wt X$, ainsi que celles de $\psi$ et $\mu_{\mbox{\tiny BM}}$ \`a
$X$.

Soit $\wt X_0$ le sous-espace ferm\'e $\Ga$-invariant des
g\'eod\'esiques de $\wt X$ d'origine dans $\Ga x_0$. Pour tout $\ell$
dans $\wt X_0$, notons $t_\ell>0$ le premier temps de retour de
l'orbite de $\ell$ sous $\wt \psi$ dans $\wt X_0$~; par l'hypoth\`ese
sur $V'T$, c'est le minimum des entiers strictement positifs $t$ tels
que $\ell(t)$ appartienne \`a $\Ga x_0$.

Notons $\wt\psi_0:\wt X_0\ra\wt X_0$ l'application de premier retour
d\'efinie par $\wt\psi_0\ell=\psi\,^{t_\ell}\ell$, qui est
$\Ga$-\'equivariante.  Notons $X_0$ l'image de $\wt X_0$ par $\pi$.
L'application $\psi_0$ est l'application induite par passage au
quotient de $\wt\psi_{0}$ \`a $\wt X_{0}$. La mesure $m_{0}$ sur $X_0$
est la mesure induite par passage au quotient de la restriction $\wt
m_0$ de ${\wt \mu}_{\mbox{\tiny BM}}$ \`a $\wt X_0$. Renormalisons les
mesures de sorte que la mesure $m_0$ soit de probabilit\'e, ce qui est
possible car la mesure ${\mu}_{\mbox{\tiny BM}}$ est suppos\'ee \^etre
de probabilit\'e.

Notons $S$ l'ensemble des $\ga$ dans $\Ga-\Ga_{x_0}$ tels que
l'intersection $]\,x_0,\ga x_0[\;\cap\;\Ga {x_0}$ soit vide.
Remarquons que $S$ est invariant par translations \`a droite et \`a
gauche par $\Ga_{x_0}$.

\'Enon\c{c}ons maintenant le r\'esultat principal de cette partie.

\btheo \label{theo:codageseul} %
Il existe un sous-d\'ecalage $(X_A,\sigma)$ sur l'alphabet $S$ et une
mesure de Markov $\mu_\Pi$ sur ce sous-d\'ecalage, de sorte que le
syst\`eme dynamique pro\-ba\-bi\-li\-s\'e $(X_0, \psi_0,m_0)$ soit un
facteur du syst\`eme de Markov $(X_A,\sigma,\mu_\Pi)$.  
\etheo

Avant de d\'emontrer ce r\'esultat, nous allons \'enoncer et
d\'emontrer ses corollaires.

\bcoro \label{coro:lachebernou}%
Si le syst\`eme dynamique mesur\'e $(\Ga\backslash\G'
T,\psi,{\mu}_{\mbox{\tiny BM}})$ est de probabilit\'e et 
m\'elangeant, alors il est l\^achement Bernoulli.  
\ecoro

\dem %
Il suffit de montrer que le syst\`eme dynamique probabilis\'e
$(X,\psi,{\mu}_{\mbox{\tiny BM}})$, qui est de mesure pleine dans
$(\Ga\backslash\G' T,\psi,{\mu}_{\mbox{\tiny BM}})$, est l\^achement
Bernoulli.

Remarquons que $X_0$ est une transversale totale du syst\`eme
dynamique pro\-ba\-bi\-li\-s\'e $(X,\psi,{\mu}_{\mbox{\tiny BM}})$,
dont $\psi_0$ est l'application de premier retour, $\Ga\ell\mapsto
t_{\ell}$ le temps de premier retour, et $m_0$ la mesure induite, et
donc $(X,\psi,\mu_{\mbox{\tiny BM}})$ est une suspension m\'elangeante
de $(X_{0},\psi_{0},m_{0})$.

Rappelons qu'une transformation $T$ d\'efinie sur $X$ est {\it
  markovienne}, s'il existe une fonction $f$ d\'efinie sur $X$ et
prenant un nombre fini ou d\'enombrable de valeurs r\'eelles telle que
le processus $(f\circ T^n)_{n\in\NN}$ soit markovien et si la plus
petite tribu qui rende mesurable tous les it\'er\'es $f\circ T^n$ de
$f$ est la tribu tout enti\`ere de l'espace sur lequel est d\'efinie
$T$.

Avec le r\'esultat de Adler, Shields et Smorodinsky \cite{smorodinsky}
dans le cas d'un d\'ecalage de Markov sur un espace d'\'etats fini et
sa g\'en\'eralisation pour un espace d'\'etats d\'enombrable (il
suffit d'utiliser une g\'en\'eralisation du th\'eor\`eme de
Perron-Frobenius, voir le chapitre 7 du livre de Kitchens \cite{Kit}),
on peut caract\'eriser les transformations markoviennes~: une
transformation markovienne est le produit direct d'une rotation sur un
espace d'\'etats fini et d'un sch\'ema de Bernoulli.

Cette propri\'et\'e est stable par passage \`a un facteur, donc un
facteur d'une transformation markovienne est encore une transformation
markovienne.
 
Comme toute suspension m\'elangeante d'une transformation markovienne
est un syst\`eme l\^achement Bernoulli (voir par exemple \cite{Tho}),
le th\'eor\`eme \ref{theo:codageseul} implique bien le corollaire
\ref{coro:lachebernou}.  \cqfd

\medskip %
Ce r\'esultat permet de d\'eduire les deux \'enonc\'es suivants.  En
prenant $N=1$ et $V'T=VT$, nous obtenons l'\'enonc\'e du
th\'eor\`eme \ref{theo:introneuf} de l'introduction, que nous
rappelons ci-dessous.

\bcoro\label{coro:thmintroneuf} %
Soit $\Gamma$ un sous-groupe discret non \'el\'ementaire de ${\rm
  Aut}(T)$ et soit $\wt\mu_{\mbox{\tiny BM}}$ une mesure de
Bowen-Margulis pour $\Gamma$ sur $\G T$, dont la mesure de
Patterson-Sullivan sur $\partial T$ ne charge aucun ensemble de points
fixes d'\'el\'ement elliptique non trivial de $\Ga$. Si le syst\`eme
dynamique mesur\'e quotient de $(\G T,\wt\varphi,\wt\mu_{\mbox{\tiny
    BM}})$ par $\Gamma$ est de probabilit\'e et m\'elangeant, alors il
est l\^achement Bernoulli.  \cqfd 
\ecoro
 
\medskip %
Nous obtenons un second r\'esultat dans le cas alg\'ebrique.

\bcoro\label{coro:thmintro} %
Avec les notations $\underline{G},\wh K,G,S,M,\TT$ de l'introduction,
pour tout r\'e\-seau $\Gamma$ de $G$, si $\G_0 \TT$ est le sous-espace
de $\G \TT$ form\'e des g\'eod\'esiques d'origine un sommet ˆ distance
paire d'un sommet donn\'e de $\TT$, alors le syst\`eme dynamique
probabilis\'e $(\Gamma\backslash\G_{0}\TT,\psi,\mu_{\mbox{\tiny BM}})$
est l\^achement Bernoulli.

Ainsi, l'action par translation \`a droite de $S\backslash M$ sur
$\Gamma\backslash G\slash M$ est l\^achement Ber\-noulli.
\ecoro

\dem %
Nous allons utiliser le corollaire \ref{coro:lachebernou} avec $N=2$,
$V'T$ l'ensemble des sommets ˆ distance paire du sommet donn\'e de $\TT$
de sorte que $\G'T=\G_0\TT$, $\delta$ l'entropie volumique de $\TT$,
et $(\mu_x)_{x\in V\TT}$ la mesure de Patterson-Sullivan de dimension
$\delta$ (unique par \cite{BM} par exemple). Rappelons que tout r\'eseau
de $G$ est g\'eom\'etriquement fini et poss\`ede la propri\'et\'e de
Selberg. Soit $\Ga'$ un sous-groupe d'indice fini de $\Ga$ dont tout
\'el\'ement de torsion est conjugu\'e ˆ un \'el\'ement d'un stabilisateur d'un
sommet de $\TT$ se projetant dans un rayon cuspidal ouvert de
$\Ga'\bac\TT$. Alors par l'exemple suivant la remarque
\ref{rem:condnecmesfixellnul}, la mesure de Patterson-Sullivan sur
$\partial T$ ne charge aucun ensemble de points fixes d'\'el\'ement
elliptique non trivial de $\Ga$. Par la proposition \ref{prop:melang}
(1) appliqu\'ee ˆ $\TT$ o\`u l'on a enlev\'e les sommets de valences $2$
s'ils existent, ou la proposition \ref{prop:melang} (2) sinon, le
syst\`eme dynamique probabilis\'e
$(\Gamma'\backslash\G_{0}\TT,\psi,\mu_{\mbox{\tiny BM}})$ (voir la
convention suivant la preuve de la proposition
\ref{prop:finitudebowmarg}) est m\'elangeant.  Donc par le corollaire
\ref{coro:lachebernou}, le syst\`eme
$(\Gamma'\backslash\G_{0}\TT,\psi,\mu_{\mbox{\tiny BM}})$ est
l\^achement Bernoulli. Comme tout facteur d'un syst\`eme l\^achement
Bernoulli l'est encore (voir par exemple \cite{Tho}), la premi\`ere
assertion en d\'ecoule.

La seconde assertion d\'ecoule de la
premi\`ere par la correspondance rappel\'ee en introduction.  
\cqfd

\medskip
Ce r\'esultat est plus faible que le th\'eor\`eme \ref{theo:intro} de
l'introduction. Celui-ci sera compl\`etement d\'emontr\'e dans la partie
\ref{sec:codagegeomfin}.

\bigskip
{\it D\'emonstration du th\'eor\`eme \ref{theo:codageseul}~:}
D\'efinissons une matrice {\it de transition} $A=(A_{\alpha,\beta})
_{\alpha,\beta\in S}$.  Pour tous $\alpha,\beta$ dans $S$, posons
$A_{\alpha,\beta}=1$ si $x_0$ appartient au segment g\'eod\'esique entre
$\alpha^{-1}x_0$ et $\beta x_0$, et $A_{\alpha,\beta}=0$ sinon.
Remarquons que $A_{\alpha,\beta}=1$ si et seulement si la r\'eunion
$[x_0,\alpha x_0]\cup[\alpha x_0,\alpha\beta x_0]$ est un segment
g\'eod\'esique.

Notons $(X_A,\sigma)$ le sous-d\'ecalage d\'efini par la matrice de
transition $A$, i.e. $$X_A=\{(\alpha_i)_{i\in\mathbb{Z}}\in
S^{\mathbb{Z}}\;:\; \forall\; i\in{\mathbb{Z}},\;
A_{\alpha_i,\alpha_{i+1}}=1\}\;,$$
muni de la restriction de la topologie produit, et $\sigma$ la
restriction \`a $X_A$ du d\'ecalage vers la gauche des suites bilat\`eres de
$S^{\mathbb{Z}}$.

Construisons une application $\Theta:X_A\ra X_0$ de la mani\`ere
suivante. Soit $(\alpha_i)_{i\in\mathbb{Z}}\in X_A$, d\'efinissons une
suite $(\ga_i)_{i\in\mathbb{Z}}$ dans $\Ga$ par r\'ecurrence, en posant
$\ga_0$ l'\'el\'ement neutre $e$ de $\Ga$, et $\ga_{i+1}=
\ga_i\alpha_{i+1}$ pour tout $i$ dans $\mathbb{Z}$. Alors les points
$(\ga_i x_0)_{i\in\mathbb{Z}}$ sont cons\'ecutivement align\'es sur une
g\'eod\'esique $\ell$. En effet, comme $A_{\alpha_i,\alpha_{i+1}}=1$, la
r\'eunion $[x_0,\alpha_i x_0]\cup[\alpha_i x_0,\alpha_i\alpha_{i+1}
x_0]$ est un segment g\'eod\'esique, donc son image par $\ga_{i-1}$ aussi
et donc $\ga_{i-1} x_0,\ga_{i} x_0,\ga_{i+1} x_0$ sont bien align\'es
dans cet ordre.  Param\'etrons $\ell$ de sorte que $\ell(0)=x_0$ et que
$\ga_i x_0$ converge vers les extr\'emit\'es $\ell_\pm$ de $\ell$ quand
$i\ra\pm\infty$. Par d\'efinition de $\wt X_0$, la g\'eod\'esique $\ell$
appartient \`a $\wt X_0$.  Posons alors $\Theta((\alpha_i)
_{i\in\mathbb{Z}}) =\pi(\ell)$.

Montrons que $\Theta$ est surjective. En effet, soit $\ov\ell$ un
\'el\'ement de $X_0$. Choisissons un relev\'e $\ell$ de $\ov\ell$ dans $\wt
X_0$ tel que $\ell(0)=x_0$.  Notons $(\ga_i x_0)_{i\in\mathbb{Z}}$ la suite
des points cons\'ecutifs de $\ell$ dans $\Ga{x_0}$, o\`u l'on peut
supposer que $\ga_0=e$.  Posons $\alpha_{i}=\ga_{i-1}^{-1}\ga_{i}$,
qui appartient \`a $S$. Puisque $\ga_{i-1} x_0,\ga_{i} x_0,\ga_{i+1}
x_0$ sont align\'es dans cet ordre sur $\ell$, nous avons
$A_{\alpha_i,\alpha_{i+1}}=1$. Par construction, nous avons alors
$\Theta((\alpha_i) _{i\in\mathbb{Z}}) =\pi(\ell)=\ov\ell$.

Par d\'efinition de la topologie produit sur $S^{\mathbb{Z}}$ et de
la topologie de $\G T$, l'application $\Theta:X_A\ra  X_0$ est
continue. Il est imm\'ediat que le diagramme suivant est commutatif~:
$$
\begin{array}{ccc}
X_A &\stackrel{\Theta}{\longrightarrow} & X_0 \smallskip\\
\sigma\downarrow\;\;\; & & \;\;\downarrow\psi_0\smallskip\\
X_A &\stackrel{\Theta}{\longrightarrow} &  X_0
\end{array}\;.  \;\;\;(*)
$$

Donc topologiquement, le syst\`eme dynamique $(X_0, \psi_0)$ est un
facteur du syst\`eme dynamique symbolique $(X_A,\sigma)$.

\medskip %
Nous allons maintenant construire une mesure de Markov sur le
sous-d\'ecalage $X_A$, dont l'image par $\Theta$ sera $m_0$. Nous
commen\c cons pour cela par d\'efinir une mesure $(\nu_\alpha)_{\alpha\in
  S}$ sur l'ensemble d\'enombrable discret $S$ et des probabilit\'es de
transitions $(\pi_{\alpha,\beta})_{\alpha,\beta\in S}$, apr\`es quelques
notations.

\medskip %
Pour tout $\alpha$ dans $S$, notons $e_\alpha^+$ (respectivement
$e_\alpha^-$) l'ar\^ete (orient\'ee) d'origine $\alpha x_0$
(respectivement $x_0$) et contenue dans $[x_0,\alpha x_0]$. Remarquons
que pour $\alpha$ dans $\Gamma_{x_{0}}$ et $\beta$ dans $S$, nous
avons $\alpha \,e_\beta^+=e_{\alpha\beta}^+$. Pour $\alpha,\beta$ dans
$S$, nous avons $A_{\alpha,\beta}=1$ si et seulement si
$\alpha\,\partial_{e_\beta^+}T= \partial_{\alpha e_\beta^+}T$ est
contenu dans $\partial_{e_\alpha^+} T$.

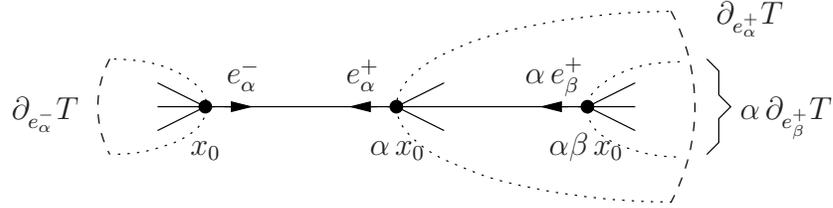
\medskip
\begin{center}
\input{fig_transition.pstex_t}
\end{center}
\medskip

Notons $\partial_0 T$ l'ensemble des points $\xi$ de $\partial T$ tels
qu'il existe une infinit\'e de points de l'orbite $\Ga x_0$ sur le rayon
g\'eod\'esique entre $x_0$ et $\xi$. Il est invariant par $\Ga$. 

Comme la g\'eod\'esique entre deux points de $\partial T$ est contenue
dans la r\'eunion des rayons g\'eod\'esiques entre $x_0$ et ces points, une
g\'eod\'esique $\ell$ de $T$ appartient \`a $\G'T(x_0)$ si et seulement si
ses extr\'emit\'es $\ell_\pm$ appartiennent \`a $\partial_0 T$. En
particulier, $\wt X=\{\ell\in\G' T\;:\; \ell_{\pm}\in \partial_0 T\}$
et $\wt X_0=\{\ell\in\G T\;:\; \ell(0)\in\Ga x_0,\ell_{\pm}\in
\partial_0 T\}$. Par le param\'etrage de Hopf, l'image de $\wt X$ est
donc contenue dans $\partial_0 T\times\partial_0T\times\mathbb{Z}$.
Comme $\wt X$ est de mesure pleine dans $\G' T$ pour la mesure de
Bowen-Margulis, et par les propri\'et\'es de celle-ci, nous en
d\'eduisons que $\partial_0 T$ est de mesure pleine dans $\partial T$
pour la mesure de Patterson-Sullivan $\mu_{x_0}$.

Comme $\Ga$ est non \'el\'ementaire, quitte \`a remplacer $T$ par son
unique sous-arbre non vide invariant minimal, nous pouvons supposer
que $T$ soit sans ar\^ete terminale, et que le support de la mesure de
Patterson-Sullivan soit \'egal \`a $\partial T$. En particulier, pour
tout sommet $x$ et ar\^ete $e$ de $T$, on a $\mu_{x}(\partial_eT)>0$.
Rappelons que le groupe $\Ga_{x_0}$ est de cardinal $|\Ga_{x_0}|$
fini.

\medskip %
Pour tout $\alpha$ dans $S$, posons maintenant
$$\nu_{\alpha}=
\frac{1}{|\Ga_{x_0}|^2}\;\;\mu_{x_0}(\partial_{e_\alpha^-}T)
\;\mu_{x_0}(\partial_{e_\alpha^+}T)>0\;.$$
Pour tous $\alpha,\beta$ dans $S$, posons $\pi_{\alpha,\beta} =0$ si 
$A_{\alpha,\beta} =0$, et sinon
$$\pi_{\alpha,\beta} =\frac{1}{|\Ga_{x_0}|}\;
\frac{\mu_{x_0}(\partial_{\alpha\,e_\beta^+}T)}
{\mu_{x_0}(\partial_{e_\alpha^+}T)}>0\;.$$

\medskip %
\noindent{\bf Remarques} %
 (1) Pour $\alpha$ dans $S$, les parties $\partial_{e_\alpha^-} T,
\partial_{e_\alpha^+} T$ ne d\'ependent que de la classe \`a droite de
$\alpha$ modulo $\Ga_{x_0}$, et pour $\alpha,\beta$ dans $S$,
l'\'egalit\'e $A_{\alpha,\beta}=1$ ne d\'epend que de la classe \`a
gauche de $\alpha$ modulo $\Ga_{x_0}$, et de la classe \`a droite de
$\beta$ modulo $\Ga_{x_0}$.

(2) Pour tous $\ga,\ga',\ga''$ dans $\Ga_{x_0}$ et tous
$\alpha,\beta$ dans $S$, par invariance de la mesure $\mu_{x_0}$
par $\Ga_{x_0}$, on a
$$\nu_{\alpha}=\nu_{\ga\alpha\ga'}\;\;\;{\rm et}
\;\;\;\pi_{\alpha,\beta} =
\pi_{\ga'\alpha\ga,\ga^{-1}\beta\ga''}\;.\;\;\;(\sharp)$$

\blemm \label{lem:constructmatstochmesstat} %
La matrice $\Pi=(\pi_{\alpha,\beta})_{\alpha,\beta\in S}$ est une
matrice stochastique sur $S$, de mesure stationnaire
$\nu=(\nu_\alpha)_{\alpha\in S}$.
\elemm

\dem %
Il s'agit de montrer que $\sum_{\alpha\in S}\nu_{\alpha}= 1$, que
$\sum_{\beta\in S}\pi_{\alpha,\beta}= 1$ pour tout $\alpha$ dans $S$,
et que $\sum_{\alpha\in S}\nu_\alpha\pi_{\alpha,\beta}= \nu_\beta$ pour
tout $\beta$ dans $S$. Nous commen\c cons la preuve par quelques
notations et r\'esultats ensemblistes.

Pour tout $x$ dans $V'T$ et toutes les parties $B_-,B_+$ de $\partial
T$ tels que toutes les g\'eod\'esiques entre un point de $B_-$ et un point
de $B_+$ passent par $x$, notons $B_+\times B_-\times\{x\}$ l'ensemble
des g\'eod\'esiques $\ell$ de $T$ telles que $\ell(0)=x,\ell_-\in
B_-,\ell_+\in B_+$. Par d\'efinition de la mesure de Bowen-Margulis,
remarquons que
$${\wt \mu}_{\mbox{\tiny BM}}(B_-\times B_+\times\{x\})=
\mu_{x}(B_-)\;\mu_{x}(B_+)\;.$$
En particulier, pour tout $\alpha$ dans $S$, on a
$$\nu_{\alpha}=\frac{1}{|\Ga_{x_0}|^2}\;\; 
\wt\mu_{\mbox{\tiny BM}}(\partial_{e_\alpha^-}T
\times\partial_{e_\alpha^+}T\times\{x_{0}\})\;.$$

\blemm \label{lem:sommdisj} %
Nous avons les r\'eunions disjointes suivantes (o\`u nous notons de la
m\^eme mani\`ere une classe et l'un de ses repr\'esentants)~:
\begin{enumerate}
\item[(1)] $\displaystyle \partial_0 T=\bigsqcup_{\alpha\in S/\Ga_{x_0}}
  (\partial_{e_\alpha^+} T\cap \partial_0 T)$,
\item [(2)] pour tout $\alpha $ dans
  $S$,
$\displaystyle \partial_{e_\alpha^+} T\cap \partial_0
  T=\bigsqcup_{\beta\in S/\Ga_{x_0}\;:\;A_{\alpha,\beta}=1} (\partial_{\alpha
    e_\beta^+} T\cap \partial_0 T)$,
\item[(3)]  pour tous
  $\beta$ dans $S$,
$\displaystyle \partial_{e_\beta^-} T\cap \partial_0
  T=\bigsqcup_{\alpha\in \;\Ga_{x_0}\backslash S\;:\;A_{\alpha,\beta}=1}
  (\partial_{\alpha^{-1} e_\alpha^-} T\cap \partial_0 T)$,
  \item[(4)] $\displaystyle X_0=\bigsqcup_{\alpha\in\;
      \Ga_{x_0}\backslash S/\Ga_{x_0}}\pi\left(
  (\partial_{e_\alpha^-} T\cap \partial_0 T)\times
  (\partial_{e_\alpha^+} T\cap \partial_0 T)\times\{x_0\}\right)$.
\end{enumerate}
\elemm

\dem %
Les objets sont bien d\'efinis, par la premi\`ere des remarques pr\'ec\'edant
l'\'enonc\'e du lemme \ref{lem:constructmatstochmesstat}.

Montrons (1). Soient $\alpha,\beta$ dans $S$, et supposons qu'il
existe $\xi$ dans $\partial_{e_\alpha^+} T\cap \partial_{e_\beta^+}
T\cap \partial_0 T$, alors le rayon g\'eod\'esique de $x_0$ \`a $\xi$ passe
par $\alpha x_0$ et par $\beta x_0$. De plus, ces points co\"{\i}ncident,
car ce sont tous les deux le premier point de $\Ga x_0-\{x_0\}$
rencontr\'e par le rayon.  Donc $\alpha \Ga_{x_0}=\beta \Ga_{x_0}$. Par
cons\'equent, la r\'eunion dans l'assertion (1) est disjointe.
 
Pour tout $\xi$ dans $\partial_0 T$, le rayon g\'eod\'esique de $x_0$ \`a
$\xi$ passe par un premier point de l'orbite $\Ga{x_0}$ diff\'erent de
$x_0$. Ce point s'\'ecrit $\alpha x_0$ avec $\alpha$ dans $S$, donc
$\xi$ appartient \`a $\partial_{e_\alpha^+} T$, ce qui montre (1).

Les autres assertions se montrent de mani\`ere analogue.  \cqfd

\medskip %
Par le lemme \ref{lem:sommdisj} (4), et le fait que l'action de
$\Ga_{x_0}$ sur $(\partial T,\mu_{x_{0}})$ soit essentiellement libre,
nous avons
\begin{eqnarray*}
\sum_{\alpha\in S}\nu_{\alpha} 
&=&
\frac{1}{|\Ga_{x_0}|^2}\;\sum_{\alpha\in S}\;
\wt \mu_{\mbox{\tiny  BM}}(\partial_{e_\alpha^-}T
\times\partial_{e_\alpha^+}T\times\{x_{0}\})\\
&=&
\frac{1}{|\Ga_{x_0}|}\;
\sum_{\alpha\in S/\Ga_{x_0}} \;
\wt\mu_{\mbox{\tiny BM}}(\partial_{e_\alpha^-}T
\times\partial_{e_\alpha^+}T\times\{x_{0}\})\\
&=&
\sum_{\alpha\in\;\Ga_{x_0}\!\backslash S/\Ga_{x_0}}
\frac{1}{|\Ga_{[x_0,\alpha x_0]}|}\;\wt
\mu_{\mbox{\tiny BM}}(\partial_{e_\alpha^-}T
\times\partial_{e_\alpha^+}T\times\{x_{0}\})\\
&=&
\sum_{\alpha\in\;\Ga_{x_0}\!\backslash S/\Ga_{x_0}}\!\!
m_0(\pi(\partial_{e_\alpha^-}T
\times\partial_{e_\alpha^+}T\times\{x_{0}\}))=
 {m}_0(X_0)=1\;.
\end{eqnarray*}

Par le lemme \ref{lem:sommdisj} (2), nous avons $\sum_{\beta\in
  S}\pi_{\alpha,\beta}= 1$.  De plus, si $A_{\alpha,\beta}=1$, puisque
toute g\'eod\'esique entre un point de $\partial_{\,e_\alpha^-}T$ et
un point de $\partial_{\alpha\,e_\beta^+}T$ passe par $\alpha x_0$, et
par invariance de $\wt\mu_{\mbox{\tiny BM}}$, alors

\begin{eqnarray*}
\nu_\alpha\pi_{\alpha,\beta}
&=&
\frac{1}{|\Ga_{x_0}|^3}\;
\mu_{x_0}(\partial_{\,e_\alpha^-}T)
\mu_{x_0}(\partial_{\alpha\,e_\beta^+}T)=\frac{1}{|\Ga_{x_0}|^3}\;
\wt\mu_{\mbox{\tiny    BM}}(\partial_{\,e_\alpha^-}T
\times\partial_{\alpha\,e_\beta^+}T\times\{\alpha x_0\})\\
&=&\frac{1}{|\Ga_{x_0}|^3}\;\wt\mu_{\mbox{\tiny
    BM}}(\partial_{\alpha^{-1}\,\,e_\alpha^-}T
\times\partial_{e_\beta^+}T\times\{x_0\})\;.
\end{eqnarray*}
 Donc 
$$\sum_{\alpha\in S} \nu_\alpha\pi_{\alpha,\beta}= \sum_{\alpha\in
  \;\Ga_{x_0}\backslash S\;:\;A_{\alpha,\beta}=1} \;
\frac{1}{|\Ga_{x_0}|^2} \;\wt\mu_{\mbox{\tiny
    BM}}(\partial_{\alpha^{-1}\,\,e_\alpha^-}T
\times\partial_{e_\beta^+}T\times\{x_0\})=\nu_\beta$$ 
par le lemme \ref{lem:sommdisj} (3). Ceci conclut la preuve du lemme
\ref{lem:constructmatstochmesstat}. \cqfd

\bigskip %
Par d\'efinition, la {\it mesure de Markov} $\mu_\Pi$ associ\'ee \`a
$(\Pi,\nu)$ est l'unique mesure bor\'elienne de probabilit\'e sur
$S^{\mathbb{Z}}$, invariante par le d\'ecalage, telle que, pour tout
$k$ dans $\NN$ et tous $\alpha_0,\dots, \alpha_k$ dans $S$, on ait
$$\mu_\Pi[X_0=\alpha_0,\dots ,X_k=\alpha_k]=
\nu_{\alpha_0}\left(\prod_{i=0}^{k-1} \pi_{\alpha_i,
    \alpha_{i+1}}\right) \;.$$

\blemm \label{lem:egalmesusi}
Les mesures  $\Theta_*\mu_\Pi$ et $m_0$ co\"{\i}ncident.
\elemm

\dem %
Pour $\alpha_1,\dots,\alpha_{n}$ dans $S$ tels que pour tout $i$ dans
$\{1,\dots, n-1\}$, nous ayons $A_{\alpha_{i},\alpha_{i+1}}=1$, notons
$U_{\alpha_1,\dots,\alpha_{n}}$ l'image par $\pi$ de l'ensemble $\wt
U_{\alpha_1,\dots,\alpha_{n}}$ des $\ell$ dans $\wt X_0$ tels que les
$n+1$ premiers points de $\Ga x_0$ rencontr\'es par $\ell$ \`a partir
de l'instant $0$ soient $x_0, \alpha_1 x_0,\dots,
\alpha_1\dots\alpha_{n} x_0$.  Les $U_{\alpha_1,\dots,\alpha_{n}}$ et
leurs images par les puissances de $\psi_0$ sont des bor\'eliens de
$X_{0}$, qui engendrent la tribu des bor\'eliens de $X_0$.

\bigskip
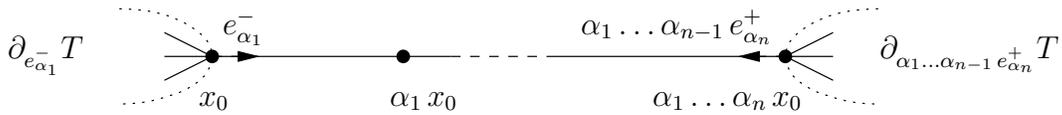
\begin{center}
$\!\!\!\!\!\!\!\!\!\!\!\!\!\!\!\!\!\!\!\!\!\!\!\!\!\!\!\!\!\!\!\!\!$
\input{fig_markovienbis.pstex_t}
\end{center}
\medskip

Par commutativit\'e du diagramme $(*)$, les mesures $\Theta_*\mu_\Pi$ et
$m_0$ sont deux mesures de probabilit\'e sur $X_0$ invariantes par
$\psi_0$. Pour montrer que $\Theta_*\mu_\Pi=m_0$, il suffit donc de
montrer que
$$\Theta_*\mu_\Pi(U_{\alpha_1,\dots,\alpha_{n}})=
m_0(U_{\alpha_1,\dots,\alpha_{n}})$$
pour tous les $\alpha_1,\dots,\alpha_{n}$ dans $S$ tels que
$A_{\alpha_{i},\alpha_{i+1}}=1$. Proc\'edons par r\'ecurrence sur $n$, le
cas $n=0$ \'etant clair, car, par les conventions usuelles, 
$\Theta_*\mu_\Pi(X_0)=m_0(X_0)=1$.

Rappelons que, par d\'efinition de $\Theta$, pour tout
$(Z_i)_{i\in\mathbb{Z}}$ dans $X_A$, on a
$\Theta((Z_i)_{i\in\mathbb{Z}})=\pi(\ell')$ o\`u les $n+1$ premiers
points de $\Ga x_0$ rencontr\'es par $\ell'$ \`a partir de l'instant
$0$ sont $x_0$, $Z_1 x_0$, $\dots$, $Z_1\dots Z_{n} x_0$. Remarquons
que $\pi(\ell')$ appartient \`a $\pi(\wt
U_{\alpha_1,\dots,\alpha_{n}})$ si et seulement s'il existe
$\beta_0,\beta_1,\dots,\beta_n$ dans $\Ga_{x_0}$ tels que
$$\forall\;i\in \{1,\dots,n\},\; \beta_0^{-1}\alpha_1\dots
\alpha_i\beta_i=Z_1\dots Z_{i}\;. \;\;\;(*\, *)$$
Donc $\Theta((Z_i)_{i\in\mathbb{Z}})$ appartient \`a
$U_{\alpha_1,\dots,\alpha_{n}}$ si et seulement s'il existe
$\beta_0,\beta_1,\dots,\beta_n$ dans $\Ga_{x_0}$ tels que
$$\forall\;i\in \{1,\dots,n\},\;\;Z_i=
\beta_{i-1}^{-1}\alpha_i\beta_i\;.\;\;\;(***)$$


Par d\'efinition de la mesure de Markov et les propri\'et\'es
d'invariance $(\sharp)$ des $\nu_\alpha,\pi_{\alpha,\beta}$, nous
avons, pour tous $\beta_0,\beta_1,\dots,\beta_n$ dans $\Ga_{x_0}$,
$$\mu_\Pi[Z_1=\beta_0^{-1}\alpha_1\beta_1,\dots,
Z_{n}=\beta_{n-1}^{-1}\alpha_n\beta_n]
=\nu_{\alpha_1}\left(\prod_{i=1}^{n-1}
  \pi_{\alpha_i, \alpha_{i+1}}\right) \;.$$

D\'efinissons une relation d'\'equivalence $\sim$ sur
$\Ga_{x_0}^{n+1}$ par $$(\beta_0,\beta_1,\dots,\beta_n)\sim
(\beta'_0,\beta'_1,\dots,\beta'_n)\;\;\Leftrightarrow\;\; \left(
  \forall\;i\in \{1,\dots,n\},\;\;\beta_{i-1}^{-1}\alpha_i\beta_i=
  {\beta'}_{i-1}^{-1}\alpha_i\beta'_i\right)\;.$$ Le cardinal de
chaque classe d'\'equivalence est \'egal \`a celui du stabilisateur
$\Ga_{[x_0,\alpha_1\dots\alpha_n x_0]}$ (car, par l'\'equivalence
entre les assertions $(*\,*)$ et $(***)$, dans une classe donn\'ee,
$\beta_0$ est uniquement d\'etermin\'e modulo ce stabilisateur, et les
$\beta_i$ pour $i>0$ sont uniquement d\'etermin\'es par $\beta_0$).
Donc le nombre $N_{\alpha_1,\dots,\alpha_{n}}$ de classes
d'\'equivalence vaut
$$N_{\alpha_1,\dots,\alpha_{n}}=
\frac{|\Ga_{x_0}|^{n+1}}{|\Ga_{[x_0,\alpha_1\dots\alpha_n x_0]}|}\;.
\;\;\;(***\,*)$$

Par les propri\'et\'es d'invariance des mesures de Patterson-Sullivan
$(\mu_x)_{x\in T}$, pour toute partie $A$ de $\partial T$ et tous les
points $x,y$ de $T$, si le rayon g\'eod\'esique entre $x$ et chaque point
de $A$ passe par $y$, alors $\mu_x(A)=e^{-\delta d(x,y)}\mu_y(A)$. En
particulier, pour $1\leq i\leq n$,
$$\mu_{x_0}(\partial_{e^+_{\alpha_i}} T)=
\mu_{\alpha_1\dots\alpha_{i-1} x_0} 
(\alpha_1\dots\alpha_{i-1}\partial_{e^+_{\alpha_i}} T)=
e^{-\delta d(x_0,\,\alpha_1\dots\alpha_{i-1}x_0)}\;
\mu_{x_0}(\alpha_1\dots\alpha_{i-1}\partial_{e^+_{\alpha_i}} T)\;.$$

 Il vient alors
\noindent\begin{eqnarray*}
 &&  \mu_\Pi(\Theta^{-1}(U_{\alpha_1,\dots,\alpha_{n}}))
  =
  N_{\alpha_1,\dots,\alpha_{n}}\;\nu_{\alpha_1}\left(\prod_{i=1}^{n-1}
    \pi_{\alpha_i, \alpha_{i+1}}\right)\\
  &\;&  \quad\quad\quad\quad\quad =  
\frac{N_{\alpha_1,\dots,\alpha_{n}}}{|\Ga_{x_0}|^{n+1}}
  \mu_{x_0}(\partial_{e^-_{\alpha_1}} T)
  \mu_{x_0}(\partial_{e^+_{\alpha_1}} T)
  \left(\prod_{i=1}^{n-1}\;\frac{
      \mu_{x_0}(\partial_{\alpha_i e^+_{\alpha_{i+1}}} T)
    }{\mu_{x_0}(\partial_{e^+_{\alpha_i}} T)}\right)\\
  &\;&  \quad\quad\quad\quad\quad=
  \frac{N_{\alpha_1,\dots,\alpha_{n}}}{|\Ga_{x_0}|^{n+1}}
  \mu_{x_0}(\partial_{e^-_{\alpha_1}} T)
  \mu_{x_0}(\partial_{e^+_{\alpha_1}} T)
  \left(\prod_{i=1}^{n-1}\;\frac{
      \mu_{x_0}(\alpha_1\dots\alpha_{i-1}\alpha_i 
      \;\partial_{e^+_{\alpha_{i+1}}} T)
    }{\mu_{x_0}(\alpha_1\dots\alpha_{i-1}\;\partial_{e^+_{\alpha_i}}
      T)}\right)\\
  &\;&  \quad\quad\quad\quad\quad= 
  \frac{N_{\alpha_1,\dots,\alpha_{n}}}{|\Ga_{x_0}|^{n+1}} \;
  \mu_{x_0}(\partial_{e^-_{\alpha_1}} T)
  \mu_{x_0}(\partial_{\alpha_1\dots\alpha_{n-1}
    \;e^+_{\alpha_{n}}} T)\\
 &\;&   \quad\quad\quad\quad\quad =
  \frac{N_{\alpha_1,\dots,\alpha_{n}}}{|\Ga_{x_0}|^{n+1}} \;
  \wt\mu_{\mbox{\tiny
      BM}}(\partial_{e^-_{\alpha_1}} T\times
\partial_{\alpha_1\dots\alpha_{n-1}
 \;e^+_{\alpha_{n}}} T\times \{x_0\})\;. 
 \end{eqnarray*}

 Comme $\Ga_{x_0}$ agit essentiellement librement sur $\partial T$, et
 que tout \'el\'ement de $U_{\alpha_1,\dots,\alpha_{n}}$ est l'image
 d'exactement $|\Ga_{[x_0,\alpha_1\dots\alpha_n x_0]}|$ \'el\'ements
 de $\partial_{e^-_{\alpha_1}} T\times
 \partial_{\alpha_1\dots\alpha_{n-1} \; e^+_{\alpha_{n}}} T\times\{x_0\}$, on~a
$$ m_0(U_{\alpha_1,\dots,\alpha_{n}})=
\frac{1}{|\Ga_{[x_0,\alpha_1\dots\alpha_n x_0]}|}\;
\wt\mu_{\mbox{\tiny
    BM}}(\partial_{e^-_{\alpha_1}} T\times
\partial_{\alpha_1\dots\alpha_{n-1}
  \;e^+_{\alpha_{n}}} T\times \{x_0\})\;.
$$

Par $(***\,*)$, on a donc $\mu_\Pi(\Theta^{-1}
(U_{\alpha_1,\dots,\alpha_{n}}))= m_0(U_{\alpha_1,\dots,\alpha_{n}})$,
ce qu'il fallait d\'emontrer.
\cqfd

\bigskip %
Le lemme \ref{lem:egalmesusi} ach\`eve de montrer que $\Theta$ est un
morphisme surjectif entre les syst\`emes dynamiques
pro\-ba\-bi\-li\-s\'es $(X_A,\sigma,\mu_\Pi)$ et $(X_0,\psi_0,m_0)$.
Ceci conclut la preuve du th\'eor\`eme \ref{theo:codageseul}.  \cqfd

\section{Codage du flot g\'eod\'esique sur un 
graphe de grou\-pes g\'eom\'etriquement fini}
\label{sec:codagegeomfin}

Le but de cette partie est de montrer le r\'esultat fin suivant,
am\'eliorant pour des groupes g\'eom\'etriquement finis le corollaire
\ref{coro:lachebernou}.

\btheo \label{theo:geodBernou} Soit $T$ un arbre localement fini,
$\Ga$ un sous-groupe g\'eom\'etriquement fini de ${\rm Aut}(T)$ ayant la
propri\'et\'e de Selberg, avec $T_{\Ga,{\rm min}}$ uniforme, sans sommet
de valence $2$. 

Si la transformation g\'eod\'esique sur $\Ga\backslash \G
T$ est m\'elangeante pour la mesure de Bowen-Margulis, alors elle est
Bernoulli d'entropie finie pour cette mesure.  

Si la transformation g\'eod\'esique sur $\Ga\backslash \G T$ n'est pas
m\'elangeante pour la mesure de Bowen-Margulis, alors son carr\'e est
Bernoulli d'entropie finie en restriction \`a $\Ga\backslash \G_0 T$,
o\`u $\G_0 T$ est le sous-espace de $\G T$ form\'e des g\'eod\'esiques
d'origine un sommet ˆ distance paire d'un sommet donn\'e de $T$.
 \etheo

\dem %
Nous pouvons supposer que $T=T_{\Ga,{\rm min}}$. 

Notons $X=\Ga\backslash T$ le graphe quotient, et $X_0$ le sous-graphe
compl\'ementaire des rayons cuspidaux maximaux ouverts.  Pour tout
$\beta$ dans $EX_0$, fixons un relev\'e $\wt\beta$ de $\beta$ dans
$ET$, de sorte que $\wt{\ov\beta}=\overline{\wt\beta}$.  Num\'erotons
de $1$ \`a $r$ les rayons cuspidaux maximaux. Notons
$(a_{i,n})_{n\in\NN}$ la suite cons\'ecutive des ar\^etes, orient\'ees
vers le bout, du $i$-\`eme rayon cuspidal maximal. On fixe un relev\'e
dans $T$ de ce rayon, dont on note $(\wt{a_{i,n}})_{n\in\NN}$ la suite
cons\'ecutive des ar\^etes, orient\'ees vers un point \`a l'infini
not\'e $\xi_i$ de $T$.  Notons $\Ga_{i,n}$ le stabilisateur du sommet
$o(\wt{a_{i,n}})$, de sorte que $\Ga_{i,n}$ est contenu dans
$\Ga_{i,n+1}$, et que $\Ga_{\xi_i}=\bigcup_{n\in\NN} \Ga_{i,n}$.

Rappelons qu'un facteur d'une transformation m\'elangeante
(resp.~Bernoulli d'en\-tro\-pie finie) l'est encore (voir par exemple
\cite{Orn}). Comme $\Ga$ poss\`ede la propri\'et\'e de Selberg, il
existe un sous-groupe d'indice fini $\Ga'$ de $\Ga$ dont tout
stabilisateur de sommet ne se projetant pas dans un rayon cuspidal
ouvert est trivial. Supposons le r\'esultat d\'emontr\'e pour $\Ga'$.
Si la transformation g\'eod\'esique de $\Ga'\backslash \G T$ est
m\'elangeante, alors elle est Bernoulli d'entropie finie, et son
facteur $\Ga\backslash \G T$ est aussi m\'elangeant et Bernoulli
d'entropie finie.  Sinon, le carr\'e de la transformation
g\'eod\'esique de $\Ga'\backslash \G_0 T$ est Bernoulli d'entropie
finie. Si $\G_0 T$ est invariant par $\Ga$, alors le facteur
$\Ga\backslash \G_0 T$ de $\Ga'\backslash \G_0 T$ est aussi Bernoulli
d'entropie finie pour le carr\'e de la transformation g\'eod\'esique,
et la transformation g\'eod\'esique de $\Ga\backslash \G T$ n'est pas
m\'elangeante. Si $\G_0 T$ n'est pas invariant par $\Ga$, alors
$\Ga\backslash \G T$ est un facteur de $\Ga'\backslash \G_0 T$, donc
le carr\'e de la transformation g\'eod\'esique de $\Ga\backslash \G T$
est Bernoulli d'entropie finie, et donc la transformation g\'eod\'esique de
$\Ga\backslash \G T$ est Bernoulli d'entropie finie (donc m\'elangeante).

\medskip 
Dans la suite, nous supposons donc que le groupe $G_x$ est
trivial pour tout sommet $x$ dans $X_0$.

Par exemple, le graphe de groupes suivant est isomorphe au graphe de
groupes quotient de l'arbre de Bruhat-Tits de ${\rm SL}_2$ sur le
corps $\FF_2((X^{-1})))$ par le r\'eseau (de congruence) $\Ga={\rm
  ker}({\rm SL}(2,\FF_2[X])\ra {\rm SL}(2,\FF_2))$, avec
$\Ga_i=\{P\in\FF_2[X]\;:\; P(0)=0, {\rm deg}\;P=i\}$.

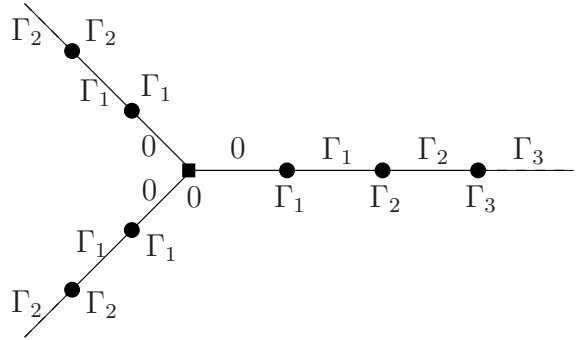
\begin{center}
\input{fig_multicusp.pstex_t}
\end{center}

En ce qui concerne le probl\`eme du codage du flot g\'eod\'esique, la
nouveaut\'e par rapport au cas cocompact est de coder les excursions
dans les rayons cuspidaux.  La remarque fondamentale (voir \cite{Ser})
est que la projection dans $\Ga\backslash T$ d'une g\'eod\'esique de $T$
est un chemin d'ar\^etes, dont les seuls aller-retours possibles sont
dans les rayons cuspidaux maximaux ouverts, et qui, s'il fait
demi-tour en montant vers un bout, est alors oblig\'e de redescendre
pour sortir compl\`etement du rayon cuspidal maximal.

\medskip 
Consid\'erons l'alphabet suivant, qui est d\'enombrable infini, sauf si
$\Ga$ est cocompact (auquel cas il est fini),
$$\A=EX\;\;\cup\;\bigcup_{1\leq i\leq r\,,\;n\in\NN}
\left((\Ga_{i,n+1}-\Ga_{i,n})/\Ga_{i,n}
  ^{\mbox{}^{\mbox{}}}\times\{+\}\right)\;\cup\;
\left(\Ga_{i,n}^{\mbox{}^{\mbox{}}}\backslash(\Ga_{i,n+1}-\Ga_{i,n})
  \times\{-\}\right)\,.$$ Nous munissons $\A$ de la topologie
discr\`ete. Notons que le passage \`a l'inverse induit une bijection
de $(\Ga_{i,n+1}-\Ga_{i,n})/\Ga_{i,n}$ sur
$\Ga_{i,n}\backslash(\Ga_{i,n+1}-\Ga_{i,n})$ et r\'eciproquement.
L'alphabet $\A$ est muni d'une involution $\beta\mapsto\ov\beta$, avec
$\ov\beta$ l'ar\^ete oppos\'ee de $\beta$ si $\beta\in EX$, et
$\overline{(g,+)}= (g^{-1},-)$ si $g\in
(\Ga_{i,n+1}-\Ga_{i,n})/\Ga_{i,n}$.

Nous allons d\'efinir une partition $\P$ de $\Ga\backslash \G T$,
index\'ee par $\A$, de la ma\-ni\`e\-re sui\-van\-te. Pour une
g\'eod\'esique $\ell$ dans $T$, notons $(\ell_i)_{i\in\ZZ}$ la suite
de ses ar\^etes cons\'ecutives, de sorte que
$t(\ell_{i-1})=o(\ell_i)=\ell(i)$. Notons qu'avec $\wt\tau$
l'application de renversement du temps, on a
$\wt\tau(\ell)_i=\overline{\ell_{-i-1}}$.  Pour tout $\beta$ dans
$\A$, d\'efinissons une partie non vide $E_\beta$ de $\G T$ par~:

\noindent\begin{minipage}{9cm}
si $\beta\in EX_0$, 
alors $E_\beta=\{\ell\in \G T\;:\; \ell_0=\wt\beta\;\}$,
\end{minipage}
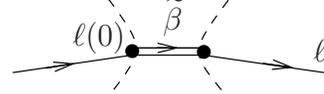
\begin{minipage}{5.8cm}
\input{fig_beta.pstex_t}
\end{minipage}

\smallskip\noindent\begin{minipage}{8cm}
si $\beta=a_{i,n}$, 
alors $E_\beta=\{\ell\in \G T\;:\; 
\ell_0=\wt{a_{i,n}},$ $\ell_1=\wt{a_{i,n+1}},\;
\ell_{-n}=\wt{a_{i,0}}\; \}$,
\end{minipage}
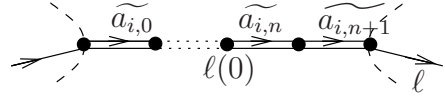
\begin{minipage}{6.9cm}$\;\;$
\input{fig_ain.pstex_t}
\end{minipage}

\smallskip\noindent\begin{minipage}{8.8cm} 
  si $\beta=(g,+)$ avec $g\in(\Ga_{i,n+1}-\Ga_{i,n})/\Ga_{i,n}$, 
  alors $E_\beta=\{\ell\in \G T\;:\; \ell_0=\wt{a_{i,n}},$
  $\ell_{1}=g\ov{\wt{a_{i,n}}},\; \ell_{-n}=\wt{a_{i,0}}\; \}$,
\end{minipage}
\begin{minipage}{6.1cm}$\;\;$
\input{fig_gplus.pstex_t}
\end{minipage}

\smallskip\noindent\begin{minipage}{8.8cm}
  si $\beta=(g,-)$ avec
  $g\in\Ga_{i,n}\backslash(\Ga_{i,n+1}-\Ga_{i,n})$, alors
  $E_\beta=\{\ell\in \G T\;:\; \ell_0=\ov{\wt{a_{i,n}}},$
  $\ell_{-1}=g^{-1}\wt{a_{i,n}},\; \ell_{n}=\ov{\wt{a_{i,0}}}\; \}$,
\end{minipage}
\begin{minipage}{6.1cm}$\;\;$
\input{ fig_gmoins.pstex_t}
\end{minipage}

\medskip
\noindent\begin{minipage}{8cm}
si $\beta=\overline{a_{i,n}}$, 
alors $E_\beta=\{\ell\in \G T\;:\; 
\ell_0=\ov{\wt{a_{i,n}}},$ $\ell_{-1}=\ov{\wt{a_{i,n+1}}},\;
\ell_{n}=\ov{\wt{a_{i,0}}}\; \}$.
\end{minipage}
\begin{minipage}{6.9cm}$\;\;$
\input{fig_ainbar.pstex_t}
\end{minipage} 

\medskip %
Les parties du second type codent les g\'eod\'esiques qui \`a
l'instant $t=0$ avancent dans un rayon cuspidal, et vont continuer
d'avancer.  Celles du troisi\`eme type codent les g\'eod\'esiques qui
\`a l'instant $t=0$ avancent dans un rayon cuspidal, et vont faire
demi-tour. Celles du quatri\`eme type codent les g\'eod\'esiques qui
\`a l'instant $t=0$ viennent de faire demi-tour. Enfin celles du
dernier type codent les g\'eod\'esiques qui \`a l'instant $t=0$
descendent dans un rayon cuspidal, et qui descendaient aussi \`a
l'instant d'avant.

\blemm %
Les parties $E_\beta$ de $\G T$, pour $\beta$ dans $\A$, sont
compactes, ouvertes et deux \`a deux disjointes. Leur r\'eunion $E$
est un domaine fondamental strict pour l'action de $\Ga$ (i.e.~$\Ga E=
\G T$ et si $\ga E\cap E$ est non vide pour un $\ga$ dans $\Ga$, alors
$\ga=1$).  
\elemm

\dem %
Comme $\Ga_{i,n}$ fixe l'ar\^ete $\wt{a_{i,n}}$ (et son ar\^ete
oppos\'ee), les parties $E_\beta$ sont bien d\'efinies. Par
d\'efinition de la topologie compacte-ouverte, et comme l'arbre $T$
est localement fini, elles sont bien compactes et ouvertes.  Comme
$g\ov{\wt{a_{i,n}}} \neq\wt{a_{i,n+1}}$ pour tout $g$ dans
$\Ga_{i,n+1}-\Ga_{i,n}$, ces parties sont deux \`a deux disjointes.
Toute g\'eod\'esique dans $T$ est \'equivalente modulo $\Ga$ \`a une
g\'eod\'esique $\ell$ telle que $\ell_0$ vaut ou bien $\wt\beta$ pour
un $\beta$ dans $EX_0$, ou bien $\wt{a_{i,n}}$ ou bien
$\ov{\wt{a_{i,n}}}$ pour un $n$ dans $\NN$ et $i$ dans
$\{1,\dots,k\}$. Dans les seconde et troisi\`eme alternatives, on a
respectivement $\pi(\ell_{-n})=a_{i,0}$ et $\pi(\ell_n)=\ov{a_{i,0}}$.
Donc, quitte \`a faire agir un \'el\'ement de $\Ga_{i,n}$, on peut
supposer respectivement que $\ell_{-n}=\wt{a_{i,0}}$ et
$\ell_{n}=\ov{\wt{a_{i,0}}}$. L'ensemble des ar\^etes d'origine
$t(\wt{a_{i,n}})$ est exactement
$$\{g\ov{\wt{a_{i,n}}}\;:\;g\in (\Ga_{i,n+1}-\Ga_{i,n})/\Ga_{i,n}\}
\cup\{\ov{\wt{a_{i,n}}},\wt{a_{i,n+1}}\}\;.$$ 
Comme une g\'eod\'esique n'a pas d'aller-retour, on obtient que $\Ga
E= \G T$.  Enfin, comme les stabilisateurs des $\wt{a_{i,0}}$ pour
$i$ dans $\{1,\dots,r\}$ et des $\wt\beta$ pour $\beta$ dans $EX_0$ sont
triviaux, la r\'eunion $E$ est un domaine fondamental strict.  \cqfd

\medskip %
En notant encore $E_\beta$ l'image de $E_\beta$ dans $\Ga\backslash \G
T$, on obtient donc une partition $\P=\{E_\beta\}_{\beta\in\A}$ de
$\Ga\backslash \G T$, par parties compactes, ouvertes et non vides.

\bigskip %
Consid\'erons la matrice de transition
$(A_{\alpha,\beta})_{\alpha,\beta\in\A}$ d\'efinie par
$A_{\alpha,\beta}=1$ si l'une des conditions suivantes est
v\'erifi\'ee

\smallskip\noindent$(1) \;\;\;\;\alpha,\beta\in EX,\;
t(\alpha)=o(\beta),\; \beta\neq\ov\alpha$\smallskip

ou

\smallskip\noindent$(2) \;\;\;\;
\exists \;i\in \{1,\dots, r\},\;\exists
\;n\in\NN,\;\exists\; g\in(\Ga_{i,n+1}-\Ga_{i,n})/\Ga_{i,n},$ 
$$\alpha\in EX, \beta=(g,+),\; t(\alpha)=o(a_{i,n}),\;
\alpha\neq\ov{a_{i,n}}$$
~~~~ou

\smallskip\noindent$(3) \;\;\;\;
\exists \;i\in \{1,\dots, r\},\;\exists
\;n\in\NN,\;\exists\; g\in(\Ga_{i,n+1}-\Ga_{i,n})/\Ga_{i,n},
\;\;\;\alpha=(g,+),\; \beta=(g^{-1},-)$\smallskip

ou

\smallskip\noindent$(4) \;\;\;\; \exists \;i\in \{1,\dots,
r\},\;\exists \;n\in\NN,\;\exists\;
g\in\Ga_{i,n}\backslash(\Ga_{i,n+1}-\Ga_{i,n}),$
$$\alpha=(g,-),\; \beta\in EX, \; o(\beta)=t(\overline{a_{i,n}}),
\;\beta\neq a_{i,n}$$
et $A_{\alpha,\beta}=0$ sinon.  Remarquons que
dans le cas (2), on a $\alpha=a_{i,n-1}$ si $n\geq 1$, et dans le cas
(4), on a $\beta=\overline{a_{i,n-1}}$ si $n\geq 1$.

Consid\'erons l'espace topologique produit $\A^\ZZ$, o\`u pour tout $x$
dans $\A^\ZZ$ on note $x_i$ la $i$-\`eme composante de $x$. Notons
$\sigma:\A^\ZZ\ra\A^\ZZ$ le d\'ecalage d\'efini par
$\sigma((x_j)_{j\in\ZZ})_i= x_{i+1}$.  Notons
$\kappa:\A^\ZZ\ra\A^\ZZ$ l'involution d\'efinie par
$\kappa((x_j)_{j\in\ZZ})_i= \overline{x_{-i-1}}$.  La matrice de transition
$(A_{\alpha,\beta})_{\alpha,\beta\in\A}$ d\'efinit un sous-d\'ecalage
(invariant par $\sigma$ et par $\kappa$)
$$X_A=\{(x_i)_{i\in\ZZ}\in \A^\ZZ\;:\; \forall \;i\in\ZZ,\;
A_{x_i,x_{i+1}}=1\}\;.$$

\bprop \label{prop:conjugtopo} %
L'application itin\'eraire $\Theta:\Ga\backslash \G T\ra X_A$,
d\'efinie par $\ell\mapsto(x_i)_{i\in\ZZ}$ o\`u pour tout $i$ dans
$\ZZ$, on a $\varphi^i(\ell)\in E_{x_i}$ est un hom\'eomorphisme,
rendant les diagrammes suivants commutatifs~:
$$\begin{array}{ccc}
\Ga\backslash \G T & \stackrel{\Theta}{\longrightarrow} &
 X_A \\
\varphi\downarrow\;\;\;\;  & & \;\;\;\downarrow \sigma \\
\Ga\backslash \G T & \stackrel{\Theta}{\longrightarrow} &
 X_A
\end{array}
\;\;\;\;\;\;\;
\begin{array}{ccc}
\Ga\backslash \G T & \stackrel{\Theta}{\longrightarrow} &
\!\!\! X_A\;\;\;\; \\
\tau\downarrow \;\;\;\; & & \;\;\;\downarrow \kappa\\
\Ga\backslash \G T & \stackrel{\Theta}{\longrightarrow} &
\!\!\!X_A \;\;\;\;
\end{array}\;.$$
\eprop

\dem On v\'erifie que pour tous les $\alpha,\beta$ dans $\A$ et $\ell$
dans $\G T$ tels que $\wt\varphi^{-1}(\ell)$ appartienne \`a $\Gamma
E_\alpha$, on a $A_{\alpha,\beta}=1$ si et seulement s'il existe une
g\'eod\'esique $\ell'$ dans $T$ telle que $\ell'$ et $\ell$
co\"{i}ncident sur $]-\infty,0\,]$ et $\ell'\in \Gamma E_\beta$. Ceci
montre que l'application $\Theta$ est bien \`a valeurs dans $ X_A$, et
qu'elle est surjective. Comme une g\'eod\'esique est d\'etermin\'ee
par la suite des ar\^etes qu'elle traverse et comme le stabilisateur
d'une ar\^ete de $X_0$ est trivial, l'application $\Theta$ est
injective. Les parties de la partition $\P$ sont d\'efinies par des
conditions ne portant que sur l'ar\^ete $\ell_0$ d'une g\'eod\'esique
$\ell$, et \'eventuellement sur un nombre fini d'ar\^etes
suppl\'ementaires. Par d\'efinition des topologies sur $\G T$ et
$\A^\ZZ$, on en d\'eduit que $\Theta$ est un hom\'eomorphisme.

\medskip %
Le premier diagramme est commutatif par d\'efinition de $\Theta$.
Pour montrer la commutativit\'e du second diagramme, on v\'erifie que
pour tout $\beta$ dans $\A$ et $\ell$ dans $\Ga\backslash \G T$, on a
(par une \'etude cas par cas pour la seconde \'equivalence)
$$\Theta\circ\tau(\ell)_0=\beta \;\Longleftrightarrow\; \tau(\ell)\in
E_\beta \;\Longleftrightarrow\; \varphi(\ell)\in E_{\ov{\beta}}
\;\Longleftrightarrow\;\kappa\circ\Theta(\ell)_0=\beta\;.$$ 
En utilisant que $\varphi\circ\tau=\tau\circ\varphi^{-1}$ et
$\sigma\circ\kappa=\kappa\circ\sigma^{-1}$, on conclut alors par la
commutativit\'e du premier diagramme.  
\cqfd

\bigskip %
\`A l'aide de la mesure de Bowen-Margulis $m_{\mbox{\tiny BM}}$,
construisons une mesure de pro\-ba\-bi\-li\-t\'e $\nu$ sur l'ensemble
discret d\'enombrable $\A$.  Posons, pour tout $\beta$ dans $\A$,
$$\nu_\beta=\nu(\{\beta\})=m_{\mbox{\tiny BM}}(E_\beta).$$
Comme $m_{\mbox{\tiny BM}}$ est une mesure de probabilit\'e sur 
$\Ga\backslash \G T$, et que $\P$ est une partition de $\Ga\backslash
\G T$, la mesure $\nu$ est bien de masse totale $1$. Comme les $E_\beta$
sont des ouverts non vides, et que $m_{\mbox{\tiny BM}}$ est une
mesure de support total, les $\nu_\beta$ sont strictement positifs.

Consid\'erons la matrice $\Pi=(\pi_{\alpha,\beta})
_{\alpha,\beta\in\A}$, d\'efinie par 
$$\pi_{\alpha,\beta}=m_{\mbox{\tiny BM}}[\varphi(\ell)\in E_\beta
\;|\;\ell\in E_\alpha]=\frac{m_{\mbox{\tiny
      BM}}(\varphi^{-1}(E_\beta)\cap E_\alpha )
}{\nu_\alpha}=\frac{m_{\mbox{\tiny BM}}(\varphi(E_\alpha)\cap E_\beta
  ) }{\nu_\alpha}\;.$$
Nous renvoyons par exemple \`a \cite{Kit} pour
les d\'efinitions sur les d\'ecalages de Markov.

La matrice $\Pi$ est une matrice stochastique sur $\A$, de mesure
stationnaire $\nu$. En effet, comme $\P$ est une partition de
$\Ga\backslash \G T$, et comme $m_{\mbox{\tiny BM}}$ est invariante
par $\varphi$, on a imm\'ediatement que $\sum_\beta \pi_{\alpha,\beta}
=1$ pour tout $\alpha$ et $\sum_\alpha \nu_\alpha\pi_{\alpha,\beta}
=\nu_\beta$ pour tout $\beta$.

Par d\'efinition, la {\it mesure de Markov} $\mu_\Pi$ associ\'ee \`a
$(\Pi,\nu)$ est l'unique mesure bor\'elienne de probabilit\'e sur
$\A^\ZZ$, invariante par $\sigma$, telle que, pour tout $k$ dans $\NN$
et tous $\alpha_0,\dots, \alpha_k$ dans $\A$, on a
$$\mu_\Pi[x_0=\alpha_0,\dots ,x_k=\alpha_k]=
\nu_{\alpha_0}\left(\prod_{i=0}^{k-1} \pi_{\alpha_i,
    \alpha_{i+1}}\right) \,.$$

\bprop \label{prop:conjugmesur} %
La mesure de Markov $\mu_\Pi$ est de support $X_A$, et
$\Theta_*m_{\mbox{\tiny BM}}=\mu_\Pi$.  
\eprop

\dem 
On remarque que $\pi_{\alpha,\beta}$ est nul si et seulement si
$A_{\alpha,\beta}$ est nul, ce qui montre, avec le fait que les
$\nu_\beta$ sont strictement positifs, que le support de $\mu_\Pi$
est $X_A$.

Puisque $m_{\mbox{\tiny BM}}$ est invariante par $\varphi$, comme
$\Theta$ est un hom\'eomorphisme qui conjugue $\varphi$ et $\sigma$,
la mesure $\Theta_*m_{\mbox{\tiny BM}}$ est une mesure bor\'elienne de
probabilit\'e, invariante par $\sigma$.  Comme $m_{\mbox{\tiny
    BM}}(E_{\alpha_0})= \nu_{\alpha_0}=\mu_\Pi[x_0=\alpha_0] $, pour
\'etablir que $\Theta_*m_{\mbox{\tiny BM}}$ et $\mu_\Pi$ sont
\'egales, il suffit de v\'erifier que $\Theta_*m_{\mbox{\tiny BM}}$
v\'erifie la propri\'et\'e de Markov, c'est-\`a-dire que pour tout $k$
dans $\NN$ et tous $\alpha_{-k},\dots, \alpha_{-1},\alpha_0,\beta$
dans $\A$, on a
$$\Theta_*m_{\mbox{\tiny BM}}[x_1=\beta|x_{0}=\alpha_{0},x_{-1}=
\alpha_{-1},\dots, x_{-k}=\alpha_{-k}]=\pi_{\alpha_{0}, \beta} \;.$$

On peut supposer, pour tout $i$ dans $\{1,\dots,k\}$, que
$A_{\alpha_{-i}, \alpha_{-(i-1)}}=1$ et que $A_{\alpha_{0}, \beta}=1$,
sinon le r\'esultat est imm\'ediat. Notons
$F_k=\varphi^{k}(E_{\alpha_{-k}})\cap\varphi^{k-1}(E_{\alpha_{-(k-1)}})\cap
\dots \cap E_{\alpha_0}$ et $F'_k=F_k\cap\varphi^{-1}(E_{\beta})$, qui
sont donc des compacts ouverts non vides. Le membre de gauche de
l'\'equation ci-dessus est \'egal \`a $m_{\mbox{\tiny
    BM}}(F'_k)/m_{\mbox{\tiny BM}}(F_k)$. Comme $\pi_{\alpha_{0},
  \beta}=m_{\mbox{\tiny BM}}(F'_0)/m_{\mbox{\tiny BM}}(F_0)$, il
suffit donc de montrer que $m_{\mbox{\tiny BM}}(F'_k)/m_{\mbox{\tiny
    BM}}(F_k)$ ne d\'epend pas de $k$.

Soit $\wt{\alpha_0}$ l'ar\^ete de $T$ d\'efinie comme le relev\'e
pr\'ec\'edemment choisi de $\alpha_0$ si $\alpha_0$ est dans $EX$, par
$\wt{\alpha_0}=\wt{a_{i,n}}$ si $\alpha_0=(g,+)$ avec $g\in
(\Ga_{i,n+1}-\Ga_{i,n})/\Ga_{i,n}$, et par
$\wt{\alpha_0}=\overline{\wt{a_{i,n}}}$ si $\alpha_0=(g,-)$ avec
$g\in\Ga_{i,n}\backslash(\Ga_{i,n+1}-\Ga_{i,n})$. Nous utiliserons le
param\'etrage de $\G T$ par $\partial_2T\times\ZZ$ d\'efini par le
point base $u=o(\wt{\alpha_0})$.

Montrons que $F_k$ (resp.~$F'_k$) est l'image injective par
l'application $\pi':\G T\ra\Ga\backslash \G T$ d'une partie $\wt F_k$
(resp.~$\wt F'_k$) de $\G T$ de la forme $V_k\times U \times \{0\}$
(resp.~$V_k\times U'\times\{0\}$), avec $U,U', V_k$ des parties de
$\partial T$, telles que $U,U'$ ne d\'ependent pas de $k$ et sont
disjointes de $V_k$, et telles que toute g\'eod\'esique entre un point
de $V_k$ et un point de $U$ (resp.~$U'$) passe par $u$.  Alors pour
tout $(\xi_-,\xi_+)$ dans $V_k\times U$ (resp.~$V_k\times U'$), on
aura $d_u(\xi_-,\xi_+)=1$. Donc, par d\'efinition de la mesure de
Bowen-Margulis (voir la partie \ref{subsec:groupgeofinmesbowmar}),
$$\frac{m_{\mbox{\tiny BM}}(F'_k)}{m_{\mbox{\tiny BM}}(F_k)}=\frac{\wt
  m_{\mbox{\tiny BM}}(\wt F'_k)}{\wt m_{\mbox{\tiny BM}}(\wt
  F_k)}=\frac{\mu_u(U')}{\mu_u(U)}\;,$$
qui ne d\'ependra pas de $k$, ce
qui montrera le r\'esultat.

Comme les $F'_k$ sont des compacts non vides d\'ecroissants en $k$, il
existe au moins une g\'eod\'esique $\ell^*$ (que l'on fixe) dans $T$ dont
l'image par $\pi'$ est un \'el\'ement de $\bigcap_{k\in\NN} F'_k$, telle
que $\ell^*_0=\wt{\alpha_0}$ (ce qui est possible, car $F'_k$ est
contenu dans $E_{\alpha_0}$).  Par d\'efinition des $E_\beta$, tout
\'el\'ement de $F_k$ (resp.~$F'_k$) se rel\`eve en au moins une g\'eod\'esique
$\ell$ de $T$ telle que $ \ell_{-k}=\ell^*_{-k},\dots, \ell_{0}=
\ell^*_{0}$ (resp.~$\ell_{-k}=\ell^*_{-k},\dots,\ell_{0}=
\ell^*_{0},\ell_{1}=\ell^*_{1}$), et en particulier $\ell$ passe \`a
l'instant $t=0$ par le point $u$. Nous allons maintenant d\'efinir $V_k$
et $U$ (resp.~$U'$) en discutant suivant les valeurs de $\alpha_{-k}$ et
$\alpha_0$ (resp.~$\beta$), voir le dessin ci-dessous.

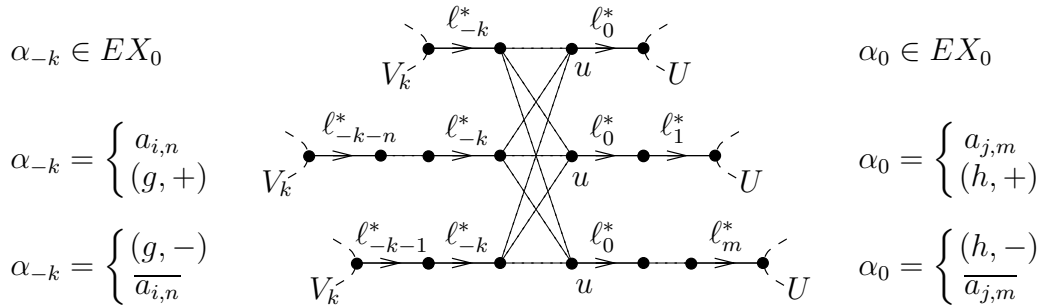
\begin{center}$\!\!\!\!\!\!\!\!\!\!\!\!\!\!\!\!
\!\!\!\!\!\!\!\!\!\!\!\!\!\!\!\!$
\input{fig_mesuinvabis.pstex_t}
\end{center}

\medskip 
\noindent
Si $\alpha_{-k}$ est dans $EX_0$, posons $V_k=\partial_{\ell^*_{-k}}T$. \\
 Si $\alpha_{-k}$ vaut $a_{i,n}$ ou $(g,+)$ avec $g\in (\Ga_{i,n+1}-
\Ga_{i,n})/\Ga_{i,n}$, alors posons $V_k= \partial_{\ell^*_{-k-n}}T$. \\
 Si $\alpha_{-k}$ vaut $(g,-)$ ou $\overline{a_{i,n}}$ avec $g\in
\Ga_{i,n}\backslash (\Ga_{i,n+1}-\Ga_{i,n})$, alors posons $V_k=
\partial_{\ell^*_{-k-1}}T$.\\
Si $\alpha_0$ est dans $EX_0$, posons $U=\partial_{\overline{\ell^*_{0}}}T$. \\
Si $\alpha_0$ vaut $a_{j,m}$ ou $(h,+)$ avec $h\in (\Ga_{j,m+1}-
\Ga_{j,m}) /\Ga_{j,m}$, alors posons $U= \partial_{\overline{\ell^*_1}}T$.\\
Si $\alpha_0$ vaut $(h,-)$ ou $\overline{a_{j,m}}$ avec
$h\in\Ga_{j,m}\backslash(\Ga_{j,m+1}-\Ga_{j,m})$, alors posons
$U=\partial_{\overline{\ell^*_m}}T$.\\
Si $\beta$ est dans $EX_0$, posons $U'=\partial_{\overline{\ell^*_{1}}}T$. \\
Si $\beta$ vaut $a_{j,m}$ ou $(h,+)$ avec $h\in(\Ga_{j,m+1}-\Ga_{j,m})
/\Ga_{j,m}$, alors posons $U'= \partial_{\overline{\ell^*_2}}T$.\\
Si $\beta$ vaut $(h,-)$ ou $\overline{a_{j,m}}$ avec
$h\in\Ga_{j,m}\backslash (\Ga_{j,m+1}-\Ga_{j,m})$, alors posons $U'=
\partial_{\overline{\ell^*_{m+1}}}T$.

\medskip %
Par construction, les parties $U$ et $U'$ ne d\'ependent pas de $k$ et
sont disjointes de $V_k$. Toute g\'eod\'esique entre un point de $V_k$
et un point de $U$ ou $U'$ passe par $u$. On v\'erifie que si $\wt
F_k=V_k\times U\times \{0\}$ (resp.~$\wt F'_k=V_k\times
U'\times\{0\}$), alors par construction des $E_\ga$, tout \'el\'ement
de $F_k$ (resp.~$F'_k$) admet un unique relev\'e dans $\wt F_k$
(resp.~$\wt F'_k$). Donc la restriction de $\pi'$ \`a $\wt F_k$
(resp.~$\wt F'_k$) est injective, d'image $F_k$ (resp.~$F'_k$).  Le
r\'esultat en d\'ecoule.  \cqfd

\bigskip %
Ainsi par les propositions \ref{prop:conjugtopo} et
\ref{prop:conjugmesur} ci-dessus, les syst\`emes dynamiques mesur\'es
$(\Ga\backslash \G T,m_{\mbox{\tiny BM}},\varphi)$ et
$(X_A,\mu_\Pi,\sigma)$ sont conjugu\'es par $\Theta$.

\bprop\label{prop:entropfini}
L'entropie de la partition $\P$  est finie.
\eprop

\dem 
Notons pour simplifier $\Ga_i=\Ga_{\xi_i}$, $\delta=
\delta_{\Ga}$, $\delta_i= \delta_{\Ga_i}$ et $u_i=o(\wt{a_{i,0}})$.

Rappelons que l'entropie $h_\P$ de la partition $\P$ pour la mesure
$m_{\mbox{\tiny BM}}$ est
$$h_\P=\sum_{\alpha\in\A}\nu_{\alpha}(-\log\nu_{\alpha})\;.$$

Calculons $\nu_{\alpha}$ pour $\alpha$ dans $\A-EX_0$ (car $EX_0$ est
fini).  Par d\'efinition, la partie $E_\alpha$ de $\G T$ est, dans le
param\'etrage de $\G T$ d\'efini par le point base $u_i$ (pour un
certain $i$ d\'ependant de $\alpha$, pr\'ecis\'e ci-dessous), de la forme
$V\times U\times \{t_\alpha\}$ pour un $t_\alpha$ dans $\ZZ$ et $U,V$
deux parties disjointes de $\partial T$, telles que toutes les
g\'eod\'esiques entre $V$ et $U$ passent par $u_i$. Donc

\smallskip\noindent\begin{minipage}{8.5cm}
si $\alpha=a_{i,n}$, 
alors $\nu_\alpha=\mu_{u_i}(\partial_{\wt{a_{i,0}}}T)
\mu_{u_i}(\partial_{\overline{\wt{a_{i,n+1}}}}T)$,
\end{minipage}
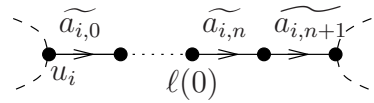
\begin{minipage}{6.4cm}$\;\;\;\;$
\input{fig_ainbis.pstex_t}
\end{minipage}

\medskip
\noindent\begin{minipage}{8.5cm}
si $\alpha=\ov{a_{i,n}}$, 
alors $\nu_\alpha=\mu_{u_i}(\partial_{\ov{\wt{a_{i,n+1}}}}T)
\mu_{u_i}(\partial_{\overline{\ov{\wt{a_{i,0}}}}}T)$,
\end{minipage}
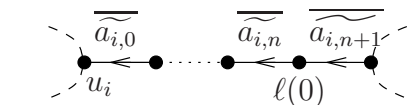
\begin{minipage}{6.4cm}$\;\;\;\;\;\;\;\;$
\input{fig_ainbarbis.pstex_t}
\end{minipage} 

\smallskip\noindent\begin{minipage}{8.8cm}
  si $\alpha=(g,+)$ avec $g\in(\Ga_{i,n+1}- \Ga_{i,n})/\Ga_{i,n}$,
  alors $\nu_\alpha=\mu_{u_i}(\partial_{\wt{a_{i,0}}}T)
  \mu_{u_i}(\partial_{\overline{g\ov{\wt{a_{i,n}}}}}T)$,
\end{minipage}
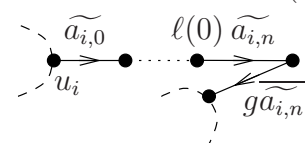
\begin{minipage}{6.1cm}$\;\;$
\input{fig_gplusbis.pstex_t}
\end{minipage}

\smallskip\noindent\begin{minipage}{8.8cm}
si $\alpha=(g^{-1},-)$ avec $g\in(\Ga_{i,n+1}- \Ga_{i,n})/\Ga_{i,n}$, 
alors $\nu_\alpha=\mu_{u_i}(\partial_{g\wt{a_{i,n}}}T)
\mu_{u_i}(\partial_{\overline{\ov{\wt{a_{i,0}}}}}T)$.
\end{minipage}
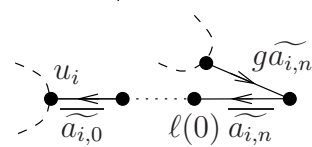
\begin{minipage}{6.1cm}$\;\;\;\;\;\;\;\;\;\;\;\;$
\input{fig_gmoinsbis.pstex_t}
\end{minipage}

\bigskip %
Par le th\'eor\`eme de la densit\'e fluctuante de Sullivan (voir
\cite[Theo.~4.1]{HP}), si $\xi_i$ est un point parabolique born\'e, si
$\Ga$ est de type divergent, et s'il existe une constante $c_0\geq 1$
telle que $$\frac{1}{c_0} \,e^{\,\delta_i n}\leq {\rm
  Card}\{\ga\in\Ga_i\;:\; d(u_i,\ga u_i)\leq n\}\leq c_0
\,e^{\,\delta_i n}\;,$$ alors il existe une constante $c_1\geq 1$
telle que pour tout $n$ dans $\NN$, on a
$$\frac{1}{c_1} e^{2(\delta_i-\delta)n} \leq \mu_{u_i}
\big(\partial_{\overline{\wt{a_{i,n}}}}T\big) \leq c_1
e^{2(\delta_i-\delta)n}\;.$$ (Pour les vari\'et\'es riemaniennes
(compl\`etes), cette propri\'et\'e a \'et\'e d\'emontr\'ee par
Sullivan dans le cas de volume fini et de courbure n\'egative
constante, et par Stratmann et Velani dans le cas g\'eom\'etriquement
fini et de courbure n\'egative constante. Dans \cite{HP}, elle est
\'enonc\'ee dans le cas de la courbure n\'egative variable, mais la
preuve est faite d'une part par des arguments de comparaison avec les
arbres, qui sont valables dans notre cas, et d'autre part par des
arguments de limites, mais une limite (pour la topologie de
Hausdorff-Gromov point\'ee) d'arbres de valences uniform\'ement
born\'ees est encore un arbre localement fini. Voir aussi
\cite[Prop.~3.2]{HP2}.) Notons que les hypoth\`eses de ce r\'esultat
sont v\'erifi\'ees, par la proposition \ref{prop:finitudebowmarg}.
Comme de plus $\delta_i=\delta/2$ par la proposition
\ref{prop:finitudebowmarg}, on~a
\begin{equation}\label{eq:asymphoroboule}
 \frac{1}{c_1} e^{-\delta n} \leq
\mu_{u_i} \big(\partial_{\overline{\wt{a_{i,n}}}}T\big)
\leq c_1 e^{-\delta n}
\;.\end{equation}

Remarquons, pour $g\in \Ga_{i,n+1}-\Ga_{i,n}$, que $o(g\wt{a_{i,0}})$
appartient \`a l'orbite de $u_i$, et que $d(u_i,o(g\wt{a_{i,0}}))=2n$.
Par le lemme de l'ombre de Sullivan (voir par exemple
\cite[Lem.~1.3]{Rob} dans un contexte plus g\'en\'eral), il existe
donc une constante $c_2\geq 1$ telle que pour tout $n$ dans $\NN$ et
$g$ dans $\Ga_{i,n+1}-\Ga_{i,n}$,
\begin{equation}\label{eq:asymphorosphere}
\frac{1}{c_2} e^{-2n\delta}\leq \mu_{u_i} 
\big(\partial_{g\wt{a_{i,0}}}T\big)\leq c_2 e^{-2n\delta}
\;.\end{equation}

Remarquons que $\partial_{\overline{g\ov{\wt{a_{i,0}}}}}T\subset
\partial_{\overline{g\ov{\wt{a_{i,n}}}}}T\subset\partial_{\overline{\wt{a_{i,n}}}}T$.  
Posons $c^-_3=\min_{1\leq i\leq r}\mu_{u_i}(\partial_{\wt{a_{i,0}}}T)$
et $c^+_3=\max_{1\leq i\leq r}\mu_{u_i}(\partial_{\wt{a_{i,0}}}T)$.
Comme les valences de $T$ sont uniform\'ement born\'ees, il existe un
$N$ dans $\NN$ qui majore le cardinal de $(\Ga_{i,n+1}-\Ga_{i,n})
/\Ga_{i,n}$ pour tout $i$ et $n$.  Donc
$$\sum_{1\leq i\leq r,\;n\in\NN} -\left(\nu_{a_{i,n}}\log
  \nu_{a_{i,n}}+\nu_\ov{a_{i,n}}\log \nu_\ov{a_{i,n}}\right)
\leq\sum_{n\in\NN} 2rc_1c^+_3e^{-\delta(n+1)} (\delta (n+1) +
\log\frac{c_1}{c^-_3}) \;,$$
et, la premi\`ere somme suivante portant sur les $g$ dans
$(\Ga_{i,n+1}-\Ga_{i,n}) /\Ga_{i,n}$, pour $1\leq i\leq r$ et $n$ dans
$\NN$, 
$$
\sum - \left(\nu_{(g,+)}\log \nu_{(g,+)}+\nu_{(g^{-1},-)}\log
  \nu_{(g^{-1},-)}\right)\leq \sum_{n\in\NN} 2Nrc_1c^+_3e^{-\delta n}
(2\delta n + \log\frac{c_2}{c^-_3}) \;.$$
Comme $\delta>0$ (voir la proposition \ref{prop:finitudebowmarg}), ces
deux sommes convergent, ce qui montre le r\'esultat.  
\cqfd

\bigskip %
Terminons maintenant la preuve du th\'eor\`eme \ref{theo:geodBernou}.
Par les propositions \ref{prop:conjugtopo} et \ref{prop:conjugmesur},
les syst\`emes dynamiques mesur\'es $(\Ga\backslash \G
T,m_{\mbox{\tiny BM}},\varphi)$ et $(X_A,\mu_\Pi,\sigma)$ sont
conjugu\'es par $\Theta$.

Supposons tout d'abord que $L_\Ga=\ZZ$. Alors $(\Ga\backslash \G
T,m_{\mbox{\tiny BM}},\varphi)$ est m\'elangeant par la proposition
\ref{prop:melang}, donc $(X_A,\mu_\Pi,\sigma)$ aussi.
L'hom\'eomorphisme $\Theta$ envoie par construction la partition $\P$
de $\Ga\backslash \G T$ sur la partition g\'en\'eratrice
$\{\{(x_i)_{i\in\ZZ}\;:\; x_0=\alpha\}\}_{\alpha\in\A}$ de $X_A$. Donc
par la proposition \ref{prop:entropfini}, le syst\`eme
$(X_A,\mu_\Pi,\sigma)$ est un d\'ecalage de Markov m\'elangeant sur un
alphabet d\'enombrable, de partition g\'en\'eratrice d'entropie finie.
Il est bien connu qu'un tel d\'ecalage est conjugu\'e \`a un
d\'ecalage de Bernoulli, voir par exemple \cite{FO} lorsque l'alphabet
est fini, et \cite[sect.~8]{Tho} pour l'extension \`a un alphabet
d\'enombrable de partition g\'en\'eratrice d'entropie finie.

Supposons maintenant que $L_\Ga=2\ZZ$. Alors $(\Ga\backslash \G
T,m_{\mbox{\tiny BM}},\varphi)$ n'est pas m\'elangeant, mais
$(\Ga\backslash \G_0 T,m_{\mbox{\tiny BM}},\varphi^2)$ l'est, par la
proposition \ref{prop:melang}. En rempla\c{c}ant les ar\^etes par les
suites de deux ar\^etes cons\'ecutives d'extr\'emit\'es deux sommets
de $T$ \`a distance paire du point base de $T$, on construit de
mani\`ere analogue \`a ce qui pr\'ec\`ede un codage de $(\Ga\backslash
\G_0 T,m_{\mbox{\tiny BM}},\varphi^2)$ par un d\'ecalage de Markov sur
un alphabet d\'enombrable, de partition g\'en\'eratrice d'entropie
finie.  Le th\'eor\`eme \ref{theo:geodBernou} en d\'ecoule. \cqfd

\bigskip 
Le th\'eor\`eme \ref{theo:intro} de l'introduction est un
corollaire du th\'eor\`eme \ref{theo:geodBernou}, par la discussion
pr\'ec\'edant la proposition \ref{prop:melang}.

\section{Le flot g\'eod\'esique sur un graphe de groupes}
\label{sec:exappendice}

Dans cette partie, nous nous int\'eressons au codage du flot
g\'eod\'esique sur un arbre quotient\'e par un sous-groupe
d'automorphismes, sans supposer le sous-groupe discret ni m\^eme
l'arbre localement fini.

\subsection{Flot g\'eod\'esique d'ordre $1$ sur un graphe 
de groupes}
\label{subseq:ordreun}

Soit $(X,G_*)$ un graphe de groupes. Posons~:
    $$\Omega=\coprod_{e,e'\in EX\,:\; o(e)=t(e')}
\rho_{e'}(G_{e'})\backslash G_{o(e)}/\rho_{\overline{e}}(G_e)\;.$$

On munit $\Omega\times EX$ de la topologie discr\`ete et
${\displaystyle\Sigma=(\Omega\times EX)^\ZZ}$ de la topologie produit.
D'ailleurs, si $(X,G_\ast)$ est un graphe fini de groupes finis, alors
$\Omega\times EX$ est fini. On note $\sigma$ le d\'ecalage \`a gauche
sur $\Sigma$.

Notons $\Sigma_0$ l'ensemble des \'el\'ements $(g_i,e_i)_{i\in \ZZ}$
de $\Sigma$ tels que, pour tout $i$ dans $\ZZ$,
$$
\left\{\begin{array}{l} (1)\;\; o(e_{i})=t(e_{i-1})\\
    (2)\;\; g_i\in\rho_{e_{i-1}}(G_{e_{i-1}})\backslash
    G_{o(e_i)}/\rho_{\overline{e_i}}(G_{e_i})
    \\
    (3) \;\;{\rm si}\; e_{i}=\overline{e_{i-1}}, \;{\rm alors}\;
    g_i\notin \rho_{\overline{e_{i}}}(G_{e_{i}})\;.
\end{array}\right.
$$
Cette derni\`ere condition dit simplement que si
$e_{i}=\overline{e_{i-1}}$, alors $g_i$ n'est pas la double classe
triviale.

\medskip 
On appelle {\it flot g\'eod\'esique d'ordre $1$ sur
  $(X,G_\ast)$} le sous-d\'ecalage $(\Sigma_0,\sigma)$.  Il est
imm\'ediat que $\sigma$ pr\'eserve $\Sigma_0$, et que les conditions
d\'efinissant $\Sigma_0$ sont ``locales''. Ainsi, lorsque $(X,G_\ast)$
est un graphe fini de groupes finis, le sous-d\'ecalage
$(\Sigma_0,\sigma)$ est de type fini.

\bprop\label{prop:flotgeod_ordun} Soit $T$ un arbre, et $\Gamma$ un
sous-groupe de ${\rm Aut}(T)$. Si $(\Sigma_0,\sigma)$ est le flot
g\'eod\'esique d'ordre $1$ du graphe de groupes quotient $\Ga\bac T$,
alors il existe une application continue surjective
$\theta:\Ga\backslash \G T \ra \Sigma_0$ qui rend le diagramme suivant
commutatif
$$\begin{array}[b]{ccc} \Ga\backslash \G T &
\stackrel{\varphi}{\longrightarrow} & \Ga\backslash \G T \\
\theta\downarrow\;\; & & \;\;\downarrow \theta \\ \Sigma_0&
\stackrel{\sigma}{\longrightarrow}& \Sigma_0
\end{array}\;.$$
\eprop

\dem 
Soit $(X,G_\ast)=\Ga\bac T$ et $\pi:T\ra X=\Ga\backslash T$ la
projection canonique.  Dans un premier temps, on construit une
application $\wt\theta:\G T\ra \Sigma_0$. Soit $f$ une application
simpliciale de $\RR$ dans $T$. Pour $i$ dans $\ZZ$, notons $f_i$
l'ar\^ete $f([i,i+1])$ d'origine $f(i)$ et posons $e_{i}=\pi(f_i)$. On a
alors $t(e_{i-1})=\pi\circ f(i)=o(e_{i})$.

On reprend les notations de la d\'efinition du graphe de groupes
quotient $\Ga\bac T$ dans la partie \ref{subsec:graphgroupflogeod},
o\`u l'on avait fix\'e des relev\'es $\wt e$ et $\wt v$ dans $T$ d'une
ar\^ete $e$ et d'un sommet $v$ de $X$, et choisi un \'el\'ement $g_e$ dans
$\Ga$ pour toute ar\^ete $e$. Comme $\pi(f_i)=\pi(\wt{e_i})=e_i$, il
existe $h_i$ dans $\Ga$ tel que $h_if_i=\wt{e_i}$. L'\'el\'ement $h_i$
est bien d\'efini modulo multiplication \`a gauche par un \'el\'ement
du fixateur de l'ar\^ete $\wt{e_{i}}$, qui est un \'el\'ement de
$G_{e_i}$. Posons
$$g_i=g_{e_{i-1}}^{-1}h_{i-1}h_i^{-1}g_{\overline{e_i}}\;.
\;\;\;\;\;\;(*)$$ L'\'el\'ement $g_i$ fixe $\wt{\pi\circ f(i)}$, il
appartient donc \`a $G_{o(e_{i})}=G_{t(e_{i-1})}$ (voir la figure
suivante).  L'\'el\'ement $g_i$ est d\'efini modulo multiplication \`a
gauche par un \'el\'ement de $g_{e_{i-1}}^{-1}G_{e_{i-1}}g_{e_{i-1}}
=\rho_{e_{ i-1}}(G_{e_{i-1}})$, ainsi que modulo multiplication \`a
droite par un \'el\'ement de
$g_{\overline{e_i}}^{-1}G_{e_i}g_{\overline{e_i}}=
\rho_{\overline{e_i}}(G_{e_i})$.  Nous noterons encore $g_i$ la double
classe de $g_i$ dans le double quotient $\rho_{e_{i-1}}(G_{e_{i-1}})
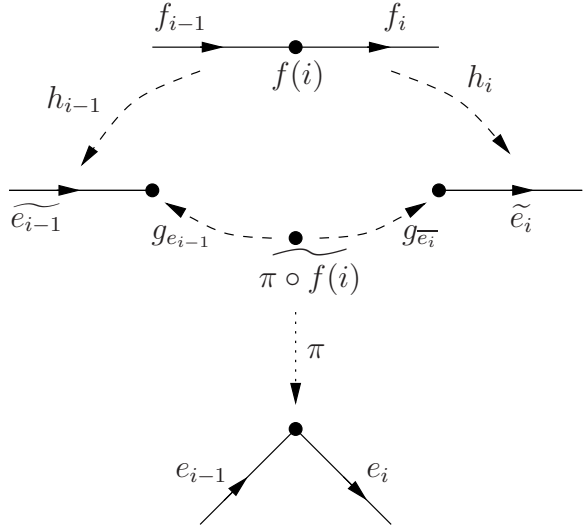
\backslash G_{o(e_i)}/\rho_{\overline{e_i}}(G_{e_i})$.

\begin{center}
\input{fig_ledrappi.pstex_t}
\end{center}

Remarquons que si $e_i=\overline{e_{i-1}}$, alors
$$\overline{f_{i-1}}=h_{i-1}^{-1}\overline{\wt{e_{i-1}}}=
h_{i-1}^{-1}\wt{\overline{e_{i-1}}}=h_{i-1}^{-1}\wt{e_{i}}\;.$$ Donc
$f$ est localement injective si et seulement si pour tout $i$ dans
$\ZZ$, on a $f_{i}\neq\overline{f_{i-1}}$, donc si et seulement si

\begin{itemize}
\item ou bien $e_i\neq\overline{e_{i-1}}$~;
\item ou bien $e_i=\overline{e_{i-1}}$ et 
$g_ig_{\overline{e_i}}^{-1}\wt{e_i}\neq 
g_{\overline{e_i}}^{-1}\wt{e_i}\;.$
\end{itemize}
Cette derni\`ere in\'egalit\'e est \'equivalente \`a la non-appartenance
de $g_i$ \`a $\rho_{\overline{e_i}}(G_{e_i})$.

Posons alors $\wt\theta(f)=(g_i,e_i)_{i\in\ZZ}$. Par construction,
l'application $\wt\theta$ est \`a valeurs dans $\Sigma_0$, et est
invariante par $\Gamma$, elle induit donc une application
$\theta:\Ga\backslash \G T \ra\Sigma_0$. Par construction, on a
$\theta\circ \varphi=\sigma\circ\theta$.

\medskip %
Montrons que $\wt\theta$ est surjective, ce qui entra\^{\i}ne que
$\theta$ est surjective. Soit $(g_i,e_i)_{i\in\ZZ}$ un \'el\'ement de
$\Sigma_0$. Pour tout entier relatif $i$, on choisit un repr\'esentant
dans la double classe $g_i$, encore not\'e $g_i$.  Posons
$f_0=\wt{e_0}$ et $h_0={\rm id}$.  Montrons par r\'ecurrence sur
$n\geq 1$ que pour tout $1\leq k\leq n-1$, il existe des ar\^etes
$f_k$ de $T$ et des \'el\'ements $h_k$ de $G$ tels que
 $$t(f_{k-1})=o(f_k), \;\;h_k^{-1}\wt{e_k}=f_k\;\;{\rm ~et~}\;\;
 g_k=g_{e_{k-1}}^{-1}h_{k-1} h_k^{-1} g_{\overline{e_k}}\;.$$ Pour
 $n=1$, on prend $h_{1}=g_{\overline{e_{1}}}g_{1}^{-1}g_{e_{0}}^{-1}$
 et $f_{1}=h_{1}^{-1}\wt e_{1}$, alors~:
 $o(f_{1})=g_{e_{0}}g_{1}g_{\overline{e_{1}}}^{-1}o(\wt
 e_{1})=g_{e_{0}}g_{1}\wt {o(e_{1})}=g_{e_{0}}\wt {o(e_{1})}=t(\wt
 {e_{0}})=t(f_{0})$. Supposons que l'hypoth\`ese de r\'ecurrence soit
 satisfaite au rang $n$. Posons~:
$$h_n=(h_{n-1}^{-1}g_{e_{n-1}}g_ng_{\overline{e_n}}^{-1})^{-1}\;.$$
Alors la relation $g_n=g_{e_{n-1}}^{-1}h_{n-1} h_n^{-1}
g_{\overline{e_n}}$ est v\'erifi\'ee, et $h_n^{-1}$ envoie $o(\wt{e_n})$
sur $t(f_{n-1})$.  En effet,
\begin{eqnarray*}
h_n^{-1}o(\wt {e_n})
&=&
h_{n-1}^{-1}g_{e_{n-1}}g_ng_{\overline{e_n}}^{-1}o(\wt{e_n})=
h_{n-1}^{-1}g_{e_{n-1}}g_n\wt {o(e_n)}\\
&=&
h_{n-1}^{-1}g_{e_{n-1}}\wt{o(e_n)}=h_{n-1}^{-1}t(\wt{e_{n-1}})
=t(f_{n-1})\;.
\end{eqnarray*}
Maintenant, posons $f_n=h_n^{-1}\wt{e_n}$. L'origine de $f_n$ est bien
$t(f_{n-1})$.  Pour $1\leq k\leq n$, les ar\^etes $f_k$ et les
\'el\'ements $h_k$ de $G$ v\'erifient donc l'hypoth\`ese de
r\'ecurrence au rang $n+1$.

On fait une construction similaire pour $n<0$. Notons $f:\RR\ra T$
l'application simpliciale telle que pour tout $i$ dans $\ZZ$,
$f([i,i+1])=f_i$. Alors $f$ appartient \`a $\G T$ (par la condition
(3) de la d\'efinition de $\Sigma_0$) et $\wt\theta(f)=
(g_i,e_i)_{i\in\ZZ}$, ce qui montre la surjectivit\'e de $\wt\theta$.

\medskip %
Montrons maintenant que $\wt\theta$ est continue, ce qui entra\^{\i}ne
que $\theta$ est continue.  Si $f$ et $f'$ sont deux g\'eod\'esiques
de $\G T$, posons $\wt\theta(f)=(g_i,e_i)_{i\in\ZZ}$ et
$\wt\theta(f')=(g'_i,e'_i)_{i\in\ZZ}$.  Supposons que $f$ et $f'$
co{\"\i}ncident sur $[-N,N]$, alors pour tout entier $i$ de $[-N,N]$,
on a $e_i=\pi (f_i)=\pi (f_i')=e'_i$, et on peut supposer que $h_i$ et
$h'_i$ co\"{\i}ncident.  La formule $(*)$ montre alors que $g_i$ et
$g'_i$ co\"{\i}ncident pour $i$ entier dans $[-(N-1),N-1]$. Ceci
montre que $\wt\theta$ est continue sur $\Gamma\backslash\G T$ par
d\'efinition de la topologie-produit.  \cqfd

\bigskip %
Malheureusement, l'application $\theta$ n'est pas toujours injective.
Elle l'est si les stabilisateurs d'ar\^etes $G_e$, pour $e$ dans $EX$,
sont triviaux. Mais quand il existe un \'el\'ement (non trivial) de
$\Ga$ qui fixe un rayon g\'eod\'esique de $T$, cela cr\'ee des
probl\`emes de propagation d'\'el\'ements dans des stabilisateurs
d'ar\^etes.

\subsection{Actions acylindriques de groupes sur les arbres}
\label{subseq:actacy}

Soit $k$ un entier strictement positif. Comme le fait Z.~Sela
\cite{Sel}, nous dirons qu'une action simpliciale (sans inversion)
d'un groupe $\Gamma$ sur un arbre $T$ est {\it $k$-acylindrique} si le
fixateur d'un chemin (localement) injectif de $k$ ar\^etes
cons\'ecutives est trivial. Nous dirons qu'une action est {\it
  acylindrique} s'il existe $k$ dans $\NN-\{0\}$ tel que l'action soit
$k$-acylindrique.  Par exemple, une action est $1$-acylindrique si
elle est \`a stabilisateurs d'ar\^ete triviaux.

\medskip \rem %
En fait, dans la d\'efinition originale, Z.~Sela demande aussi que
l'action soit {\it minimale} (i.e.~sans sous-arbre invariant non vide
propre) et {\it r\'eduite} (i.e.~si le stabilisateur d'un sommet $v$
admet exactement deux orbites d'ar\^etes d'origine $v$, alors chacune
de ces deux orbites est de cardinal au moins $2$.) Ici nous n'aurons
pas besoin de ces deux hypoth\`eses.

\medskip \rem Beaucoup d'actions de groupes sur des arbres sont
acylindriques, voir \cite{Sel}, mais il existe des arbres $T$
localement finis et des sous-groupes discrets cocompacts $\Ga$ de
${\rm Aut}(T)$, tels que l'action de $\Ga$ sur $T$ ne soit pas
acylindrique. Par exemple, l'action du groupe fondamental du graphe
fini de groupes finis suivant \\
\centerline{\input{fig_contrexacyl.pstex_t}}\\
sur son arbre de Bass-Serre n'est pas acylindrique.

\medskip
Avant de donner dans la partie \ref{subseq:ordreka} un th\'eor\`eme de
codage pour des actions acylindriques, nous donnons quelques exemples.
Le r\'esultat suivant est sans doute bien connu, nous n'en donnons une
preuve que par souci de compl\'etude.

\bprop \label{prop:ratrapage} Soit $T$ un arbre localement fini, $\Ga$
un sous-groupe de ${\rm Aut}(T)$, tel que $\Ga\backslash T$ soit fini
et les stabilisateurs d'ar\^ete soient finis. Alors l'action de $\Ga$ sur
$T$ n'est pas acylindrique si et seulement s'il existe  des \'el\'ements
$h,g$ de $\Ga$ tels que
\begin{itemize}
\item$h$ est non trivial et admet un point fixe dans $T$,
\item $g$ n'admet pas de point fixe dans $T$,
\item $h$ et $g$ commutent.
\end{itemize}
\eprop

\dem %
Supposons qu'il existe $h$ et $g$ comme dans l'\'enonc\'e.  Comme $g$
n'a pas de point fixe, alors il admet un axe de translation $A_g$. Cet
ensemble $A_{g}$ est (l'image d')une g\'eod\'esique de $T$, et est
invariant par $g$.  Sur $A_{g}$, l'application $g$ induit une
translation de distance de translation $\lambda=\inf_{x\in
  T}d(x,gx)>0$.  D'ailleurs, $A_{g}$ est l'ensemble des points $x$ de
$T$ tels que $d(x,gx)=\lambda$ (voir par exemple \cite{Ser}).

Comme $h$ commute avec $g$, l'application $h$ pr\'eserve $A_g$.
L'ensemble des points fixes de $h$ (qui est non vide) rencontre tout
sous-arbre invariant par $h$ et non vide de $T$ (voir par exemple
\cite{Ser}). Donc $A_{g}$ contient au moins un point fixe $x$ de $h$,
et pour tout $n$ dans $\ZZ$, l'isom\'etrie $h$ fixe $g^{n }x$. Donc
$h$ fixe le chemin g\'eod\'esique de longueur $n\lambda$ compris entre
$x$ et $g^{n}x$.  Donc l'action n'est pas acyclindrique.

R\'eciproquement, soit $N$ dans $\NN$ tel que $N-1$ soit le maximum
des cardinaux des stabilisateurs des ar\^etes de $T$ (ces
stabilisateurs sont en nombre fini modulo conjugaison). Soit $N'$ le
nombre d'ar\^etes (tenant compte des deux orientations possibles) de
$\Ga\backslash T$.  Soit $c$ un chemin d'ar\^etes g\'eod\'esique dans
$T$, fix\'e par un \'el\'ement non trivial $h'$ de $\Ga$, dont la
longueur est strictement sup\'erieure \`a $k=N'(N+1)$. Alors il existe
au moins $N+1$ ar\^etes de $c$, not\'ees $(e_i)_{i\in\{0,\dots, N\}}$,
orient\'ees compatiblement le long de $c$ et num\'erot\'ees dans
l'ordre o\`u on les rencontre en parcourant $c$, dont les images par
la projection canonique $\pi:T\ra \Ga\backslash T$ co\"{\i}ncident (en
tenant compte aussi de l'orientation). Il existe donc $g_1,\dots, g_N$
dans $\Ga$ tels que $g_i$ envoie $e_0$ sur $e_i$ en pr\'eservant
l'orientation.  Comme $h'$ fixe $c$, il existe $h_i$ dans le
stabilisateur de $e_0$ tel que $h'=g_ih_ig_i^{-1}$.  Par d\'efinition
de $N$, il existe $i\neq j$ tel que $h_i=h_j$.  On peut toujours
supposer (quitte \`a permuter) que l'ar\^ete $e_j$ est entre l'ar\^ete
$e_0$ et l'ar\^ete $e_i$ sur $c$. On pose alors $h=h_i=h_j$.
L'\'el\'ement $h$ est bien non trivial (car $h'$ l'est) et fixe un
point de $T$.

Rappelons (voir par exemple \cite{Ser}) que si une isom\'etrie
$\gamma$ d'un arbre envoie une ar\^ete $e$ sur une ar\^ete distincte
$e'$ (en tenant compte aussi des orientations) et s'il existe un
chemin d'ar\^etes orient\'e g\'eod\'esique contenant $e$ et $e'$ avec
orientations compatibles, alors $\gamma$ n'a pas de point fixe, et
admet un axe de translation contenant $e$ et $e'$.

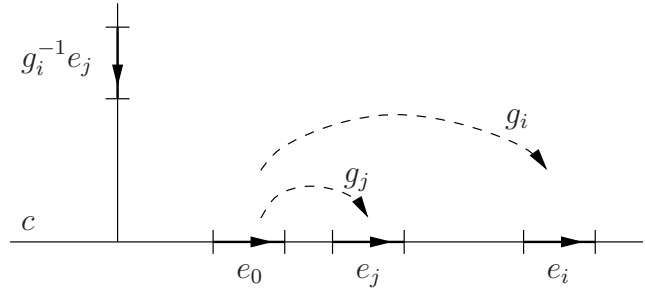
\begin{center}\input{fig_cylindre.pstex_t}\end{center}

Posons $g=g_j^{-1}g_i$. Comme l'ar\^ete $e_j$ est contenue dans le
sous-chemin d'ar\^etes orient\'e de $c$ entre $e_0$ et $e_i$, et comme
$g_i$ envoie $e_0$ sur $e_i$ (avec les orientations), on en d\'eduit
que les ar\^etes $g_i^{-1}e_j$ et $e_0$ sont contenues avec
orientations compatibles dans un chemin d'ar\^etes orient\'e
g\'eod\'esique.  Donc $g$, qui envoie $g_i^{-1}e_j$ sur $e_0$ (avec
orientations), n'a pas de point fixe. De plus, on a
$h'=g_ih_ig_i^{-1}= g_jh_jg_j^{-1}$, donc les \'el\'ements $h$ et $g$
commutent.  \cqfd

\medskip 
\rem 
Un automorphisme sans point fixe d'un arbre est d'ordre
infini. Rappelons qu'un groupe est {\it localement fini} si tout
sous-groupe de type fini est fini.  Si $\Ga\backslash T$ et les
stabilisateurs d'ar\^e\-te sont finis, alors, par la proposition
\ref{prop:ratrapage} pr\'ec\'edente, l'action de $\Gamma$ sur
$T$ est acylindrique si l'un des deux cas suivants est v\'erifi\'e~:
\begin{itemize}
\item le centralisateur de tout \'el\'ement d'ordre infini est
  cyclique,
\item 
le centralisateur de tout \'el\'ement non trivial ayant un point
  fixe est localement fini.
\end{itemize}

\bprop \label{prop:quotientsuspensionacylindBernoulli} %
Soit $T$ un arbre localement fini, et $\Ga$ un sous-groupe
g\'eom\'etri\-que\-ment fini de ${\rm Aut}(T)$. Alors, il existe un
arbre $T'$ et une action de $\Ga$ sur $T'$ qui v\'erifient~:
  \begin{itemize}
\item $\Ga\backslash T'$ est fini,
\item les stabilisateurs d'ar\^ete de $T'$ sont finis,
\item les stabilisateurs de sommet de $T'$ sont localement finis,
\item il existe une application simpliciale \'equivariante de $T$ dans
  $T'$.
\end{itemize}
De plus, la transformation g\'eod\'esique de $\Ga$ pour $T$, en
restriction \`a un $G_\delta$-dense de mesure de Bowen-Margulis totale
de $\Ga\backslash \G T_{\Ga,{\rm min}}$, s'obtient par suspension
topologique de la transformation g\'eod\'esique de $\Ga$ pour $T'$.

En outre, si $\Ga$ poss\`ede la propri\'et\'e de Selberg, alors il existe
un sous-groupe d'indice fini $\Gamma'$ de $\Gamma$ dont
l'action sur $T'$ est  acylindrique.  
\eprop

\rem
Lorsque $\Gamma$ est un r\'eseau uniforme, ceci est \'evident. En effet,
  $T'=T$ convient et, par des arguments de \cite{Ser}, le groupe
$\Ga$ contient un sous-groupe libre d'indice fini, agissant librement
sur $T$, donc de mani\`ere $1$-acylindrique.

\medskip \dem %
Comme $T$ se r\'etracte de mani\`ere $\Ga$-\'equivariante sur
$T_{\Ga,{\rm min}}$, nous pouvons supposer que $T=T_{\Ga,{\rm min}}$.
Par d\'efinition, le graphe de groupes quotient $(X,G_*)=\Ga\bac T$
est r\'eunion d'un graphe fini de groupes finis et d'un nombre fini de
rayons cuspidaux (maximaux), not\'es $(R_1,G_*), \dots, (R_r,G_*)$.
Pour $i$ dans $\{1,\dots, r\}$, notons $R'_i$ le rayon $R_i$ priv\'e
de son sommet origine et de son ar\^ete ouverte initiale, et
$(v_{i,n})_{n\in\NN}$ les sommets cons\'ecutifs de $R_i$. Notons
$\pi:T\ra X$ la projection canonique.  Soit $\sim$ la relation
d'\'equivalence sur $T$ engendr\'ee par $x\sim y$ si et seulement s'il
existe $i$ dans $\{1,\dots, r\}$ et une composante connexe de
$\pi^{-1}(R'_i)$ contenant $x$ et $y$.  Notons $T'$ le graphe quotient
de $T$ par $\sim$. Comme les composantes connexes de $\pi^{-1}(R'_i)$
sont des sous-arbres deux \`a deux disjoints de $T$, le graphe $T'$
est un arbre.  Puisque la relation d'\'equivalence $\sim$ est
invariante par $\Ga$, si l'on munit $T'$ de l'action quotient de
$\Ga$, alors la projection canonique $f:T\ra T'$ est
$\Ga$-\'equivariante.

Pour tout sommet $v$ de $T$, si $v$ n'appartient pas \`a
$\pi^{-1}(\bigcup_{i=1}^{r} R'_i)$, alors le stabilisateur dans $\Ga$
du sommet $f(v)$ de $T'$ est fini (\'egal au stabilisateur de $v$ dans
$T$). Si, par contre, $v$ appartient \`a une composante connexe de
$\pi^{-1}(R'_i)$, alors le stabilisateur de $f(v)$ est conjugu\'e au
groupe localement fini qui est la limite inductive des groupes des
sommets de $R_i$.

Notons ${\cal G}^\sharp T$ le $G_\delta$-dense de ${\cal G}T$ form\'e
des g\'eod\'esiques dont aucun bout  n'est un point
parabolique born\'e pour $\Ga$. Posons $X_1=\overline{X- \pi^{-1}
  \left(\bigcup_{i=1}^{r} R' _{\!i}\right)}$.  Notons ${\cal G}_1T$ le
sous-espace $\Ga$-invariant de ${\cal G}^\sharp T$ form\'e des
g\'eod\'esiques $\ell$ telles que d'une part $\ell(0)$ appartienne \`a
$X_1$, et d'autre part si $\ell(0)$ est dans
$\pi^{-1}(\{v_{i,1}\;:\;i\in\{1,\dots,r\}\})$, alors $\ell(-1)$ ne
soit pas dans l'int\'erieur de $\pi^{-1}(\bigcup_{i=1}^{r} R'_{\!i})$.

Pour $\ell$ dans ${\cal G}_1T$, notons $\tau(\ell)$ l'entier non nul,
valant $1$ si $\ell(1)\in X_1$ et valant $1+\min\{n\in\NN-\{0\}\;:\;
\ell(n)\in X_1\}$ sinon. Il s'agit du temps de premier retour dans
${\cal G}_1T$ de l'orbite de $\ell$ par la transformation
g\'eod\'esique $\varphi$. Notons $\psi:{\cal G}_1T\ra{\cal G}_1T$
l'application de premier retour, d\'efinie par $\psi:\ell\mapsto
\varphi^{\tau(\ell)}(\ell)$.

Soit $\ell$ une g\'eod\'esique appartenant \`a ${\cal G}^\sharp T$.
Alors chaque composante connexe de $\ell(\RR)\cap
\pi^{-1}(\bigcup_{i=1}^{r} R'_i)$ est un segment compact.  La
restriction \`a $\ell(\RR)$ de $f$ \'ecrase en un point chacune de ces
composantes connexes. L'image par $f:T\ra T'$ de $\ell(\RR)$ est donc
l'image d'une g\'eod\'esique $\theta\ell$ de $T'$, d'origine
$\theta\ell(0)=f(\ell(0))$, de sorte que
$f:\ell(\RR)\ra\theta\ell(\RR)$ pr\'eserve l'orientation.  Notons
$\varphi':{\G}T'\ra{\G}T'$ la transformation g\'eod\'esique pour $T'$.

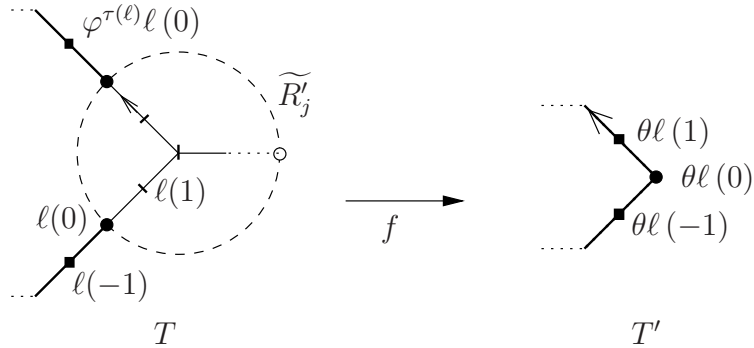
\medskip
\begin{center}
\input{fig_ecrashor.pstex_t}
\end{center}

L'application $\theta:{\G}_1T\ra {\G}T'$ d\'efinie par $\ell\mapsto
\theta\ell$ est un hom\'eomorphisme qui rend le diagramme suivant
commutatif:
$$\begin{array}[b]{ccc}
{\G}_1T & \stackrel{\theta}{\longrightarrow} &{\G}T' \\
\psi\downarrow\;\;&   & \;\;\downarrow \varphi'\\
{\G}_1T & \stackrel{\theta}{\longrightarrow} &{\G}T' 
\end{array}\;.$$
De plus, ce diagramme commute avec l'action de $\Ga$.  Donc la
restriction \`a $\Ga\backslash\G^\sharp T$ de la transformation
g\'eod\'esique sur $\Ga\backslash\G T$ est topologiquement conjugu\'ee \`a une
suspension de la transformation g\'eod\'esique pour $\Ga\backslash\G T'$.
La mesure de Patterson-Sullivan de tout point parabolique born\'e est
nulle (voir \cite{DOP} dans le cadre des vari\'et\'es, et \cite{Rob} plus
g\'en\'eralement). Comme il n'y a qu'un ensemble d\'enombrable de points
paraboliques born\'es pour $\Ga$, le
sous-espace $\G^\sharp T$ est donc de mesure pleine.

Si $\Ga$ poss\`ede la propri\'et\'e de Selberg, alors il existe un
sous-groupe $\Ga'$ d'indice fini de $\Ga$, tel que le graphe de
groupes quotient $\Ga'\bac T$ soit r\'eunion d'un graphe fini de groupes
triviaux, et d'un nombre fini de rayons cuspidaux. Tout \'el\'ement non
trivial de $\Ga'$, fixant un point de $T'$, fixe donc un unique point de
$T'$, dont le stabilisateur est localement fini. Par la seconde
assertion de la remarque qui pr\'ec\`ede la proposition
\ref{prop:quotientsuspensionacylindBernoulli}, l'action de $\Ga'$ sur
$T'$ est donc acylindrique.  
\cqfd

\brema 
Soit $T$ un arbre localement fini, et $\Ga$ un
sous-groupe de ${\rm Aut}(T)$ tel que $\Ga\bac T$ soit r\'eunion de $n$
rayons cuspidaux de groupes, recoll\'es en leurs origines. S'il existe
$m$ dans $\NN-\{0\}$ tel que l'action de $\Gamma$ sur les $m$-uplets
de points deux \`a deux distincts de $\partial T$ soit libre, alors
l'action de $\Gamma$ sur $T'$ (construit dans la preuve ci-dessus) est
$(2m-1)$-acylindrique.  
\erema

Remarquons que le graphe de groupes $\Ga\bac T'$ est alors une \'etoile
de groupes de la forme
\\ \centerline{\input{fig_etoilegroupes.pstex_t}}\\
avec, lorsque $\Ga$ est discret, $G,G'_1,\dots, G'_n$ des groupes
finis, et $G_1,\dots, G_n$ des groupes localement finis.

\medskip 
\dem %
La projection dans $\Ga\backslash T'$ d'un chemin d'ar\^etes
g\'eod\'esique $c$ dans $T'$, de longueur $2m-1$, comporte $m-1$
allers-retours plus une ar\^ete, donc est la projection d'un chemin
d'ar\^etes g\'eod\'esique de $T$ p\'en\'etrant dans $m-1 +1$
horoboules distinctes de la famille $\Gamma$-invariante maximale
d'horoboules d'int\'erieurs disjoints. Tout \'el\'ement de $\Ga$
fixant $c$ fixe donc dans $\partial T$ les $m$ points \`a l'infini de
ces horoboules.  
\cqfd

\medskip %
Les hypoth\`eses de cette remarque sont v\'erifi\'ees, avec $n=1$ et
$m=3$, lorsque l'on consid\`ere le groupe $\Ga={\rm PGL}(2,\FF_q[X])$
agissant sur l'arbre de Bruhat-Tits $T$ de $({\rm
  PGL}_2,\FF_q((X^{-1})))$.  En particulier, comme annonc\'e dans
l'introduction, l'action de $\Ga$ sur $T$ n'est pas acylindrique, mais
celle de $\Ga$ sur l'arbre $T'$ associ\'e \`a $T$ par la preuve de la
proposition \ref{prop:quotientsuspensionacylindBernoulli} est
$5$-acylindrique.

\subsection{Flot g\'eod\'esique d'ordre $k$ sur un 
graphe de groupes}
\label{subseq:ordreka}

Soit $(X,G_\ast)$ un graphe de groupes et $k>0$ un entier. Un
{\it $k$-chemin} de $X$ est une suite $(e_0,e_1,\dots,e_k)$ d'ar\^etes
de $EX$ telles que $t(e_{i-1})=o(e_i)$ pour $i=1,\dots,k$. Pour tout
$k$-chemin, le groupe produit $G_{e_0}\times G_{e_1}\times \dots\times
G_{e_k}$ agit sur l'ensemble $G_{o(e_1)}\times G_{o(e_2)}\times
\dots\times G_{o(e_k)}$ par
$$((\alpha_0,\alpha_1,\dots,\alpha_k),(g_1,\dots,g_k)\mapsto
(\rho_{e_0}(\alpha_0)\,g_1\, \rho_{\overline{e_1}}(\alpha_1^{-1}),
\,\dots\;, \rho_{e_{k-1}}(\alpha_{k-1})\,g_k\,
\rho_{\overline{e_k}}(\alpha_k^{-1}))\;.$$ Notons
$G(e_0,e_1,\dots,e_k)$ l'ensemble quotient.  Nous noterons
$[g_1,\dots,g_k]$ la classe d'un \'el\'ement $(g_1,\dots,g_k)$ de
$G_{o(e_1)}\times G_{o(e_2)}\times \dots\times G_{o(e_k)}$.

Consid\'erons les trois applications $$(g_1,...,g_k)\mapsto
(g_2,...,g_k), \;\;(g_1,...,g_k)\mapsto(g_1,...,g_{k-1}) {\rm
  ~~et~~}(g_1,...,g_k)\mapsto g_k$$ (les deux premi\`eres n'\'etant
d\'efinies que si $k\geq 2$). Elles induisent des surjections ${\rm
  pr}_1:G(e_0,e_1,\dots,e_k)\ra G(e_1,\dots,e_k)$, ${\rm
  pr}_2:G(e_0,e_1,\dots,e_k)\ra G(e_0,\dots,e_{k-1})$ et \\ $ {\rm
  pr}_3:G(e_0,e_1,\dots,e_k)\ra G(e_{k-1},e_k).$

Remarquons d'ailleurs
que $G(e_{k-1},e_k)=\rho_{e_{k-1}}(G_{e_{k-1}})\backslash
G_{o(e_k)}/\rho_{\overline{e_k}}(G_{e_{k}})$.

Notons $\Omega'$ la r\'eunion disjointe des ensembles
$G(e_0,e_1,\dots,e_k)$ pour tous les $k$-chemins $(e_0,e_1,\dots,e_k)$
de $X$.  Munissons l'ensemble $\Omega'\times EX$ de la topologie discr\`ete
et l'ensemble ${\displaystyle\Sigma'=(\Omega'\times EX)^\ZZ}$ de la
topologie produit.  Remarquons que si $(X,G_\ast)$ est un graphe fini
de groupes finis, alors $\Omega'\times EX$ est fini. Notons encore
$\sigma$ le d\'ecalage sur $\Sigma'$.

Notons $\Sigma'_0$ l'ensemble des \'el\'ements
$([g_{i,1},g_{i,2},\dots,g_{i,k}],e_i)_{i\in \ZZ}$ de $\Sigma'$ tels
que pour tout entier $i$,
$$
\left\{\begin{array}{l} (1')\;\; o(e_{i})=t(e_{i-1})\\
       (2')\;\; [g_{i,1},g_{i,2},\dots,g_{i,k}]\in G(e_{i-k},
       \dots,e_{i-1},e_{i})
       \\
       (3') \;\;{\rm pr}_{1}([g_{i,1},g_{i,2},\dots,g_{i,k}])=
       {\rm pr}_{2}([g_{i+1,1},g_{i+1,2},\dots,g_{i+1,k}])
       \\
      (4') \;\;{\rm si}\; e_{i}=\overline{e_{i-1}}, \;{\rm alors}\;
       {\rm pr}_{3}([g_{i,1},g_{i,2},\dots,g_{i,k}])\notin
       \rho_{\overline{e_{i}}}(G_{e_{i}})\;.
\end{array}\right.
$$

D\'efinissons alors le {\it flot g\'eod\'esique d'ordre $k$} sur
$(X,G_\ast)$ comme le sous-d\'ecalage $(\Sigma'_0,\sigma)$.  Lorsque
$k=1$, on retrouve la d\'efinition de la partie \ref{subseq:ordreun}.

Il est imm\'ediat que $\sigma$ pr\'eserve $\Sigma'_0$, et que les
conditions d\'efinissant $\Sigma'_0$ sont ``locales'' (au temps $i$,
elles ne d\'ependent que du temps $i$ et des $k$ termes
pr\'ec\'edents). En particulier, $(\Sigma'_0,\sigma)$ est un
sous-d\'ecalage de type fini lorsque $(X,G_\ast)$ est un graphe fini
de groupes finis.

Montrons maintenant le r\'esultat suivant~:

\btheo \label{prop:flotgeod_ordqqc} Soit $T$ un arbre et $\Ga$ un
sous-groupe de ${\rm Aut}(T)$. Supposons que l'action de $\Gamma$ sur
$T$ soit $k$-acylindrique. Alors la transformation g\'eod\'esique
quotient est topologiquement conjugu\'e \`a un sous-d\'ecalage.

Plus pr\'ecis\'ement, si $(\Sigma'_0,\sigma)$ est le flot
g\'eod\'esique d'ordre $k$ du graphe de groupes quotient $\Ga\bac T$,
alors il existe un hom\'eomorphisme $\theta':\Ga\backslash \G T \ra
\Sigma'_0$ qui rend le diagramme suivant commutatif
$$\begin{array}[b]{ccc} \Ga\backslash \G T &
\stackrel{\varphi}{\longrightarrow} & \Ga\backslash \G T \\
\theta'\downarrow\;\; & & \;\;\downarrow \theta' \\ \Sigma'_0&
\stackrel{\sigma}{\longrightarrow}& \Sigma'_0
\end{array}\;.$$
\etheo

Lorsque $T$ est localement fini, et $\Ga$ un sous-groupe discret
cocompact, ce r\'esultat est d\'eja connu~: voir \cite{Moz} pour un
cadre alg\'ebrique (rappel\'e en introduction), qui s'\'etend en rang
sup\'erieur; voir \cite{CP} pour le cadre plus g\'en\'eral des groupes
hyperboliques, par des m\'ethodes utilisant leur dynamique sur leur
espace \`a l'infini.

\medskip
\dem %
Consid\'erons l'application $\wt{\theta'}:\G T\ra \Sigma'_0$, d\'efinie par
$$\wt{\theta'}(f)=([g_{i-k+1},\dots,g_{i-1},g_{i}],e_i)_{i\in \ZZ}\;,$$
en reprenant dans la d\'emonstration de la proposition
\ref{prop:flotgeod_ordun} les objets $f_i, e_i, h_i, g_i$ associ\'es
\`a $f$. En effet, la classe de $(g_{i-k+1},\dots,g_{i-1},g_{i})$ dans
$G(e_{i-k}, \dots,e_{i-1},e_{i})$ ne d\'epend pas du choix des $h_i$.
On v\'erifie imm\'ediatement que les conditions ($1'$)-($4'$)
ci-dessus sont satisfaites pour $\wt{\theta'}(f)$.

De la m\^eme mani\`ere que pour $\wt\theta$, l'application
$\wt{\theta'}:\G T\ra \Sigma'_0$ est continue et $\Ga$-invariante.
L'application
$$\rho:([g_{i,1},g_{i,2},\dots,g_{i,k}],e_i)_{i\in \ZZ}\mapsto ({\rm
pr}_{3}([g_{i,1},g_{i,2},\dots,g_{i,k}]),e_i)_{i\in \ZZ}$$ est une
application continue de $\Sigma'_0$ sur $\Sigma_0$. Il est imm\'ediat
que $\rho\circ\wt{\theta'}=\wt{\theta}$.

\medskip Montrons que $\wt{\theta'}$ est surjective. Soit
$([g_{i,1},g_{i,2},\dots,g_{i,k}],e_i)_{i\in \ZZ}$ un \'el\'ement de
$\Sigma'_0$.  Montrons que pour tout entier $i$, il existe un
\'el\'ement $g_i$ de $G_{o(e_i)}$ tel que
$$[g_{i-k+1},\dots,g_{i-1},g_i]=[g_{i,1},g_{i,2},\dots,g_{i,k}]$$
dans $G(e_{i-k+1},\dots,e_{i-1},e_i)$. Posons $g_{-k+1}=g_{0,1}$,
\dots, $g_{-1}=g_{0,k-1}$ et $g_0=g_{0,k}$. Cons\-trui\-sons les $g_n$
pour $n\geq 1$ par r\'ecurrence sur $n$ (et de mani\`ere similaire les
$g_{-n}$ pour $n\geq k$).

Soit $n\geq 0$, et supposons construits $g_{-k+1},\dots, g_n$. En
particulier, l'\'egalit\'e
$$[g_{n-k+1},\dots,g_{n-1},g_n]=[g_{n,1},g_{n,2},\dots,g_{n,k}]$$
est v\'erifi\'ee dans $G(e_{n-k+1},\dots,e_{n-1},e_n)$. Par la
condition $(3')$,
$${\rm pr}_{1}[g_{n-k+1},\dots,g_{n-1},g_n]=
       {\rm pr}_{2}[g_{n+1,1},g_{n+1,2},\dots,g_{n+1,k}]\;.$$
Donc il existe $h_1$ dans $G_{e_{n-k+1}}$, \dots, $h_k$ dans
$G_{e_{n}}$ tels que
$$g_{n-k+2}=\rho_{e_{n-k+1}}(h_1)g_{n+1,1}\rho_{\ov
e_{n-k+2}}(h_2)^{-1}, \dots,
g_{n}=\rho_{e_{n-1}}(h_{k-1})g_{n+1,k-1}\rho_{\ov
e_{n}}(h_k)^{-1}\;.$$

Posons maintenant $g_{n+1}=\rho_{e_{n}}(h_k)
g_{n+1,k}$. C'est un \'el\'ement de $G_{e_{n+1}}$. De plus,
$$[g_{n-k+2},\dots,g_{n},g_{n+1}]=
[g_{n+1,1},g_{n+1,2},\dots,g_{n+1,k}]\;,$$
ce qui termine la r\'ecurrence.

Par la surjectivit\'e de $\wt{\theta}$, il existe $f\in \G T$ tel que
$\wt{\theta}(f)=(g_i,e_i)_{i\in\ZZ}$. Alors, par cons\-truc\-tion de
$\wt{\theta'}$, il est imm\'ediat que
$\wt{\theta'}(f)=([g_{i,1},g_{i,2},\dots,g_{i,k}],e_i)_{i\in \ZZ}$, ce
qui montre la surjectivit\'e de $\wt{\theta'}$.

\medskip %
D'apr\`es ce qui pr\'ec\`ede, l'application $\wt{\theta'}$ induit une
application continue surjective $\theta':\Ga\backslash \G T\ra
\Sigma'_0$.

\medskip %
Montrons que $\theta'$ est injective. Soient $f,f'$ dans $\G T$ tels
que $\wt{\theta'}(f)= \wt{\theta'}(f')$. Montrons qu'il existe $\ga$
dans $\Ga$ tel que $f'=\ga f$. Notons $f'_i, e'_i, h'_i, g'_i$ les
objets correspondants pour $f'$ \`a ceux introduits dans la
d\'emonstration de la proposition \ref{prop:flotgeod_ordun} pour $f$.
Comme $\wt{\theta'}(f)= \wt{\theta'}(f')$, on a $e'_i=e_i$.  Quitte
\`a remplacer $f$ et $f'$ par leurs images par des \'el\'ements de
$\Ga$, nous pouvons supposer que $f_0=f'_0=\wt{e_0}$, et dans la
construction, nous pouvons alors prendre $h_0=h'_0={\rm id}$.

Comme $\wt{\theta'}(f)= \wt{\theta'}(f')$, il existe pour tout  $i$,
une suite $(\alpha_{i,-k},\alpha_{i,-k+1},\dots,\alpha_{i,0})$ dans
$G_{e_{i-k}}\times \dots\times G_{e_{i-1}}\times G_{e_i}$ telle que,
pour tout $j=1,\dots,k$,
$$g'_{i-k+j}= \rho_{e_{i-k+j-1}}(\alpha_{i,j-1-k})\;g_{i-k+j}\;
\rho_{\,\overline{e_{i-k+j}}}\,
(\alpha_{i,j-k}\mbox{}^{-1})\;.\;\;\;(**)$$

La d\'emonstration du lemme suivant est technique. Ceci s'explique par
la possibilit\'e de propagation de fixateurs d'ar\^etes dans des
suites d'ar\^etes. C'est pour lui que nous avons besoin de
l'hypoth\`ese que l'action est $k$-acylindrique.

\blemm
Pour tout $i$ dans $\ZZ$, on a $\alpha_{i,-1}=\alpha_{i-1,0}$.
\elemm

\dem %
Le terme de gauche de l'\'egalit\'e $(**)$ \'etant inchang\'e en
rempla\c{c}ant simultan\'ement $i$ par $i-1$ et $j$ par $j+1$, il en
est de m\^eme pour le terme de droite. Donc, pour tout
$j=1,\dots,k-1$,
$$\rho_{e_{i-k+j-1}}(\alpha_{i,j-1-k})g_{i-k+j}
\rho_{\,\overline{e_{i-k+j}}}\,(\alpha_{i,j-k}^{-1})=
\rho_{e_{i-k+j-1}}(\alpha_{i-1,j-k})g_{i-k+j}
\rho_{\,\overline{e_{i-k+j}}}\,(\alpha_{i-1,j+1-k}^{-1}).$$
Ceci \'equivaut \`a
$$\rho_{e_{i-k+j-1}}(\alpha_{i-1,j-k}^{-1}\alpha_{i,j-1-k})=
g_{i-k+j}\;\rho_{\,\overline{e_{i-k+j}}}\, 
(\alpha_{i-1,j+1-k}^{-1}\alpha_{i,j-k})
\;g_{i-k+j}^{-1}\;.\;\;\;\;(***)$$
Comme $\alpha_{m,\ell}$ appartient \`a $G_{e_{m+\ell}}$, l'\'el\'ement
$\alpha_{i-1,j-k}^{-1}\alpha_{i,j-1-k}$ de $\Ga$ fixe l'ar\^ete
$\wt{e_{i-k+j-1}}$.

Notons $\varepsilon_{i-k}=g_{e_{i-k}}^{-1}\wt{e_{i-k}}$, et pour
$j=1,\dots,k-1$,
\begin{eqnarray*}
\varepsilon_{i-k+j}
&=&
(g_{i-k+1}\;
g_{\,\overline{e_{i-k+1}}}^{-1}\;g_{e_{i-k+1}})(g_{i-k+2}\;
g_{\,\overline{e_{i-k+2}}}^{-1}\;g_{e_{i-k+2}})\\
&&\dots
(g_{i-k+j-1} \;g_{\,\overline{e_{i-k+j-1}}}^{-1}\;g_{e_{i-k+j-1}})\,
g_{i-k+j}\;g_{\,\overline{e_{i-k+j}}}^{-1}\;\wt{e_{i-k}}.
\end{eqnarray*}
Comme $g_{\ell}$ fixe $\wt{t(e_{\ell-1})}=\wt{o(e_{\ell})}$, la suite
$(\varepsilon_{i-k},\varepsilon_{i-k+1},\dots,\varepsilon_{i-1})$ est
un $(k-1)$-chemin de $T$ (voir la figure 
ci-dessous).

\medskip
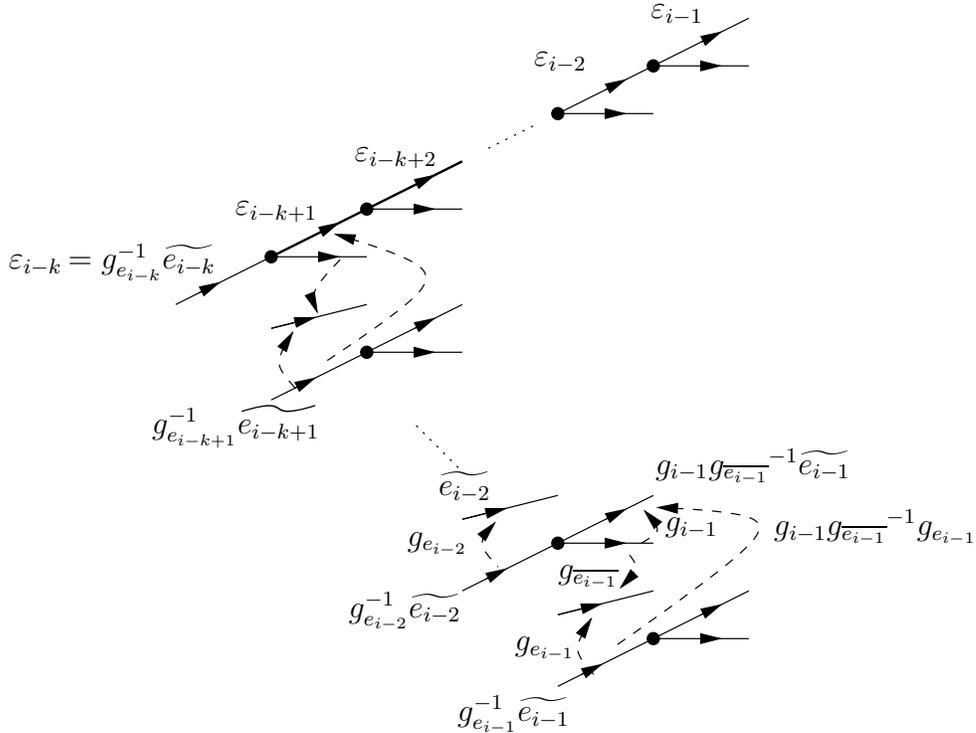
\input{fig_ledrappi2.pstex_t}

\medskip 
Remarquons que, pour $\ell=0,\dots,k-1$, l'\'el\'ement
$\rho_{e_{i-k+\ell}}(\alpha_{i-1,\ell+1-k} ^{-1}\alpha_{i,\ell-k})$
fixe l'ar\^ete $g_{e_{i-k+\ell}}^{-1}\wt{e_{i-k+\ell}}$.

L'\'el\'ement $\rho_{e_{i-k}}(\alpha_{i-1,1-k}^{-1} \alpha_{i,-k})$,
qui fixe $g_{e_{i-k}}^{-1}\wt{e_{i-k}}= \varepsilon_{i-k}$, est
\'egal, d'apr\`es l'\'egalit\'e \mbox{$(***)$} pour $j=1$, \`a un
\'el\'ement fixant $$g_{i-k+1} g_{\,\overline{e_{i-k+1}}}^{-1}
\wt{e_{i-k+1}}= (g_{i-k+1} g_{\,\overline{e_{i-k+1}}}^{-1}
g_{e_{i-k+1}})g_{e_{i-k+1}}^{-1}\wt{e_{i-k+1}}=
\varepsilon_{i-k+1}\;.$$ Donc
$\rho_{e_{i-k}}(\alpha_{i-1,1-k}^{-1}\alpha_{i,-k})$ fixe
$(\varepsilon_{i-k},\varepsilon_{i-k+1})$. D'apr\`es la condition
$(4')$, cette suite de deux ar\^etes cons\'ecutives n'est pas un
aller-retour.  Une r\'ecurrence imm\'ediate montre que l'\'el\'ement
$\rho_{e_{i-k}}(\alpha_{i-1,1-k} ^{-1} \alpha_{i,-k})$ fixe le
$(k-1)$-chemin $(\varepsilon_{i-k}, \varepsilon_{i-k+1},
\dots,\varepsilon_{i-1})$, qui est localement injectif.  Il vaut donc
l'identit\'e, car l'action de $\Ga$ est $k$-acylindrique.  Puisque les
morphismes $\rho_{e_n}$ sont injectifs, on a donc
$\alpha_{i-1,1-k}=\alpha_{i,-k}$. Une r\'ecurrence imm\'ediate
utilisant les \'egalit\'es $(***)$ pour $j=1,\dots,k-1$ montre alors
que
$$\alpha_{i-1,\ell-k}=\alpha_{i,\ell-1-k}$$
pour $\ell=1,\dots,k$, ce
qui montre le lemme pour $\ell=k$.
\cqfd

\medskip %
Reprenons la d\'emonstration de l'injectivit\'e de $\theta'$.  Par
l'\'egalit\'e $(**)$ pour $j=k$, et par le lemme pr\'ec\'edent, on a
$$g'_i=\rho_{e_{i-1}}(\alpha_{i,-1})\;g_{i}\;
\rho_{\,\overline{e_{i}}}\,(\alpha_{i,0}^{-1})=
\rho_{e_{i-1}}(\alpha_{i-1,0})\;g_{i}\;
\rho_{\,\overline{e_{i}}}\,(\alpha_{i,0}^{-1})\;.$$
La formule $(*)$ de la d\'emonstration de la proposition
\ref{prop:flotgeod_ordun} donne donc
$$g_{e_{i-1}}^{-1} h'_{i-1} (h'_i)^{-1}g_{\overline{e_i}} =
\rho_{e_{i-1}}(\alpha_{i-1,0})\;g_{e_{i-1}}
^{-1}h_{i-1}h_i^{-1}g_{\overline{e_i}}\;
\rho_{\,\overline{e_{i}}}\,(\alpha_{i,0}^{-1}),$$
ce qui \'equivaut \`a
$$h'_i(h'_{i-1})^{-1}=\alpha_{i,0}h_ih_{i-1}
^{-1}\alpha_{i-1,0}^{-1}\;.$$
Posons $f''=\alpha_{0,0}^{-1}f'$.  Alors, par naturalit\'e, pour
construire $\wt\theta(f'')$, on peut prendre $h''_i=h'_i\alpha_{0,0}$.
Donc, pour tout  $i$ dans $\ZZ$,
$$h''_i(h''_{i-1})^{-1}=\alpha_{i,0}h_ih_{i-1}
^{-1}\alpha_{i-1,0}^{-1}.$$
Montrons par r\'ecurrence sur $n\geq 0$ que $h''_n=\alpha_{n,0}h_n$.
Comme $h_0=h'_0={\rm id}$, on a $h''_0=h'_0\alpha_{0,0}=
\alpha_{0,0}h_0$. Supposons la formule vraie au rang $n-1$, alors
$$h''_n=\alpha_{n,0}h_nh_{n-1}
^{-1}\alpha_{n-1,0}^{-1}h''_{n-1}
=\alpha_{n,0}h_n\;.$$
Un raisonnement analogue pour $n\leq 0$ montre que pour tout $i$ dans
$\ZZ$, on a encore $h''_i=\alpha_{i,0}h_i$.  Donc, puisque
$\alpha_{i,0}$ fixe l'ar\^ete $\wt{e_i}$, on a
$$f''_i=(h''_i)^{-1}\wt{e_i}=h_i^{-1}\alpha_{i,0}^{-1}\wt{e_i}=
h_i^{-1}\wt{e_i}=f_i\;.$$
Donc $f''=f$, ce qui montre bien que $f$ et $f'$ sont dans la m\^eme
orbite par $\Ga$. Donc $\theta'$ est injective.

Cette d\'emonstration montre aussi que si les $i$-\`emes termes des
suites $\wt{\theta'}(f)$ et $\wt{\theta'}(f')$ sont \'egaux pour $i$ dans
$[-N-2,N+2]\;\cap\;\NN$, alors $f$ et $f'$ co\"{\i}ncident sur
$[-N,+N]$. Donc $(\theta')^{-1}$ est continue. Ceci termine la
d\'emonstration du th\'eor\`eme \ref{prop:flotgeod_ordqqc}.
\cqfd


\begin{thebibliography}{BGS}
     {\small

\bibitem[Alp]{Alp}
 R.~Alperin,  {\it An elementary account of Selberg's lemma},
 L'Ens. Math. {\bf 33}  (1987)  269-373.

\bibitem[ASS]{smorodinsky}%
R.L.~Adler, P.~Shields, M.~Smorodinsky, {\it Irreducible Markov shifts},
Ann. Math. Statistics,~{\bf 43} (1972) 1027-1029.

\bibitem[BL]{BL} 
H.~Bass, A.~Lubotzky, {\it Tree lattices}, Prog.
Math. {\bf 176}, Birkh\"auser, 2001.

\bibitem[BT]{BT} 
F.~Bruhat, J.~Tits, {\it Groupes r\'eductifs sur un
corps local (donn\'ees radicielles valu\'ees)},  Pub.~Math.~I.H.E.S.    
{\bf 41} (1972), 5-252.

\bibitem[Bou]{Bou}
M.~Bourdon, {\it Structure conforme au bord et flot g\'eod\'esique
d'un CAT($-1$) espace},  L'Ens. Math. {\bf 41} (1995) 63-102.

\bibitem[BH]{BH}
M.R.~Bridson, A.~Haefliger, 
{\it Metric spaces with non-positive curvature}, 
Grund. math.~Wiss. {\bf 319}, Springer Verlag (1998).

\bibitem[BP]{BP} 
A.~Broise-Alamichel, F.~Paulin, {\it Dynamique sur le rayon
modulaire et fractions continues en caract\'eristique $p$},
Pr\'epublication, Univ. Orsay 2002.

\bibitem[BM]{BM} 
M.~Burger, S.~Mozes, {\it CAT($-1$) spaces, divergence 
groups and their commensurators}, 
J. Amer. Math. Soc {\bf 9} (1996) 57-94.

\bibitem[Coo]{Coo} 
M.~Coornaert, {\it Mesures de Patterson-Sullivan
sur le bord d'un espace hyperbolique au sens de Gromov},
Pacific J. Math. {\bf 159} (1993) 241--270.

\bibitem[CP]{CP}
M.~Coornaert, A.~Papadopoulos, {\it Symbolic dynamics and hyperbolic
groups}, Lect. Notes Math. {\bf 1539}, Springer Verlag, 1993.

\bibitem[DOP]{DOP}
F.~Dal'Bo, J.-P.~Otal, M.~Peign\'e, 
{\it S\'eries de Poincar\'e des groupes g\'eom\'etriquement finis},
Israel J. Math. {\bf  118} (2000) 109--124.

\bibitem[FO]{FO} 
N.~Friedman, D.~Ornstein, {\it On isomorphism of weak Bernoulli 
transformations}, Adv. Math. {\bf 5} (1971) 365-394.

\bibitem[GL]{GL}
D.~Gaboriau, G.~Levitt,   {\it  The rank of actions on $\RR$-trees},  
Ann. Scien. Ec. Norm. Sup. (4) {\bf 28} (1995) 549-570.

\bibitem[HK]{HK}
B.~Hasselblatt, A.~Katok,
eds., {\it Handbook of Dynamical Systems}, Elsevier, 2002.

\bibitem[HP1]{HP} 
S.~Hersonsky, F.~Paulin, {\it Counting orbit points
in covering of negatively curved manifolds and Hausdorff dimension of
cusp excursions}, Erg. Theo. Dyn. Sys., {\bf 24}, (2004), 803-824.

\bibitem[HP2]{HP2} S.~Hersonsky, F.~Paulin, {\it A logarithm law for
    tree automorphism groups}, en pr\'eparation.

\bibitem[Kai]{Kai} 
V.A.~Kaimanovich, {\it Bowen-Margulis and Patterson
    measures on negatively curved compact manifolds}. in ``Dynamical
  systems and related topics'' (Nagoya, 1990), 223--232, Adv. Ser.
  Dyn. Syst., {\bf 9}, World Sci. Pub., 1991.

\bibitem[Kit]{Kit} 
B.~Kitchens, {\it Symbolic dynamics: one-sided, two sided and
  countable state Markov shifts}, Universitext, Springer Verlag, 1998.

\bibitem[LP]{LP}
F.~Ledrappier, M.~Pollicot, {\it Distribution results for lattices in
SL$(2,\QQ_p)$}, Bul. Braz. Math. Soc. {\bf 36} (2005) 143-176.

\bibitem[LW]{LW} 
E.~Lindenstrauss, B.~Weiss, {\it On sets invariants
    under the action of the diagonal group},  Erg.  Theo. Dyn. Sys.
    {\bf 21} (2001) 1481-1500.

\bibitem[Lub]{Lub}
A.~Lubotzky, {\it  Lattices in rank one Lie groups over local fields},  
GAFA {\bf 1} (1991) 405-431.

\bibitem[Mar1]{Mar1}
G.~Margulis,  {\it Discrete subgroups of semi-simple groupes}, 
Ergeb. Math. Grenz. {\bf 17}, Springer Verlag, 1991.

\bibitem[Mar2]{Mar2} 
G.~Margulis, {\it Problems and conjectures in rigidity theory}, 
in ``Mathematics: frontiers and perspectives 2000'',
161-174, Amer. Math. Soc. 2000.

\bibitem[Moz]{Moz}
S.~Mozes, {\it Actions of Cartan subgroups},
Israel J. Math. {\bf 90} (1995) 253-294.

\bibitem[Orn]{Orn} 
D.~Ornstein, {\it Factors of  Bernoulli shifts are Bernoulli shifts}, 
Adv. Math. {\bf 5} (1971) 349-364.

\bibitem[Pau]{Pau} 
   F.~Paulin, {\it Groupes g\'eom\'etriquement finis
    d'automorphismes d'arbres et approximation diophantienne dans les
    arbres}, Manuscripta Math. {\bf 113} (2004) 1-23.

\bibitem[Rob]{Rob} 
T.~Roblin, {\it Ergodicit\'e et  \'equidistribution en courbure 
n\'egative},  M\'emoires Soc. Math. France, {\bf 95}, (2003). 

\bibitem[Sel]{Sel}
Z.~Sela, {\it Acylindrical accessibility for groups},
Inv. Math. {\bf 129} (1997) 528-565.

\bibitem[Ser]{Ser}
J.-P.~Serre, {\it Arbres, amalgames, SL$_2$},
Ast\'erisque  {\bf 46}, Soc. Math. France (1983).

\bibitem[Tom]{Tom}
G.~Tomanov, {\it Actions of maximal tori on homogeneous spaces}, in
``Rigidity in dynamics and geometry'' (Cambridge, 2000), M. Burger,
A. Iozzi eds, Springer Verlag (2002), 407-424.

\bibitem[Tho]{Tho}
J.-P..~Thouvenot, {\it Entropy, isomorphism and equivalence in
ergodic theory}, in Hand. Dyn. Sys. Vol.~1A, B.~Hasselblatt, A.~Katok
eds., Elsevier, 2002, 205-238.

\bibitem[Zim]{Zim} 
R.J.~Zimmer,  {\it Ergodic theory and semisimple groups}, 
Birkhauser, 1984.

}
\end{thebibliography}
\end{document}